\documentclass[11pt]{article}%

\RequirePackage{amsthm, amsmath, amsfonts, amssymb}%

\usepackage{amsmath,mathtools}%
\usepackage{amssymb}%
\usepackage[dvipsnames,svgnames]{xcolor}%
\usepackage[titletoc,title]{appendix}%
\usepackage[utf8]{inputenc}%
\usepackage{dsfont}%
\usepackage{caption}
\usepackage{float}
\usepackage{mathrsfs}

\usepackage{graphicx, color}%

\usepackage{tikz}%
\usepackage{tkz-base}%
\usepackage{tkz-fct}%
\usepackage{amsmath}%
\usepackage{pgfplots}%
\usetkzobj{all}%
\usetikzlibrary{calc}%
\usetikzlibrary{arrows,positioning}

\usepackage[linesnumbered,ruled,vlined]{algorithm2e}

\usepackage[a4paper, left=3cm, right=3cm, top=3cm, bottom=3cm]{geometry}

\definecolor{darkblue}{rgb}{0.0,0.0,0.7}

\RequirePackage[%
colorlinks = true,%
linkcolor = darkblue,%
citecolor = darkblue,%
urlcolor = darkblue, %
]{hyperref}%

\hypersetup{%
  pdfauthor = {St\'ephane Ga\"iffas, Agathe Guilloux},%
  pdftitle = {},%
  pdfcreator = {pdflatex},%
  pdfproducer = {pdflatex}}

\newcommand{\Sum} {\displaystyle \sum\limits}              

\newcommand{\R}{{\mathbb{R}}}
\newcommand{\E}{\mathds{E}}

\newcommand{\cE}{\mathcal E}
\newcommand{\cT}{\mathcal T}

\newtheorem{thm}{Theorem}
\newtheorem{cor}{Corollary}
\newtheorem{lem}{Lemma}
\newtheorem{prop}{Proposition}

\newtheorem{defn}{Definition}

\usepackage{datenumber}

\newcommand{\inr}[1]{\langle #1 \rangle}
\newcommand{\ind}[1]{{\mathds{1}}_{{#1}}}%
\newcommand{\norm}[1]{\|#1\|}

\renewcommand{\P}{\mathds{P}}
\DeclareMathOperator{\argmin}{argmin}
\DeclareMathOperator{\TV}{TV}
\DeclareMathOperator{\prox}{prox}
\DeclareMathOperator{\sgn}{sign}
\DeclareMathOperator{\minimize}{minimize}
\DeclareMathOperator{\mise}{MISE}

\newcommand{\bN}{\mathbf N}%
\newcommand{\bT}{\mathbf T}

\theoremstyle{assumption}%
\newtheorem{assumption}{Assumption}{\bf}{\rm}%

\begin{document}

\title{Learning the intensity of time events with change-points}

\author{Mokhtar Z. Alaya$^{1}$ \and St\'ephane Ga\"iffas$^{2}$ \and
  Agathe Guilloux$^{3}$}

\maketitle

\begin{abstract}
  We consider the problem of learning the inhomogeneous intensity of a
  counting process, under a sparse segmentation assumption. We
  introduce a weighted total-variation penalization, using data-driven weights that correctly scale the penalization along the observation interval. We prove that this leads to a sharp tuning of the convex relaxation of the segmentation prior, by stating oracle inequalities with fast rates of convergence, and
  consistency  for change-points detection. This
  provides first theoretical guarantees for segmentation with a convex
  proxy beyond the standard i.i.d signal + white noise setting. We introduce a fast algorithm to solve this convex problem. Numerical experiments illustrate our approach on simulated and on a high-frequency genomics dataset. \\

  \noindent%
  
\end{abstract}

\footnotetext[1]{Sorbonne Universités, UPMC Univ Paris 06, F-75005, Paris, France, \emph{email}: \texttt{elmokhtar.alaya@upmc.fr}}

\footnotetext[2]{Centre de Mathématiques Appliquées, \'Ecole Polytechnique, and CNRS UMR 7641, 91128 Palaiseau, France, \emph{email}:
  \texttt{stephane.gaiffas@cmap.polytechnique.fr}}

\footnotetext[3]{Sorbonne Universités, UPMC Univ Paris 06, F-75005, Paris, and Unité INSERM 762 ``Instabilit\'e des
  Microsatellites et Cancers'', France \emph{email}:
  \texttt{agathe.guilloux@upmc.fr}}


\section{Introduction}

Counting processes are widely used in engineering to describe systems
where stochastic events occur, such as genomics, biology,
econometrics, communications and networks,
see~\cite{AndBorGilKei-93}. In these problems, the aim is to estimate
the intensity function, which determines the instantaneous rate of
occurrence of an event. In the statistical literature, this topic has
been extensively discussed in several previous works. Procedures based
on kernel estimation~\cite{Ram-83}, cross-validation~\cite{Gre-93},
wavelet methods~\cite{PatWoo-04}, local polynomial
estimators~\cite{CheYipLam-11}, model
selection~\cite{Rey-03}  are
considered for the non-parametric estimation of the intensity.

In this paper, we want to recover the intensity $\lambda_0(t)$ of a
counting process $\{ N(t), t \in [0,1] \}$ from $n$ observations of
$N$. We work under the assumption that $\lambda_0$ can be
well-approximated by a piecewise constant function, and we deal with
this problem with a signal segmentation point-of-view, where the goal
is to find the unknown times of abrupt changes in the dynamic of the
signal. This is referred to \emph{multiple change-point problem} in
statistical literature, see~\cite{KhoAsg-08} for a recent review with
interesting references. A change-point is a time or position where the
structure of the object changes and the goal of change-point detection
is to estimate these positions.

Several examples of practical importance fulfill the model of multiple
change-points. A particularly interesting example comes from the
next-generation sequencing (NGS) DNA process. Indeed, an important
application of NGS technologies is the study of the transcriptome and
the resulting experiment is called RNA-seq. In a typical RNA-seq
experiment, a sample of RNA is amplified, shattered, and converted to
a library of a cDNA fragments. Then, it is sequenced on a
high-throughput platform which is available commercially. Finally,
the raw data result in large amounts of DNA fragments sequences called
reads. These reads are then mapped to the reference genome by an
appropriate algorithm, that tells us the region from which each read
comes from. RNA-seq can be modelled mathematically as replications of
an inhomogeneous counting process with a piecewise constant
intensity~\cite{SheZhan-12}. The counting process counts the number of
reads whose first base maps to the left base of a given chromosome's
location. In~\cite{SheZhan-12}, a Bayesian approach for the detection
of change-points is considered. Other approaches based on Bayesian
model-based clustering and segmentation are given in
\cite{picard2007segmentation}.

In the present paper, we consider the estimation of $\tau_{0,\ell}$
and $\beta_{0,\ell}$ in the following model:
   \begin{equation}
       \label{rate:std}
        \lambda_0(t) = \sum_{\ell=1}^{L_0} \beta_{0,\ell}
        \ind{(\tau_{0,\ell-1}, \tau_{0,\ell}]}(t)
   \end{equation}
for $0 \leq t \leq 1$, with the convention $\tau_{0,0} = 0$ and
$\tau_{0, L_0} =1$. Our approach consists in reframing this task as a
variable selection task. We introduce a penalized least-squares
criterion with a data-driven total-variation penalization, which is
$\ell_1$-penalization of the discrete gradient of the parameter.

This convex proxy for segmentation with an extra $\ell_1$-penalization
for sparsity, called \emph{fused Lasso}, is introduced
in~\cite{TibRosZhuKni-05}. Theoretical guarantees for this procedure
are given in~\cite{HarLev-10} in the white noise setting, for the
segmentation of a one-dimensional signal. A group fused Lasso is
introduced in~\cite{BleVer-11} for the detection of multiple
change-points shared by a set of co-occurring one-dimensional signals,
and an algorithm is derived to solve the corresponding convex
problem. The determination of the number of structural changes
in multitask learning via the group fused Lasso is considered
in~\cite{QiaSu-13}.

Beyond the one-dimensional setting, total-variation penalization is
well-known and commonly used in image denoising, deblurring and
segmentation, see for
instance~\cite{chambolle2009total} and~\cite{chambolle2010introduction}. In this
context, one needs to define a graph of neighboring nodes (pixels),
and the problem can be solved efficiently by reformulating it
as a min-cut problem and solving it using a max-flow
algorithm~\cite{hochbaum2001efficient}.

Other close references are the following:~\cite{GaiGui-12} proves sharp oracle inequalities
for the Lasso in hazards models, ~\cite{Ciu-12} studies Lasso-type estimators in a
linear regression model with multiple change-points,~\cite{Rin-09}
considers denoising of a sparse and block
signal,~\cite{BoyKemLieMunWit-09} studies the asymptotics for
jump-penalized least squares regression aiming at approximating a
regression function by piecewise constant functions.  An algorithm of
majorization-minimization for high dimensional fused Lasso regression
is proposed in \cite{YuWonLeeLimYoo-13}, a testing approach for the
segmentation of the hazard function is given in~\cite{GooLiTiw-11}.


The papers \cite{Rey-03}, \cite{TibRosZhuKni-05}, \cite{HarLev-10}, \cite{BleVer-11}, \cite{QiaSu-13}, are most relevant to our work. In~\cite{Rey-03}, a model selection procedure is introduced to estimate the intensity function. 
In \cite{HarLev-10} and \cite{TibRosZhuKni-05}, the authors propose an adaptation of the Lasso algorithm to detect  change-points in the standard i.i.d signal + Gaussian white noise framework. 
In \cite{BleVer-11} and \cite{QiaSu-13}, the authors use group fused Lasso to solve the structural change-points in linear regression problems. This paper is different from these works in the following aspects. 
First, a main feature of our results is that they are derived for a signal in continuous time, as compared to  \cite{HarLev-10}, \cite{QiaSu-13} and \cite{TibRosZhuKni-05}. Namely, we aim at detecting change-points in the intensity function. 
Hence, this problem is prone to an unavoidable non-parametric bias of approximation by a piecewise constant function, which makes our mathematical analysis very different.  
A second main feature of our results is that we introduce a weighted total-variation penalization, using data-driven weights that correctly scale the penalization along the observation interval. 
This is not necessary in the Gaussian and discrete signal + noise setting from~\cite{HarLev-10} for instance. As a side product, we are able to use the same tuning parameters both for  consistency in oracle inequalities, see Theorems~\ref{thm1} and~\ref{thm2}, and  detection of change-points, see Theorems~\ref{thm3} and~\ref{thm4}.
A third main feature of our approach is that we use a convex surrogate for the sparsity of the discrete gradient of the signal, that can be solved numerically very efficiently, see Section~\ref{section:numerical-experiments-agg}, even for a large signal (using many bins). This is not the case for the approach described in~\cite{Rey-03}, which is based on $\ell_0$ model-selection techniques. Furthermore, our oracle inequalities are sharp in the sense that the leading constant in front of the bias terms is equal to one.


The rest of the paper is organized as follows. In Section~\ref{section:notations-agg}, we
provide basic notations. Then, we present our
estimation procedure.  Section~\ref{section:oracle-inequalities-agg} develops oracle inequalities for the estimator, see Theorems~\ref{thm1} and~\ref{thm2}. Section~\ref{section:change-point-agg} gives results in change-points detection, see Theorems~\ref{thm3} and~\ref{thm4}. Section~\ref{section:numerical-experiments-agg} describes a fast algorithm to solve the convex problem studied in the paper. The proofs of the main statements are gathered in Sections~\ref{sec:proofs-thm1-and-2},~\ref{sec:proof-thm3} and~\ref{sec:proof-thm4}.

\section{Counting processes with a sparse segmentation prior}
\label{section:notations-agg}

Let $(\Omega, \mathcal{F}, \P)$ be a probability space and
$(\mathcal{F}_t)_{0 \leq t \leq 1}$ a filtration satisfying the usual
conditions~\cite{LipShi-89}: increasing, right-continuous and
complete. A \emph{counting process} is a stochastic process $\{ N(t)
\}_{0 \leq t \leq 1}$ which is $(\mathcal{F}_t)$-adapted to the
filtration, with right-continuous and piecewise constant paths almost
surely (a.s.), with jump of size +1 at event times such that
$N(0)= 0$ and $N(t) < \infty$ a.s. The term counting process is
natural: $N(t) - N(s)$ corresponds to the number of events of a
certain type occurring in the interval $(s,t]$. The Poisson process is
the most common example of a counting process, where the jumps occur
randomly and independently of each other on disjoint intervals, see
for instance~\cite{Bre-81} and~\cite{karr1991point} for references on
point processes and their statistical estimation.

Since $N$ is increasing, it is a submartingale, so it follows from the
Doob-Meyer decomposition theorem~\cite{Aal-78}. Namely, $N = \Lambda_0 + M,$
where $\Lambda_0$ is a predictable increasing process called the
compensator of $N$ and $M$ is a $(\mathcal{F}_t)$-martingale. We
assume in the following that
\begin{equation}
  \label{comp}
  \Lambda_0(t) = \E[N(t)] = \int_0^t \lambda_0(s)ds
\end{equation}
for $0 \leq t \leq 1$, where $\lambda_0$ is a non-negative
right-continuous function with left-hand limits called \emph{intensity
  rate} of $N$. Under this assumption, $M(t) = N(t) - \int_0^t
\lambda_0(s)ds$ is a local square-integrable martingale with
\emph{quadratic variation} given by $\inr{M}(t) = \int_0^t
\lambda_0(s) ds$ and \emph{optional variation} $[M](t) = \int_0^t
\lambda_0(s) d N(s)$.

\subsection{Sparse segmentation assumption}

We work under the assumption that the intensity is piecewise constant,
over unknown inhomogeneous intervals of time. From now on,
$\ind{A}$ stands for the indicator function of a set $A$. For
some results in the paper, we will use
\begin{assumption}
  \label{ass:intensity}
  We assume that the intensity writes
  \begin{equation}
    \label{rate}
    \lambda_0(t) = \sum_{\ell=1}^{L_0}\beta_{0,\ell}
    \ind{J_\ell}(t),\, 0 \leq t \leq 1,
  \end{equation}
  with $L_0 \geq 1$, $\beta_{0, \ell}$ are positive coefficients, and
  where $J_0 = \{ 0 \}$, $J_\ell = (\tau_{0,\ell-1}, \tau_{0,\ell}]$
  for $\ell=1,\ldots, L_0$ and $\tau_{0,0} = 0 < \tau_{0,1} < \cdots <
  \tau_{0, L_0-1} < \tau_{0,L_0} = 1$.
\end{assumption}

Assumption~\ref{ass:intensity} means that $L_0 - 1$ changes affect the
value of $\lambda_0$ at unknown instants $\tau_{0,\ell}$. The number
of change-points $L_0 - 1$ is unknown. In this setting, we want to
recover the intensity $\lambda_0$, by jointly estimating $L_0,
\tau_{0,\ell}$ and $\beta_{0,\ell}$, for $\ell=1,\ldots, L_0
-1$. Throughout the paper, we will assume the following.
\begin{assumption}
  \label{ass:iid}
  We observe $n$ i.i.d copies of $N$ on $[0, 1]$, denoted $N_1,
  \ldots, N_n$.
\end{assumption}
The assumption that the process is in $[0, 1]$ is for the sake of
simplicity. Assumption~\ref{ass:iid} is equivalent to observing a
single process $N$ with intensity $n \lambda_0$, which is only used to
have a notion of growing observations with an increasing~$n$.

\subsection{A procedure based on total-variation penalization}

Fix $m = m_n \geq 1$, an integer that shall go to infinity as $n
\rightarrow \infty$. Let us define the set of nonnegative piecewise
constant functions on $[0,1]$ given by
\begin{equation}
  \label{eq:Lambda_m_def}
  \Lambda_m = \Big\{\lambda_\beta =
  \sum_{j=1}^m\beta_{j,m}\lambda_{j,m} : \beta = [\beta_{j,m}]_{1 \leq
    j \leq m} \in \mathbb{R}_+^m \Big\},
\end{equation}
where
\begin{equation*}
  \lambda_{j,m} = \sqrt{m} \ind{I_{j,m}} \quad \text{ and } \quad
  I_{j,m} = \Big( \frac{j-1}{m}, \frac{j}{m} \Big].
\end{equation*}
The linear space $\Lambda_m$ is endowed by the norm $\norm{\lambda} =
(\int_0^1 \lambda^2(t) dt)^{1/2}$. 
 We introduce the
least-squares functional
\begin{equation*}
  R_n(\lambda) = \int_0^1 \lambda(t)^2 dt -
  \frac{2}{n} \sum_{i=1}^n \int_0^1 \lambda(t) dN_i(t),
\end{equation*}
which is the goodness-of-fit criterion to be used in this setting, see
among others~\cite{Rey-03}. Note that $\{ \lambda_{j,m} : j=1,
\ldots, m\}$ produces an orthonormal basis of $\Lambda_m$, it implies that 
\begin{equation*}
  R_n(\lambda_\beta) = \sum_{j=1}^m \beta_{j,m}^2
  -\frac{2 \sqrt{m}}{n} \sum_{j=1}^m \sum_{i=1}^n \beta_{j,m} N_i(I_{j,m})
\end{equation*}
for any $\beta \in \mathbb{R}_+^m$. Now, let us introduce the
\emph{weighted total-variation} penalization
\begin{equation}
  \label{eq:total-variation}
  \norm{\beta}_{\TV, \hat w} = \sum_{j=2}^{m} \hat w_j |\beta_{j} -
  \beta_{j-1}|
\end{equation}
for $\beta = [\beta_j]_{1\leq j \leq m}\in \R^m$, where $\hat w =
[\hat w_j]_{1\leq j \leq m}$ is a positive vector of weights
(eventually depending on data) to be defined later on, with $\hat w_1
= 0.$  The data-driven weights $\hat w $ will allow to design sharp
tuning of the total-variation penalization. Then, given $m \geq 1$ and
a weights vector $\hat w$, we introduce
\begin{equation}
  \label{eq:hat-beta}
  \hat \beta = \argmin\limits_{\beta \in \R_+^m}  \big\{ R_n(\lambda_\beta) +
  \norm{\beta}_{\TV, \hat w} \big\},
\end{equation}
hence an estimator of $\lambda_0$ is given by $\hat \lambda =
\lambda_{\hat \beta}$. An estimation of the change-point locations is
obtained from the support of the discrete gradient of $\hat \beta$. Namely,
define
\begin{equation}
  \label{eq:hat_S}
  \hat{S} = \big\{j:\, \hat{\beta}_{j,m} \neq \hat{\beta}_{j-1,m}
  \text{ for } j=2, \ldots, m \big\},
\end{equation}
and denote by $\hat L = |\hat S|$ the estimated number of change-points.

We denote the mean counting process $\bar N_n = n^{-1} \sum_{i=1}^n
N_i$, and  the unweighted TV penalization by
$\norm{\beta}_{\TV} = \sum_{j=2}^{m} |\beta_{j} - \beta_{j-1}|$ for $\beta
\in \R^m$. We use also the notation $\bar N_n(I) = \int_I d \bar
N_n(t)$ for any $I \subset [0, 1]$.

\section{Sharp oracle inequalities}
\label{section:oracle-inequalities-agg}

In this section we address the statistical properties of $\hat \lambda$
stated in~\eqref{eq:hat-beta}, by proving two oracle
inequalities. Theorem~\ref{thm1} below is an oracle inequality of
``slow-type'' \cite{BicRitTsy-09} that holds in full generality, while
Theorem~\ref{thm2} is a fast oracle inequality, that holds under the
assumption that the number of the estimated change-points is upper bounded by a
known constant $L_{\max}$. Both oracle inequalities are sharp in the
sense that the constant term in front of the oracle term $\inf_{\beta}
\norm{\lambda_\beta - \lambda}$ is equal to one.

\begin{thm}
  \label{thm1}
  Fix $x > 0$ and  introduce the data-driven weights, 
  \begin{equation*}
    \hat{w}_j = 5.66 \sqrt{\frac{m (x + \log m +
        \hat{h}_{n,x,j}) \hat V_j}{n}} + 9.31 \frac{\sqrt{m} (x + 1 +
      \log m + \hat{h}_{n,x,j})}{n},
  \end{equation*}
  where $\hat V_j = \bar N_n \big( \big( \frac{j-1}{m}, 1 \big] \big)$
  and
  \begin{equation*}
    \hat h_{n,x,j} = 2 \log \log \Big( \frac{6 e n \hat V_j + 14 e  (x
      + \log m)}{28 (x + \log m )  } \vee e \Big).
  \end{equation*}
  Then, if $\hat \lambda$ is given by~\eqref{eq:hat-beta}, we have
  \begin{equation}
    \label{slowrate}
    \norm{\hat \lambda - \lambda_0}^2 \leq  \inf_{\beta \in \R_+^m }
    \Big(\norm{\lambda_\beta - \lambda_0}^2 + 2 \norm{\beta_{}}_{\TV,
      \hat w} \Big)
  \end{equation}
  with a probability larger than $1 - 12.85 e^{-x}$.
\end{thm}

The proof of Theorem~\ref{thm1} is postponed in
Section~\ref{sec:proofs-thm1-and-2}. We define $\beta_{0,m} =
[\beta_{0,j,m}]_{1\leq j\leq m}$ the coefficients vector of the
projection of $\lambda_0$ on $\Lambda_m$ and $\Delta_{\beta, \max} =
\max\limits_{1 \leq \ell, \ell' \leq L_0} |\beta_{0,\ell} -
\beta_{0,\ell'}|$, which is the maximum jump size of
$\lambda_0$. Under Assumption~\ref{ass:intensity}, a control of the
approximation term leads to the following.
\begin{cor}
  \label{cor:slowrate}
  Given Assumption~\ref{ass:intensity}, and under the same assumptions
  as the ones from Theorem~\ref{thm1}, we have
  \begin{equation}
    \norm{\hat \lambda - \lambda_0}^2 \leq \frac{2 (L_0 - 1)
      \Delta_{\beta, \max}^2}{m} + 2 \norm{\beta_{0,m}}_{\TV} \max\limits_{1\leq j\leq m} \hat w_j.
  \end{equation}
\end{cor}
The proof of Corollary~\ref{cor:slowrate} is given in
Section~\ref{sec:proofs-thm1-and-2}. Theorem~\ref{thm1} uses a
data-driven weighting of the TV penalization, based on weights roughly
given by
\begin{equation}
\label{data-driven}
  \hat{w}_j \approx \sqrt{\frac{m \log m}{n} \bar N_n\Big(
    \Big(\frac{j-1}{m}, 1 \Big] \Big)}.
\end{equation}
This exhibits a new scaling of the TV penalization, which is natural
and of importance in this setting. The shape of this data-driven
weighting comes from a Bernstein's concentration with data-driven
variance, necessary for the control of the noise term (a martingale
with jumps), given in Proposition~\ref{prop1} below, see
Section~\ref{sec:proof-of-thm1}.
\begin{thm}
  \label{thm2}
  Fix $x > 0$ and let $\hat \lambda$ be the same as in
  Theorem~\ref{thm1}. Assume that the estimated number of
  change-points $\hat L$ satisfies $\hat L \leq L_{\max}$. Then, we
  have
  \begin{equation}
    \label{fastrate}
    \begin{split}
      \norm{\hat \lambda - \lambda_0}^2 \leq& \inf_{\beta \in
        \R_+^m}\big\|\lambda_\beta - \lambda_0\big\|^2+
      6({L_{\max}+2(L_0-1)})\max_{1\leq j\leq m}\hat{w}_j^2 \\
      & \quad + K_1 {\frac{\|\lambda_0\|_{\infty}
          \big(x +L_{\max} (1+ \log m) \big)}{n}} \\
      & \quad + K_2 \frac{m \big(x + L_{\max} (1+ \log m)
        \big)^2}{n^2},
    \end{split}
  \end{equation}  
  with a probability larger than $1 - L_{\max}e^{-x}$, with
  $\|\lambda_0\|_\infty = \sup_{ t\in[0,1]}\lambda_0(t),$  $K_1 =
  1670.89,$ and $K_2 = 6683.53$.
\end{thm}
The proof of Theorem~\ref{thm2} is  provided in
Section~\ref{sec:proofs-thm1-and-2}. This results proves that our
procedure has a fast rate of convergence of order
\begin{equation*}
  \frac{(L_{\max} \vee L_0) m \log m}{n},
\end{equation*}
which scales in $m / n$.
\begin{cor}
  \label{cor:fastrate}
  Given Assumption~\ref{ass:intensity}, and under the same assumptions
  as the ones from Theorem~\ref{thm2}, we have
  \begin{equation}
    \begin{split}
      \norm{\hat \lambda - \lambda_0}^2 \leq&
      \frac{2(L_0-1)\Delta_{\beta,\max}^2} {m}
      +6({L_{\max}+2(L_0-1)})\max_{1\leq j\leq m}\hat{w}_j^2 \\
      & \quad + K_1 {\frac{\|\lambda_0\|_{\infty}
          \big(x +L_{\max} (1+ \log m) \big)}{n}} \\
      & \quad + K_2 \frac{m \big(x + L_{\max} (1+ \log m)
        \big)^2}{n^2},
    \end{split}
  \end{equation}
  with a probability larger than $1 - L_{\max}e^{-x}$, with the same
  notations as in Theorem~\ref{thm2}.
\end{cor}
The proof of Corollary~\ref{cor:fastrate} is presented in
Section~\ref{sec:proofs-thm1-and-2}. A consequence of
Corollary~\ref{cor:fastrate} is that an optimal tradeoff between
approximation and complexity is given by the choice $m \approx
n^{1/2}$. Note that we are able to use the same procedure in
Theorems~\ref{thm1} and~\ref{thm2}, namely for the slow and fast rate,
while it is not the case in the signal + white noise considered
in~\cite{HarLev-10} for instance.

\section{Change-point detection}
\label{section:change-point-agg}

In this section we prove that the proposed total-variation with
data-driven weights procedure is consistent for the estimation of
 the change-point positions. Note that, however, the context
considered here is quite different from the more standard signal +
white noise setting: here we aim at detecting change-points in the
intensity function, hence this problem is prone to an unavoidable
non-parametric bias of approximation by a piecewise constant function.
This means that we will not be able to recover the exact position of
two change-points if they lie on the same interval
$I_{j,m}$. Therefore, we assume
\begin{assumption}
  \label{ass:min-dist}
  Grant Assumption~\ref{ass:intensity} and assume that there is a
  positive constant $c \geq 8$ such that
  \begin{equation}
  \label{ass:min-dist-ineq}
    \min_{1\leq \ell \leq L_0} |\tau_{0,\ell} - \tau_{0,{\ell-1}}| >
    \frac{c}{m}.
  \end{equation}
\end{assumption}
This assumption entails that the change-points of $\lambda_0$ are
sufficiently far apart, and that, in particular, there cannot be more
than one change-point in the ``high-resolution'' intervals $I_{j,
  m}$. Under Assumption~\ref{ass:min-dist}, the procedure will be able
to recover the (unique) intervals $I_{j_\ell, m}$, for $\ell = 0,
\ldots, L_0$, where the change-point belongs. Hence, we define the
\emph{approximate change-points sequence} $[j_\ell]_{0\leq \ell\leq
  L_0}$ as follows.

\begin{defn}
  \label{def:approx_change_points}
  The \emph{approximate change-points sequence} $[j_\ell]_{0\leq \ell \leq L_0}$ relative to the level of resolution $m$ is defined
  as the right-hand side boundary of the unique interval $I_{j_\ell,m}$ that
  contains the change-point $\tau_{0, \ell}$, namely
  \begin{equation}
\label{positon-of-true-change}
    \tau_{0,\ell} \in \Big(\frac{j_{\ell} - 1}{m}, \frac{j_\ell}{m}
    \Big]
  \end{equation}
  for $\ell=1, \ldots, L_{0} - 1$, where we put $j_0 = 0$ and $j_{L_0} =
  m$ by convention.
\end{defn}

Given the support $\hat S = \{ \hat j_1, \ldots, \hat j_{\hat L} \}$
with $\hat j_{1} < \cdots < \hat j_{\hat L}$ of the discrete gradient of $\hat
\beta$ defined in~\eqref{eq:hat_S}, and introducing $\hat j_0 = 0$ and
$\hat j_{\hat L + 1} = m$, we define simply
\begin{equation}
  \label{eq:hat_tau_def}
  \hat{\tau}_\ell = \frac{\hat{j}_\ell}{m}
\end{equation}
for $\ell = 0, \ldots, \hat L + 1$. In order to be able to prove a
consistency results for change-points detection, we need a set of
assumptions that quantifies the asymptotic interplay between several
quantities:
\begin{itemize}
\item $\Delta_{j,\min} = \displaystyle\min_{1\leq \ell \leq L_0-1}|j_{\ell+1} -
  j_{\ell}|$, which is the minimum distance between two consecutive
  terms in the change-points of $\lambda_0.$
\item $\Delta_{\beta,\min} = \min\limits_{1\leq q \leq m-1}|\beta_{0,q+1,m} -
  \beta_{0,q,m}|$, which is the smallest jump size of the projection
  $\lambda_{0,m}$ of $\lambda_0$ onto $\Lambda_m$.
\item $(\varepsilon_n)_{n \geq 1}$, a non-increasing and positive
  sequence that goes to zero as $n \rightarrow \infty$, and such that
  $m \varepsilon_n \geq 6$ for any $n \geq 1$.
\end{itemize}

\begin{assumption}
  \label{ass:consistency}  
  We assume that $\Delta_{j,\min}$, $\Delta_{\beta,\min}$ and
  $(\varepsilon_n)_{n \geq 1}$ satisfy
  \begin{align}
    \label{ass:consistency-1}
    \frac{\sqrt{n m} \varepsilon_n \Delta_{\beta, \min}}{\sqrt{\log m}} &\rightarrow \infty \\
    \label{ass:consistency-2}
    \frac{\sqrt{n} \Delta_{j,\min} \Delta_{\beta,\min}}{\sqrt{m \log m}}
    &\rightarrow \infty    
  \end{align}
  as $n \rightarrow \infty$.
\end{assumption}

This assumption controls the rate $(\varepsilon_n)$ of convergence of
$\hat{\tau}_\ell$ towards $\tau_{0,\ell}$. The logarithmic factor is
due to concentration inequalities for the control of the noise (the
martingale $M$ obtained by compensation of $N$). The next Theorem
proves the consistency of our procedure for the detection of
change-points, under the assumption that the estimated number of
change-points is the correct one.

\begin{thm}
  \label{thm3}
  Under Assumptions~\ref{ass:min-dist} and~\ref{ass:consistency}, and
  if $\hat L = L_0 - 1$, then the change-points estimators
  $\{\hat{\tau}_1, \ldots, \hat{\tau}_{\hat L}\}$ given
  by~\eqref{eq:hat_tau_def} satisfy
  \begin{equation}
    \label{estchgpt}
    \P \Big[ \max_{1\leq \ell\leq L_0 - 1} |\tau_{0,\ell} -
    \hat{\tau}_\ell| \leq {\varepsilon_n} \Big] \rightarrow 1
  \end{equation}
  as $n \rightarrow \infty$.
\end{thm}

The proof of Theorem~\ref{thm3} is quite involved and is presented in
Section~\ref{sec:proof-thm3} and Appendix~\ref{app:proof-thm3-case2}. It builds upon some techniques
developed in~\cite{HarLev-10}, based on a careful inspection of the
Karush-Kuhn-Tucker (KKT) optimality conditions, see for
instance~\cite{BoyVan-04}, for the solutions to the convex
problem~\eqref{eq:hat-beta}. The proof depends also heavily on a
data-driven Bernstein's inequality for the control of the martingale
errors, see Proposition~\ref{prop1} from
Section~\ref{sec:proofs-thm1-and-2}.

Let us give examples of scaling for the quantities $\Delta_{j,\min}$,
$\Delta_{\beta,\min}$ and $(\varepsilon_n)_{n \geq 1}$ that meet
Assumption~\ref{ass:consistency}. Assume for simplicity that
\begin{equation*}
  \varepsilon_n = n^{-\alpha} \quad \text{ and } \quad
  \Delta_{\beta,\min} = n^{-\gamma}
\end{equation*}
for some constants $\alpha, \gamma > 0$.
\begin{itemize}
\item If $m = n^{1 / 3}$ then Theorem~\ref{thm3} holds with any
  $\alpha, \gamma > 0$ satisfying $0 < \gamma < 1 / 3$ and $0 < \alpha
  + \gamma < 2/ 3$, and if $\Delta_{j,\min} \geq 6$.
\item If $m = n^{1 / 2}$ then Theorem~\ref{thm3} holds with any $0 <
  \gamma < 1 / 4$ and $0 < \alpha + \gamma < 3 / 4$ and if
  $\Delta_{j,\min} \geq 6$.
\end{itemize}

In order to prove change-point consistency without the assumption that
the estimated number of change-points is the correct one, we need to
relax a little bit the statement of the result given in
Theorem~\ref{thm3}. Namely, we evaluate a non-symmetrized Hausdorff
distance $\mathcal{E}(\hat{\mathcal{T}} \| \mathcal{T}_0)$ between the
set of estimated change-points
\begin{equation*}
  \hat{\mathcal{T}} = \big\{ \hat{\tau}_1, \ldots, \hat{\tau}_{\hat L}
  \big\}
\end{equation*}
and the set of true change-points 
\begin{equation*}
  \mathcal{T}_0 = \big\{\tau_{0,1}, \ldots, \tau_{0,L _0 - 1} \big\}, 
\end{equation*}
 where for two sets $A$ and $B$, the
quantity $\cE(A \| B)$ is given by
\begin{equation*}
  \cE(A\|B) = \sup_{b\in B} \inf_{a\in A}|a - b|.
\end{equation*}
Note that $\cE(A\|B) \vee \cE(B\|A)$ is the Hausdorff distance between
$A$ and $B$. When $\hat L = L_0 - 1$, Theorem~\ref{thm3} implies that
\begin{equation}
  \label{estchgpt}
  \P \Big[ \cE \big( \hat \cT \| \cT_0 \big) \leq \varepsilon_n,
  \cE(\cT_0 \|\hat \cT) \leq \varepsilon_n \Big] \rightarrow 1
\end{equation}
as $n \rightarrow \infty$. When $\hat L > L_0 - 1$, we prove in
Theorem~\ref{thm4} below that $\cE (\hat \cT \| \cT_0) \leq
\varepsilon_n$ with a probability going to $1$ as $n \rightarrow
\infty$. This means that change-point consistency holds for our
procedure whenever the estimated number of change-points is not less
than the true one.
\begin{thm}
  \label{thm4}
  Under Assumptions~\ref{ass:min-dist} and~\ref{ass:consistency}, and
  if $\hat L \geq L_0 - 1$, we have
  \begin{equation}
    \label{nbrchgpt}
    \P\Big[ \cE (\hat \cT \| \cT_0 ) \leq \varepsilon_n\Big]
    \rightarrow 1
  \end{equation}
  as $n \rightarrow \infty$.
\end{thm}

Theorem~\ref{thm4} ensures that even when the number of change-points
is over-estimated, each true change-point is close to the estimated one.
The proof of Theorem~\ref{thm4} is given in
Section~\ref{sec:proof-thm4}. It is based, as for the proof of
Theorem~\ref{thm3}, on a repeated utilization of the KKT optimality
conditions of problem~\eqref{eq:hat-beta}.

Note that a difference with~\cite{HarLev-10} is that we are able to
use the same regularization parameters $\hat{w}_j$ given
by~\eqref{data-driven} in Theorems~\ref{thm3} and~\ref{thm4}.
Besides, we don't need an upper bound
on the estimated number of change-points in Theorem~\ref{thm4}, while it is necessary in~\cite{QiaSu-13}.

\section{Numerical experiments}
\label{section:numerical-experiments-agg}

In this section we propose a fast algorithm for solving the optimization problem~\eqref{eq:hat-beta} and apply it on simulated and real datasets from genomics.

\subsection{Algorithm}
\label{sec:algorithm}

A concept of importance for convex optimization in machine learning is the proximal operator~\cite{bach2012optimization, BauCom-11}.
The proximal operator $\prox_{f}$ of a proper, lower semicontinuous, convex function $f: \R^m \rightarrow (-\infty,  \infty],$ is defined as
\begin{equation*}
\prox_{ f}(v) = \argmin_{x \in \R^m}\Big\{\frac{1}{2}\norm{v - x}_2^2 + f(x) \Big\}, \textrm{ for all } v \in \R^m.
\end{equation*}
In this section, we provide a fast algorithm to solve the optimization problem~\eqref{eq:hat-beta}, that computes the proximal operator of the weighted total-variation.

We observe $n$ i.i.d observations of $N$ over the interval
$[0,1]$. Recall that $\bar N_n = n^{-1} \sum_{i=1}^n N_i,$ and $\bar N_n(I) = \int_I d \bar
N_n(t)$ for any $I \subset [0, 1]$.
We also recall that $\hat \lambda(t) = \sum_{j=1}^m \hat \beta_j
\lambda_{j,m}(t)$, where $\hat \beta = [\hat \beta_{1}, \ldots, \hat
\beta_{m}]$ is given by~\eqref{eq:hat-beta}. Hence, we have
\begin{equation}
  \label{estvec2}
  \hat \beta = \argmin_{\beta \in \R_+^m} \Big\{ \frac{1}{2} \norm{\bN
    - \beta}_2^2 + \norm{\beta}_{\TV, \hat w} \Big\},
\end{equation}
where $\bN =[\bN_j]_{1\leq j \leq m}\in \R_+^m$ is given by
\begin{equation*}
  \bN =
  \begin{bmatrix}
    \sqrt{m} \bar N_n(I_{1, m}) \\
    \vdots \\
    \sqrt{m} \bar N_n(I_{m, m})
  \end{bmatrix}
  .
\end{equation*}
Therefore, we see that  \eqref{estvec2} is equivalent to 
\begin{equation*}
   \hat \beta = \prox_{\norm{\cdot}_{\TV, \hat w}} (\bN).
\end{equation*}
Next, we develop an algorithm that computes $\prox_{\norm{\cdot}_{\TV, \hat w}}$, which is an extension of~\cite{Cond-13} to weighted total-variation. Towards this end, we introduce the following $(m-1)\times m$  bidiagonal matrix
\begin{equation*}
   D_{\hat w}= \begin{bmatrix}
   -\hat w_2 & \hat w_2 & 0 & \cdots &0 \\
   0 &-\hat w_3 & \hat w_3 & \ddots &\vdots\\
   \vdots& \ddots& \ddots &\ddots &0 \\
   0& \cdots& 0& -\hat w_m & \hat w_m 
\end{bmatrix}.
\end{equation*}
Then, one can  express the primal problem \eqref{estvec2} as follows:
\begin{equation}
  \label{primalbeta}
  \hat \beta = \argmin_{\beta \in \R_+^m} \Big\{ \frac{1}{2}  \norm{\bN
    - \beta}_2^2 + \norm{D_{\hat w }\beta}_{1} \Big\}.
\end{equation}
Essentially, problem \eqref{primalbeta} is difficult to analyse directly because the nondifferentiable $\ell_1$ norm is composed with a linear transformation of $\beta.$  When solving \eqref{primalbeta} we may consider its Fenchel dual form\cite{BauCom-11}. First, we rewrite the primal problem as 
\begin{equation*}
\begin{split}
&\minimize\limits_{ \beta \in \R^m,\, z \in \R^{m-1}} \frac{1}{2} \norm{\bN - \beta}_2^2 + \norm{z}_1\\ 
&\textrm { subject to } D_{\hat w} \beta = z, 
\end{split}
\end{equation*}
whose Lagrangian is
\begin{equation*}
   \mathscr{L}_{}( \beta, z, u) = \frac{1}{2} \norm{ \bN - \beta}_2^2 +   \norm{z}_1 + u^\top (D_{\hat w}\beta - z),
\end{equation*}
and to derive a dual problem, we minimize this over $\beta, z.$  A straightforward computation gives
\begin{equation*}
   \min\limits_\beta\Big\{ \frac{1}{2} \norm{ \bN - \beta}_2^2  + u^\top  D_{\hat w}\beta \Big\} = - \frac{1}{2} \norm{ \bN - D_{\hat w}^\top u}_2^2,
\end{equation*}
while
\begin{equation*}
  \min\limits_z\Big\{\norm{z}_1 - u^\top z\Big\} = \left\{ \begin{array}{ll}
  0,  & \mbox{if }  \norm{u}_\infty  \leq 1,\\
  -\infty,   & \mbox{otherwise.}
  \end{array}
  \right.
\end{equation*}
Introducing  $u_0 = u_m = 0,$ we proved that a dual problem of \eqref{primalbeta} is given by
\begin{equation*}
\begin{split}
 &\minimize_{u\in \R^{m+1}} \frac{1}{2}\sum_{k=1}^m \big( \bN_k- \hat w_{k+1} u_k+ \hat w_{k}u_{k-1}\big)^2, \\
 &\textrm{ subject to }|u_j| \leq 1, \textrm{ for } k =1, \ldots, m, \textrm{ and } u_0=u_m=0.
\end{split} 
\end{equation*}
If we have a feasible dual variable $\hat u,$ we can compute the primal solution $ \hat \beta$ using 
\begin{equation}
\label{fencheldual}
\hat\beta_k = \bN_k- \hat w_{k+1} \hat u_k+ \hat w_{k}\hat u_{k-1}, \textrm{ for } k =1, \ldots, m.
\end{equation}
For this problem, strong duality holds, see~\cite{BoyVan-04}, meaning that the duality gap is zero. The KKT optimality conditions characterize the unique solutions $\hat \beta$ and $\hat \theta_k := \hat w_{k+1} \hat u_k.$ They yield, in addition to \eqref{fencheldual}:
\begin{equation}
\label{system:KKT-wTV-agg}
\hat \theta_0 = \hat \theta_m=0,\,\textrm{and } \forall \, k =1, \ldots, m-1,\\
 \left\{
    \begin{array}{lll}
        \hat \theta_k \in [-\hat w_{k+1}, \hat w_{k+1}], & \mbox{if } \hat \beta_k= \hat\beta_{k+1}, \\
        \hat \theta_k=-\hat w_{k+1},& \mbox{if } \hat \beta_k< \hat \beta_{k+1},\\
        \hat \theta_k=\hat w_{k+1}, & \mbox{if } \hat \beta_k> \hat \beta_{k+1}.
    \end{array} 
\right.
\end{equation}
Therefore, the proposed algorithm consists in running forwardly through the samples $[\bN_k]_{1\leq k\leq  m}.$ Using \eqref{system:KKT-wTV-agg}, at location $k,$  $\hat \beta_k$ stays constant where $|\hat \theta_k| < \hat w_{k+1}.$ If this is not possible, it goes back to the last location where a jump can be introduced in $\hat \beta$, validates the current segment until this location, starts a new segment, and continues. This algorithm is described precisely in Algorithm~\ref{algorithm:weighted-TV-agg}.

\LinesNotNumbered
\begin{algorithm}
\SetNlSty{textbf}{}{.}
\DontPrintSemicolon
 \caption{ $ \hat \beta = \prox_{\norm{\cdot}_{\TV,\hat w}}(\bN)$}
\label{algorithm:weighted-TV-agg}

 \KwIn{$\bN=\big(\bN_1, \ldots, \bN_m\big) ^\top\in \R^m; \hat w = (\hat w_1,\ldots, \hat w_m) \in \R^{m}_{+}.$}

 \KwOut{ $\big(\hat \beta_1,\ldots, \hat \beta_m\big)^\top.$}

 \nl\textbf{Set} {$k=k_0=k_-=k_+ \gets 1;$\\$\qquad \beta_{\min}  \gets \bN_1- \hat w_2;\, \beta_{\max} \gets \bN_1+ \hat w_2;$\\$\qquad \theta_{\min}\gets \hat w_2;\, \theta_{\max} \gets -\hat w_2;$}\\

 \nl \label{step2}{\If{$k=m$}
{  $\hat \beta_m \gets \beta_{\min}+ \theta_{\min};$}}

\nl \label{step3} \If (\tcc*[f]{negative jump}){$\bN_{k+1} + \theta_{\min} < \beta_{\min} - \hat w_{k+2}$}
{$\hat \beta_{k_0}= \cdots =\hat \beta_{k_-} \gets \beta_{\min};$\\
 $k=k_0=k_-=k_+\gets k_- + 1;$\\
 $ \beta_{\min} \gets \bN_{k} - \hat w_{k +1}+ \hat w_{k };\, \beta_{\max}   \gets \bN_{k} + \hat w_{k +1}+ \hat w_{k};$\\
 $\theta_{\min}   \gets \hat w_{k +1};\, \theta_{\max}    \gets -\hat w_{k+1};$
}
 \nl \ElseIf (\tcc*[f]{positive jump}){$\bN_{k+1} + \theta_{\max} > \beta_{\max}+ \hat w_{k+2}$}
{$\hat \beta_{k_0}= \ldots =\hat \beta_{k_+}  \gets \beta_{\max};$\\
$k=k_0=k_-=k_+\gets k_+ + 1;$\\
$\beta_{\min} \gets \bN_{k} - \hat w_{k+1} - \hat w_{k };\, \beta_{\max}\gets \bN_{k} +\hat w_{k+1}-  \hat w_{k};$\\
$ \theta_{\min}   \gets \hat w_{k+1};\, \theta_{\max} \gets -\hat w_{k+1};$
}

 \nl \Else (\tcc*[f]{no jump}){ \textbf{set }$k \gets k+1; $\\$\qquad \theta_{\min} \gets  \bN_{k} + \hat w_{k+1} - \beta_{\min};$\\$\qquad \theta_{\max} \gets  \bN_{k} - \hat w_{k+1} - \beta_{\max};$

  \If{$\theta_{\min}  \geq \hat w_{k+1}$}{ $\beta_{\min}\gets \beta_{\min}+ \frac{\theta_{\min} - \hat w_{k+1}}{k-k_0+1};$\\
$\theta_{\min} \gets \hat w_{k+1};$\\
$k_- \gets k;$}

  \If{$\theta_{\max} \leq -\hat w_{k+1}$}{$\beta_{\max}\gets \beta_{\max}+ \frac{\theta_{\max} + \hat w_{k+1}}{k-k_0+1};$\\
$\theta_{\max} \gets  -\hat w_{k+1};$\\
$k_+ \gets k;$}
}

\nl \If{$ k< m$ } { go to \textbf{\ref{step3}.};}

\nl \If{$\theta_{\min} < 0$ } {  $\hat \beta_{k_0}= \cdots =\hat \beta_{k_-} \gets \beta_{\min};$\\$k=k_0=k_- \gets k_-+1;$\\
$\beta_{\min}\gets \bN_k - \hat w_{k+1}+ \hat w_{k};$\\
$\theta_{\min}\gets \hat w_{k+1};\, \theta_{\max} \gets \bN_k+ \hat w_{k} - v_{\max};$
\\go to \textbf{\ref{step2}.};}

 \nl \ElseIf{$\theta_{\max} > 0$ } { $ \hat \beta_{k_0}= \cdots =\hat \beta_{k_+}\gets \beta_{\max};$\\
$k=k_0=k_+ \gets k_++1;$\\
$\beta_{\max} \gets \bN_k +\hat w_{k+1} -  \hat w_{k};$\\
$\theta_{\max}\gets  -\hat w_{k+1};\,\theta_{\min}\gets \bN_k- \hat w_{k} - \theta_{\min};$ \\ 
go to \textbf{\ref{step2}.};}

\nl \Else{$\hat \beta_{k_0}= \cdots =\hat \beta_m\gets \beta_{\min} + \frac{\theta_{\min}}{k-k_0+1};$  }
\end{algorithm}

\subsection{Simulated data}

\begin{figure}[htbp]
  \centering
  \includegraphics[width=0.45\textwidth]{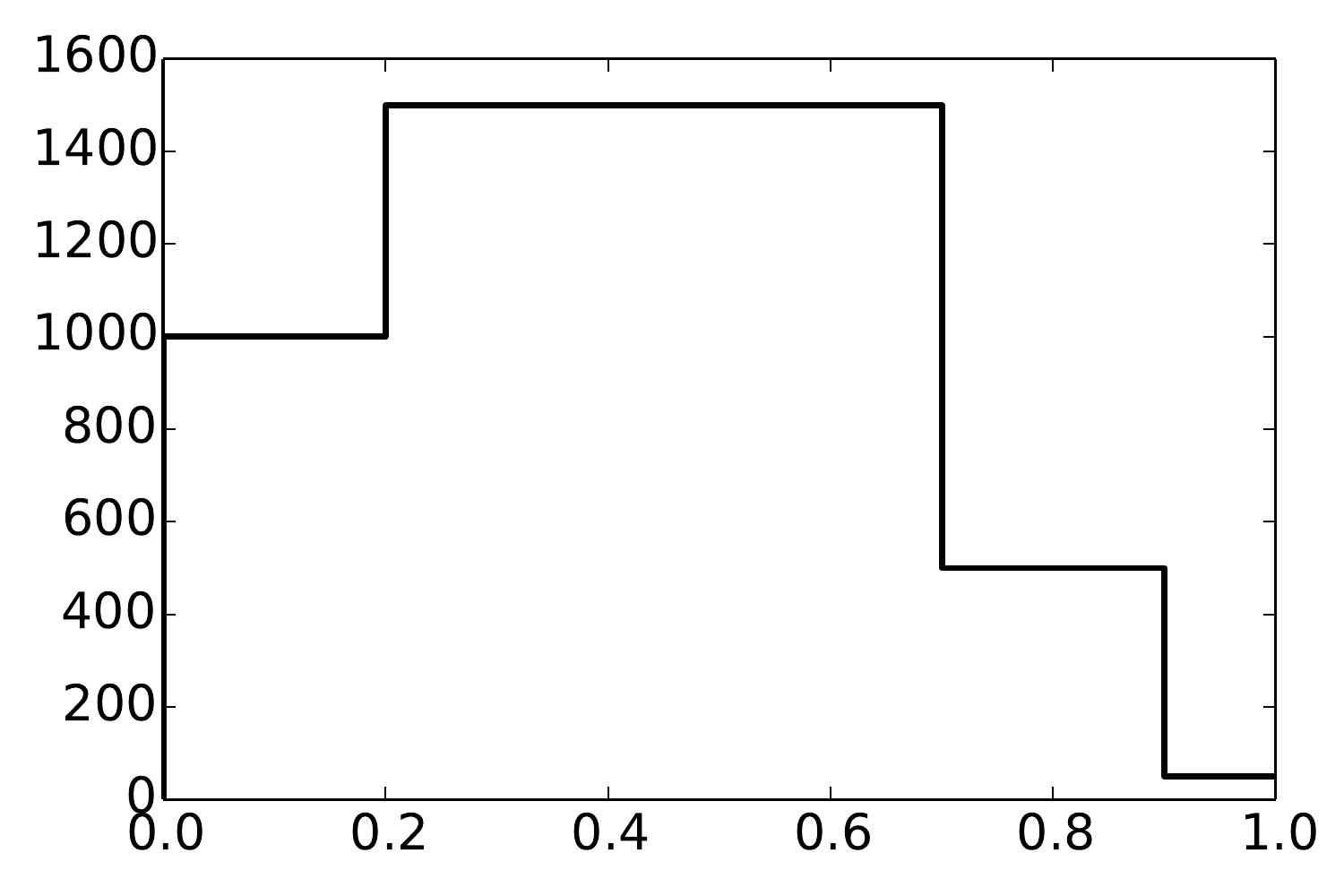}
  \includegraphics[width=0.45\textwidth]{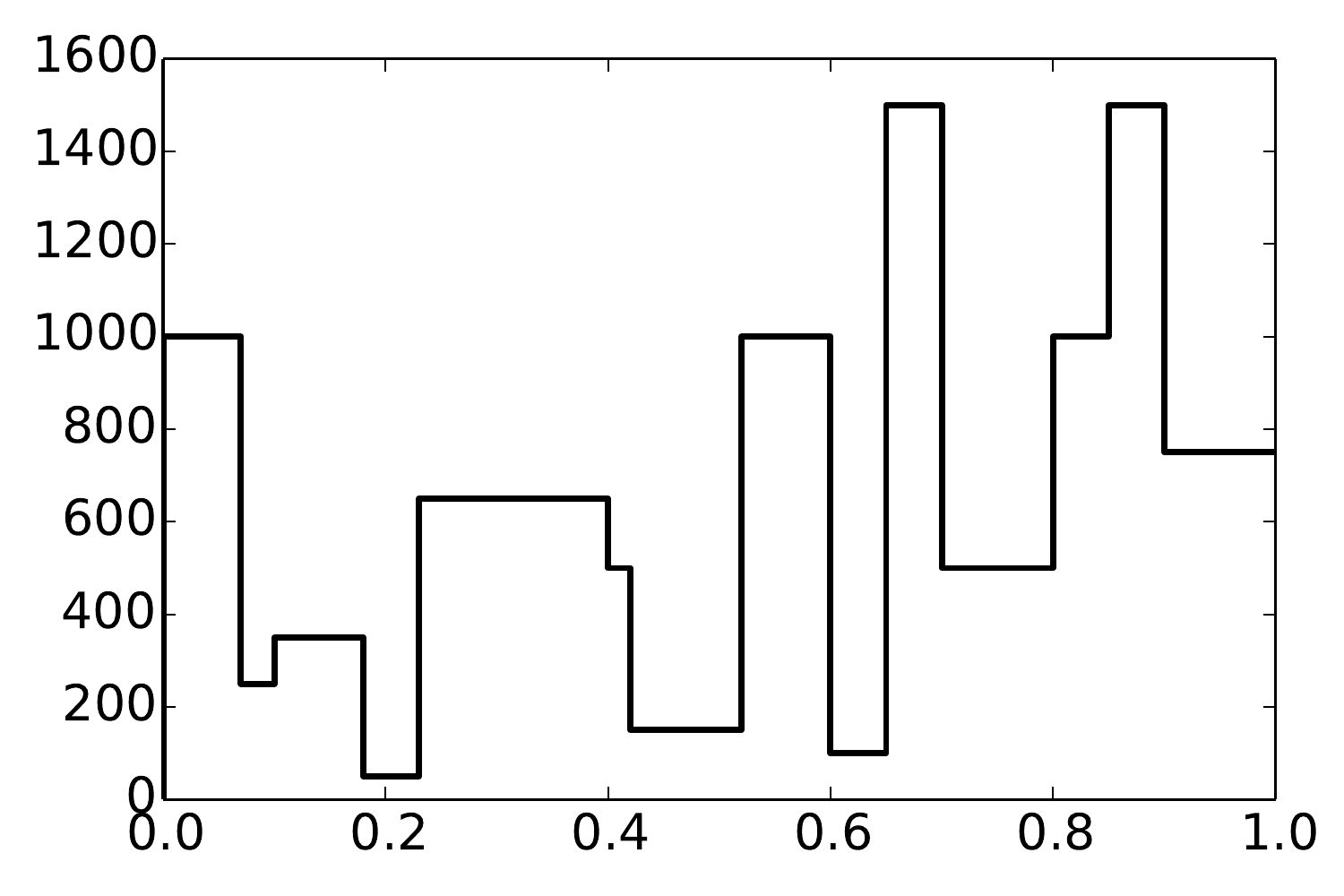}
  \caption{Intensities used for Example~1 (left) and Example~2 (right), respectively with $5$ and $15$ change-points}
  \label{fig:simu_chpts}  
\end{figure}

\begin{figure}[htbp]
  \centering
 \includegraphics[width=0.49\textwidth]{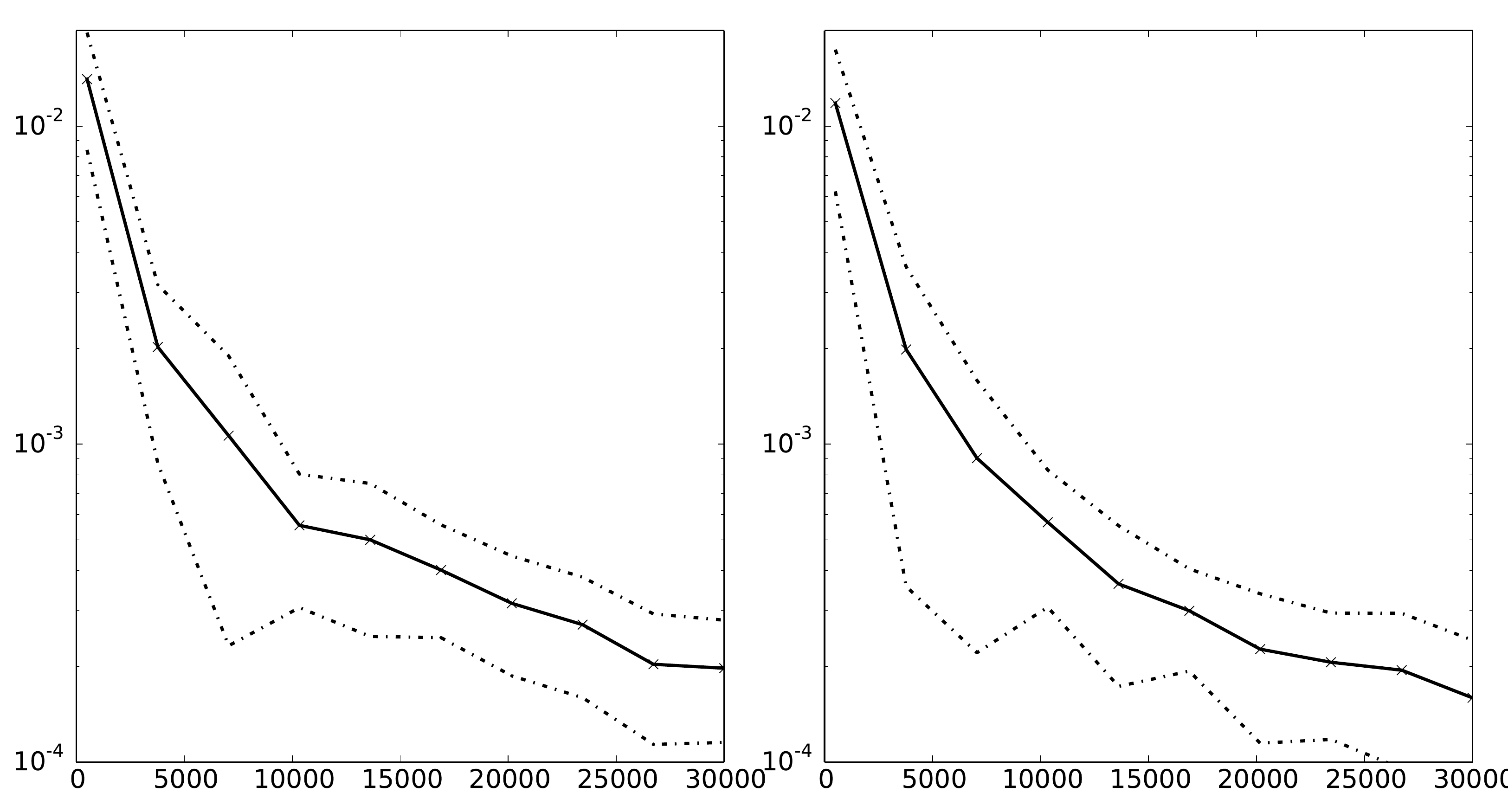}%
 \includegraphics[width=0.49\textwidth]{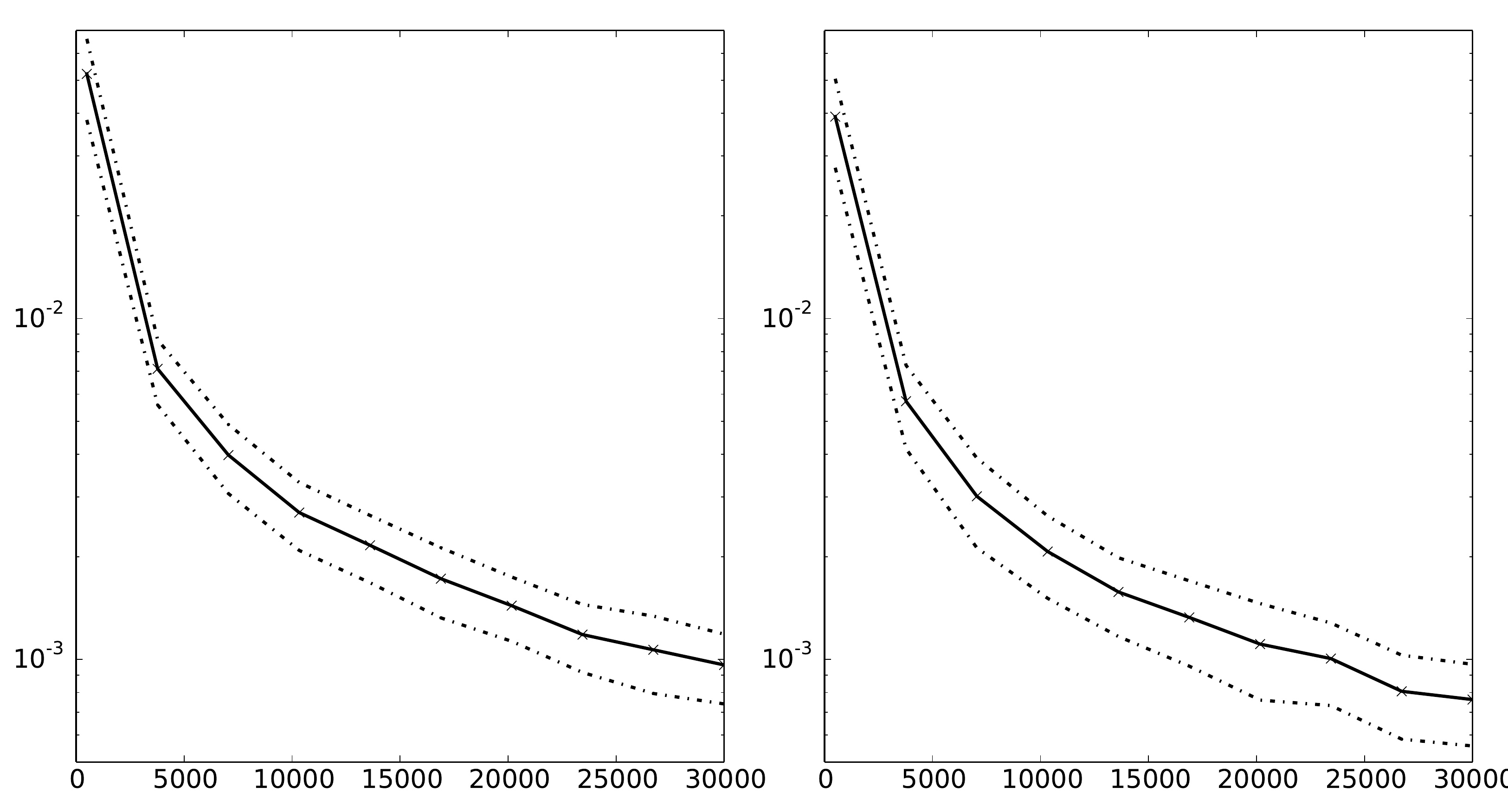}
  \caption{Average MISEs (bold lines) over 100 Monte-Carlo experiments and standard deviations of the MISEs (dashed lines). First: weighted TV for Example~1; Second: non-weighted TV for Example~1; Third: weighted TV for Example~2; Fourth: non-weighted TV for Example~2}
  \label{fig:mises_TV}  
\end{figure}

We conduct simulations on 2 examples of intensities.
We simulate counting processes with inhomogeneous piecewise intensities $\lambda_0,$ with $5$ and $15$ change-points, see Figure~\ref{fig:simu_chpts}, with an increasing sample size $n$.
In order to assess the performance of the total-variation procedure $\hat \lambda,$ we use a Monte-Carlo averaged mean integrated squared error ($\mise$) as a performance measure, given by
\begin{equation*}
\mise(\hat \lambda, \lambda_0) = \E \int_0^1 (\hat \lambda(t) - \lambda_0(t))^2 dt.
\end{equation*}
We run $100$ Monte-Carlo experiments, for an increasing sample size between $n = 500$ and $n=30000$, for each 2 examples.
In Figure~\ref{fig:mises_TV}, we plot the MISEs of the weighted and the unweighted total-variation (namely $\hat w \equiv 1$), for the 2 examples, as a function of the sample size. 
We observe in Figure~\ref{fig:mises_TV} that the estimation error is always decaying with the sample size, and that both procedures behave similarly. Differences can be observed below, using a genomics datasets.
On each simulated dataset, we perform a 10-fold cross-validation to select the best constant to use in front of the weights $\hat w_j$ (both for the weighted and unweighted total-variation). 
Cross-validation in this context is achieved by choosing uniformly at random a label between $1$ and $10$ for each point, and by using points with label $k$ in the $k$-th testing fold and removing these points for the $k$-th training fold. 
The estimated intensity is accordingly corrected, by this amount (as removing uniformly a fraction of points from a counting process biases downwards the intensity by the same fraction).


\subsection{Real data}

\begin{figure}[htbp]
  \centering

  \includegraphics[width=0.45\textwidth]{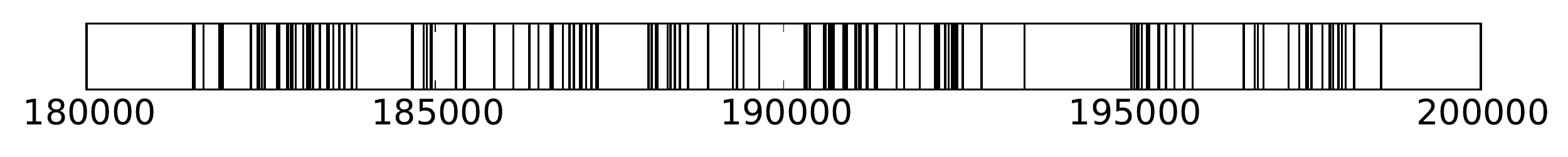} \hspace{1cm} %
  \includegraphics[width=0.45\textwidth]{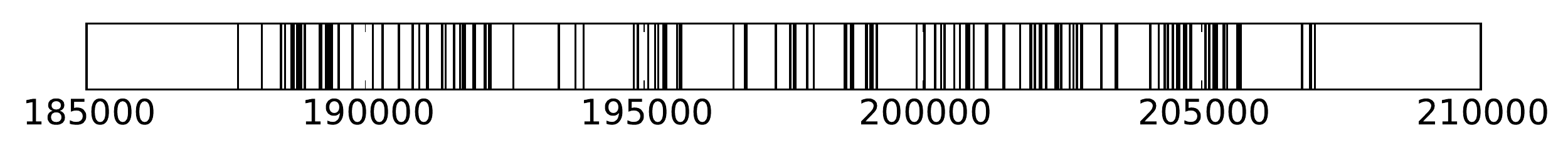}

  \caption{A zoom into the sequence of reads for normal (left) and tumor (right) data}
  \label{fig:zoom_reads}  
\end{figure}

\begin{figure}[htbp]
  \centering

  \includegraphics[width=0.45\textwidth]{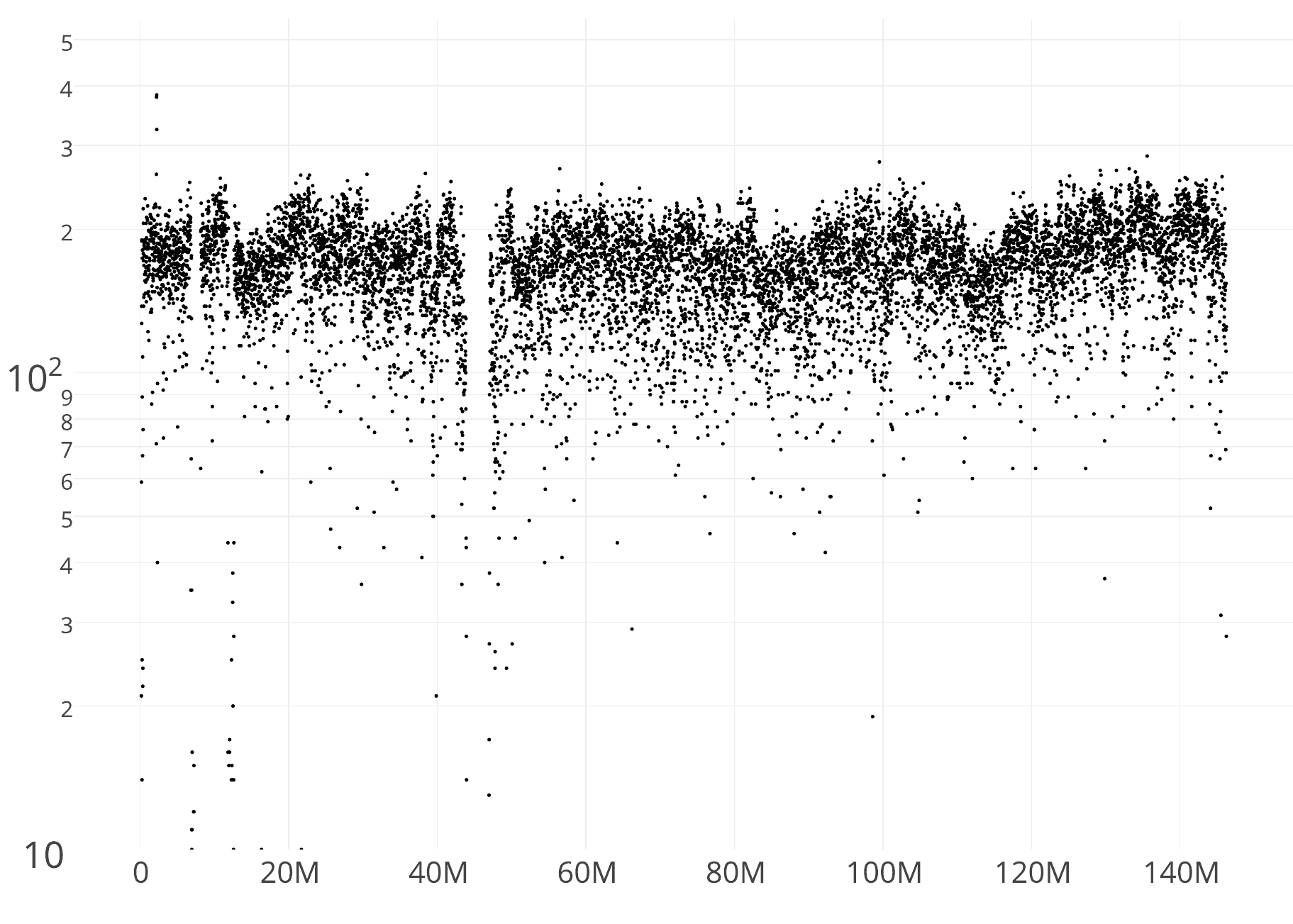} \hspace{1cm} %
  \includegraphics[width=0.45\textwidth]{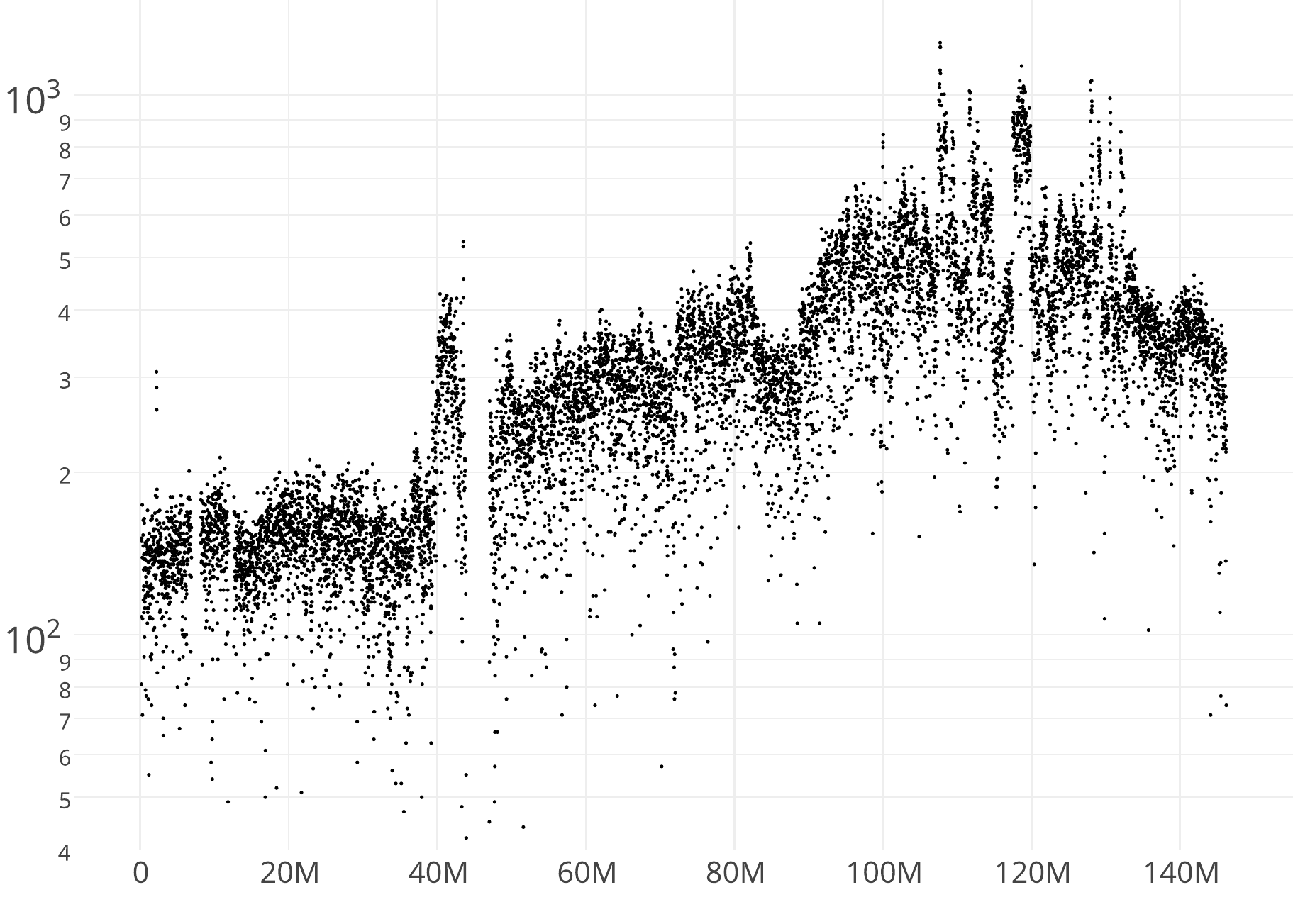}

  \caption{Binned counts of reads (log-scale) of the normal (left) and tumor (right) data} 
  \label{fig:bin_reads}
\end{figure}

Our method is illustrated on NCI-60 tumor and normal cell lines, HCC1954 and BL1954. 
This dataset was produced and investigated by~\cite{chiang2009high} using the Illumina platform, where the reads are 36bp long. 
After cleaning of this data, there are 7.72 million reads for the tumor (HCC1954) and 6.65 million reads for the normal (BL1954) samples respectively. 
A description of the sampling process for such data is described in Introduction. 
We show in Figures~\ref{fig:zoom_reads} and~\ref{fig:bin_reads} both tumor and cell lines data. 
This data consists of a list of reads number, see Figure~\ref{fig:bin_reads}, where we plot a zoomed sequence of reads.
For visualization purposes, we give in Figure~\ref{fig:bin_reads} the binned counts of reads over 10000 intervals equispaced on the range of reads.

\begin{figure}[htbp]
  \centering
   
  \includegraphics[width=0.45\textwidth]{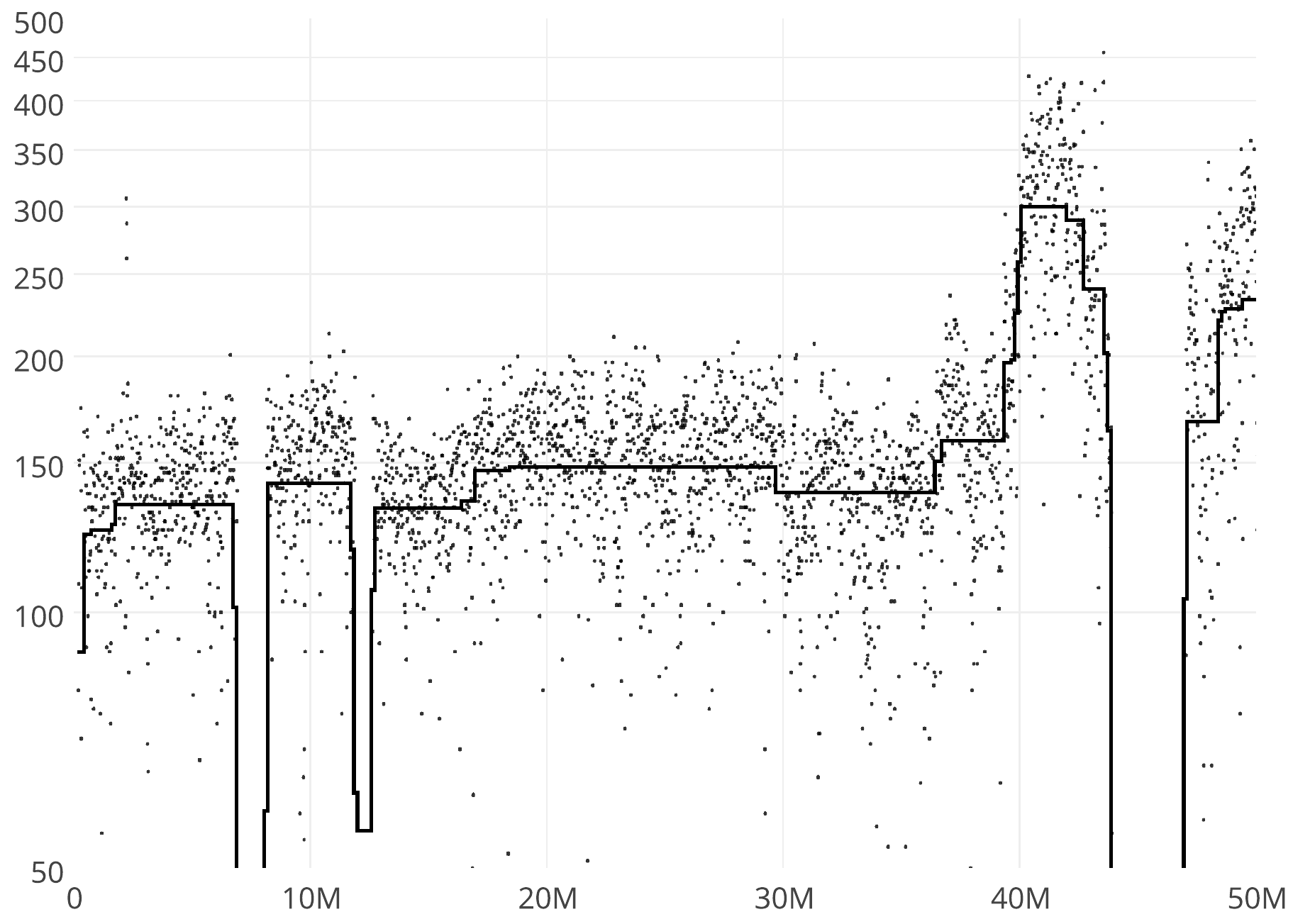} \hspace{1cm} %
  \includegraphics[width=0.45\textwidth]{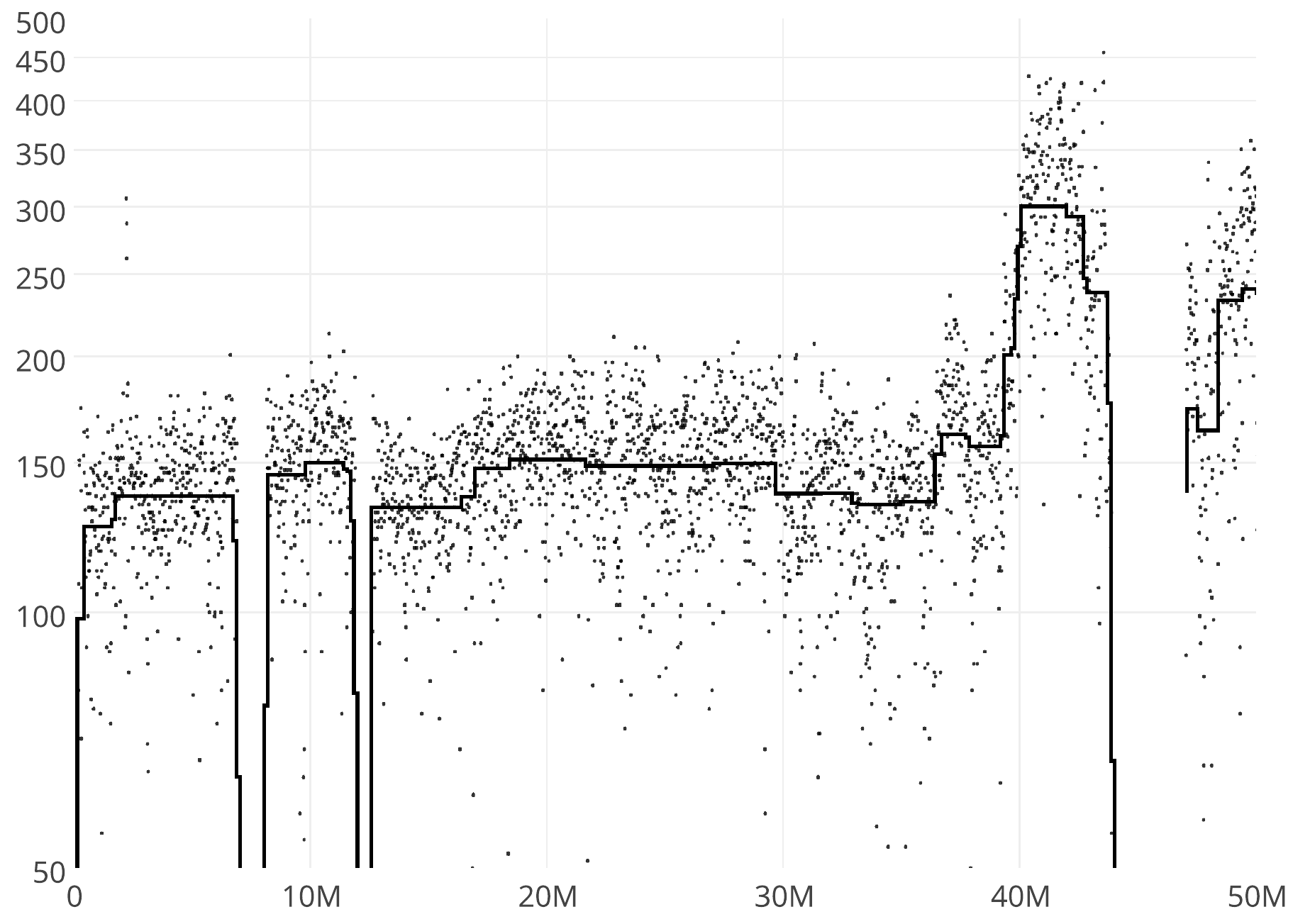} 
   \includegraphics[width=0.45\textwidth]{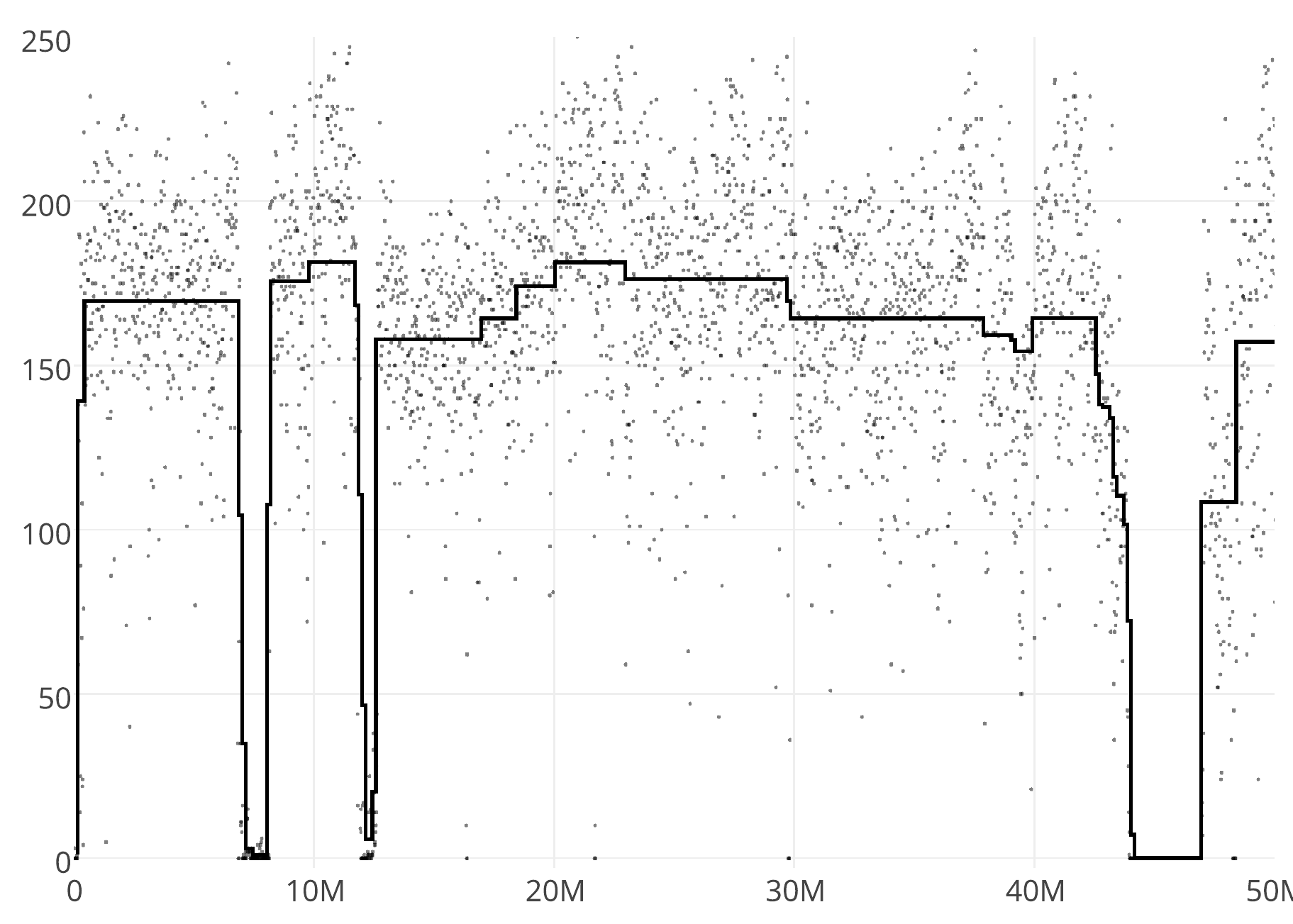} \hspace{1cm} %
  \includegraphics[width=0.45\textwidth]{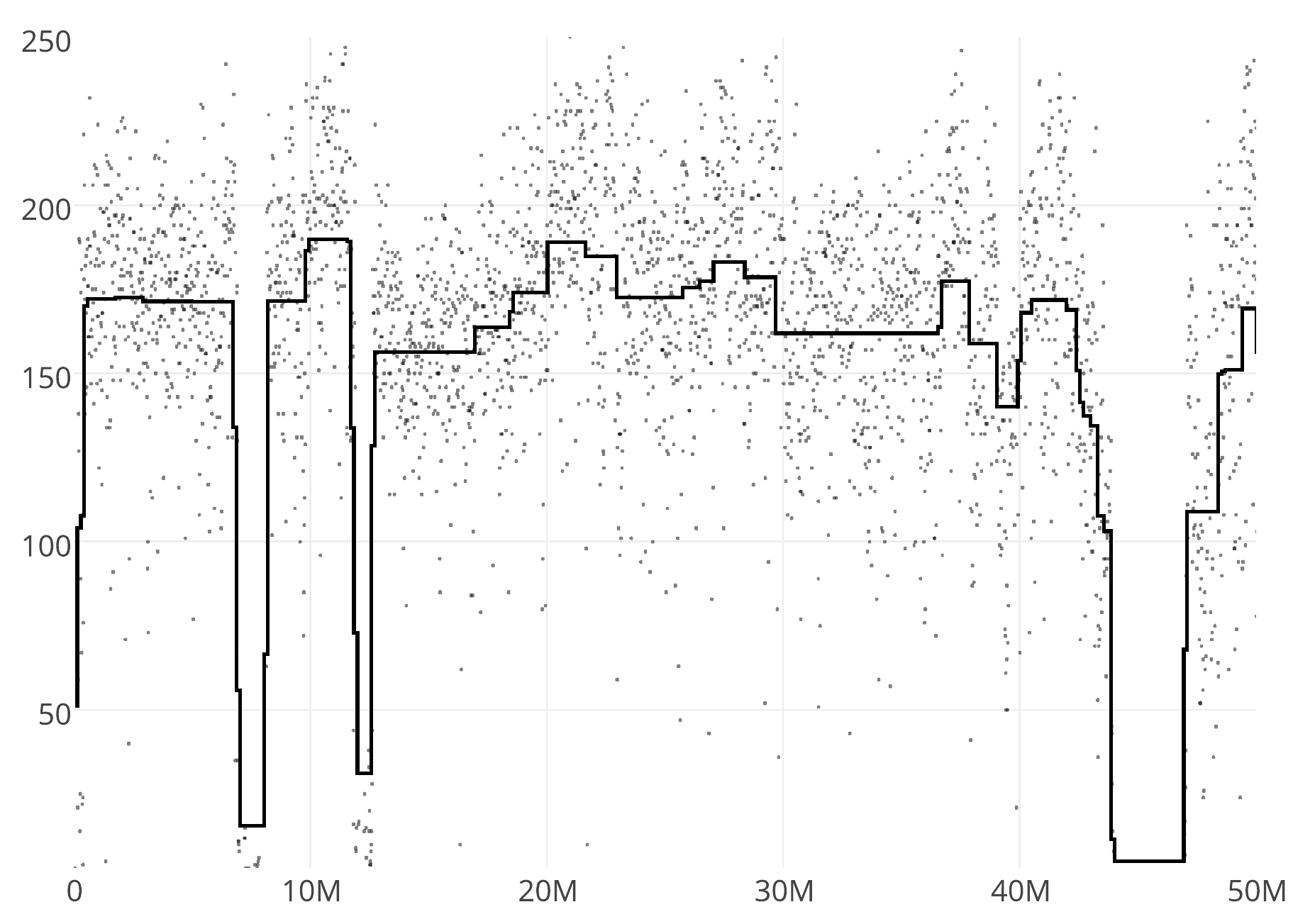} 
  
  \caption{A zoom between reads number 0 and 50M of the weighted (left) and unweighted (right) total-variation estimators applied to the tumor (top) and normal (bottom) data}
  \label{fig:reads_tumor_normal_w_unw}  
\end{figure}

In Figure~\ref{fig:reads_tumor_normal_w_unw}  we plot the best solution of the weighted and unweighted ($\hat w_j=1$) total-variation estimators on the normal and tumor reads data. 
For easier visualization we plot a zoom of the reads sequence.
We perform a 10-fold cross-validation to select the best constant to use in front of the weights $\hat w_j$ (both for the weighted and unweighted total-variation), as explained above.
We observe in this figure that the weighted total-variation gives sharper results: the piecewise constant intensity is smoother, and the obtained change-points locations seem, at least visually, better.
An important fact is that the runtime of Algorithm~\ref{algorithm:weighted-TV-agg} is extremely fast: a solution is obtained in less than one millisecond, on a modern laptop (implementation is done using python with a C extension). This is due to the fact that Algorithm~\ref{algorithm:weighted-TV-agg} is typically linear in the signal size.

\section{Conclusion}

In this work, we prove that convex optimization for the detection of
change-points in the intensity of a counting process is a powerful
tool. We introduce a data-driven weighted total-variation penalization
for this problem, with sharply tuned regularization parameters, and
prove two families of theoretical results: oracles inequalities for
the prediction error, and consistency in the estimation of
change-points. We illustrate numerically our approach via simulations and a genomics dataset application. Future directions for this work are  the study of maximum likelihood estimation instead of least-squares, and a multivariate extension of the proposed algorithm.

\section{Proof of Theorems~\ref{thm1} and \ref{thm2}}
\label{sec:proofs-thm1-and-2}

Introduce $\mu = [\mu_{j}]_{1\leq j \leq m}\in \R^m$ given by
$\mu_{1} = \beta_{1}$ and $\mu_{j} = \beta_{j} -
\beta_{j-1}$ for $j=2, \ldots, m.$ Then, we have $\beta = \bT
\mu$, where $\bT$ is the $m \times m$ lower triangular matrix with
entries $(\bT)_{j, k} = 0$ if $j < k$ and $(\bT)_{j, k} = 1$
otherwise. Note that  $\hat \beta = \bT \hat \mu$,
where
\begin{equation}
  \label{sparvec}
  \hat \mu = \argmin_{\mu \in \R^m} \Big \{ \frac{1}{2} \norm{\bN -
    \bT \mu}_2^2 + \sum_{j = 2}^{m} \hat w_j |\mu_{j}| \Big\}.
\end{equation}

\subsection{Proof of Theorem \ref{thm1}}
\label{sec:proof-of-thm1}

This proof follows a standard argument for proving slow oracle
inequalities, see for instance~\cite{BicRitTsy-09}. 
Due to the Doob-Meyer decomposition theorem, we have
\begin{equation*}
  R_n(\lambda) = \norm{\lambda - \lambda_0}^2 - \norm{\lambda_0}^2 -
  \int_0^1 \lambda(t) d \bar M_n(t),
\end{equation*}
which leads  to
\begin{equation}
  \hat \lambda = \lambda_{\hat \beta} = \argmin_{\beta \in \R_+^m}
  \Big( \norm{\lambda_\beta -
    \lambda_0}^2 - 2 \int_0^1 \lambda_\beta(t) d \bar M_n(t) +
  \norm{\beta}_{\TV, \hat w} \Big).
\end{equation}
Then, using~\eqref{eq:hat-beta}, it implies that
\begin{equation}
  \label{ineqestm}
  \norm{\hat \lambda - \lambda_0}^2 \leq  \inf_\beta\norm{\lambda_{\beta} -
    \lambda_0}^2 + \frac{2}{n} \nu_n(\hat \lambda - \lambda_{\beta}) +
  \norm{\beta}_{\TV, \hat w} - \norm{\hat \beta}_{\TV, \hat w},
\end{equation}
where $\nu_n(\lambda) = \sum_{i=1}^n \int_0^1\lambda(t) dM_i(t)$ is a
centered empirical process. Note that
\begin{align}
  \label{empproc}
    \frac 1n \nu_n(\hat \lambda - \lambda_{\beta}) &= \sum_{j=1}^m
    (\hat{\beta}_{{j,m}} - \beta_{j,m}) \int_0^1
    \lambda_{j, m} (t) d \bar M_n(t) \nonumber\\
    &= \sum_{j=1}^m ( (\bT \hat \mu)_{j,m} - (\bT \mu)_{j,m})
    \int_0^1 \lambda_{j, m}(t) d \bar M_n(t)\nonumber \\
    &= \sum_{j=1}^m (\hat{\mu}_{j,m} - \mu_{j,m} ) \sum_{q = j}^m
    \int_0^1 \lambda_{q, m}(t) d \bar M_n(t).
\end{align}
Define the event $\Omega_{n}$ by
\begin{equation*}
  \Omega_{n} = \bigcap_{j=1}^m \Big\{ \Big|
  \sum_{q=j}^m \int_0^1 \lambda_{q,m} (t) d \bar M_n(t)\Big| \leq
  \frac{\hat w_j}{2} \Big\}.
\end{equation*}
The probabilistic  control of $\Omega_n$ is given in Proposition~\ref{prop1}
from Section~\ref{sec:proofs-thm1-and-2} below. It relies on a
slight modification of an empirical Bernstein inequality
from~\cite{GaiGui-12}, see also \cite{Rey-03}. On $\Omega_{n},$ we
have using~\eqref{empproc}
\begin{equation*}
  \frac 2n \nu_n(\hat \lambda - \lambda_{\beta}) \leq \sum_{j=1}^m
  \hat{w}_j \big| \hat{\mu}_{j,m} - \mu_{j,m}\big|,
\end{equation*}
Using~\eqref{ineqestm}, we obtain
\begin{align*}
  \norm{\hat \lambda - \lambda_0}^2 &\leq \norm{\lambda_{\beta} -
    \lambda_0}^2 + \sum_{j=1}^{m} \hat w_j | \hat \mu_{j,m}
  -\mu_{j,m} | + \sum_{j=1}^{m} \hat w_j ( |\mu_{j,m}| - |\hat
  \mu_{j,m}|) \\
  &\leq \norm{\lambda_{\beta} - \lambda_0}^2 + 2
  \sum_{j=1}^{m}\hat{w}_j | \mu_{j,m}|\\
&=  \norm{\lambda_{\beta} - \lambda_0}^2 + 2  \norm{\beta}_{\TV, \hat w}.
\end{align*}
 Then, on $\Omega_n$, \eqref{slowrate} in Theorem~\ref{thm1} holds true . It
remains now to control $\P(\Omega_n^\complement)$. We have, recalling
$\lambda_{j, m}(t) = \sqrt m \ind{(\frac{j-1}{m}, \frac{j}{m}]}(t)$,
that
\begin{equation*}
  \P[\Omega_n^\complement] \leq \sum_{j=1}^m \P \Big[ \Big|
  \sqrt{m} \int_0^1 \ind{(\frac{j-1}{m}, 1]} (t) d \bar M_n(t) \Big| >
  \frac{\hat w_j}{2} \Big],
\end{equation*}
so we need to control the tails of
\begin{equation*}
  U_j = \int_0^1 \ind{(\frac{j-1}{m}, 1]} (t) d \bar M_n(t),
\end{equation*}
which is the goal of the next proposition.
\begin{prop}
\label{prop1}
  For any  numerical constants $c_h >1$, $\varepsilon >0$ and $c_0 >0$ such that $ec_0 > 2(4/3 + \varepsilon)c_h,$ the following holds for any $z>0:$

\begin{equation*}
  \P \bigg[ |U_j| \geq c_{1,\varepsilon} \sqrt{\frac{z+ \hat h_{n, z,
        j}}{n} \hat V_j} + c_{3,\varepsilon} \frac{z + 1 + \hat
    h_{n,z,j}}{n} \bigg] \leq c e^{-z}
\end{equation*}
where

\begin{equation*}
 \hat{h}_{n,z,j}= c_h\log\log\Bigg(\frac{2en\hat{V}_j+ 2e(\frac{4}{3}+\varepsilon)z}{ec_0(z+1) - 2(\frac{4}{3}+\varepsilon )c_h}\vee e \Bigg),
\end{equation*}
 $ c_{1,\varepsilon}= 2\sqrt{1+\varepsilon},  c_{3,\varepsilon}= \sqrt{2\max\big(c_0, 2(1+\varepsilon)(\frac{4}{3}+\varepsilon)\big)} + \frac{1}{3},$ and $c = 6 + 4\big(\log(1+\varepsilon)\big)^{-c_h} \sum_{q\geq
  1}q^{-c_h}.$
\end{prop} 
The proof of Proposition~\ref{prop1} is given in
Appendix~\ref{app:proof-prop1}.
Choosing $z = x+ \log m $, it yields that
\begin{align*}
&\sum_{j=1}^m\P\Bigg[ |U_j|\geq c_{1,\varepsilon} \sqrt{\frac{x+\log m+\hat{h}_{n,x,j}}{n}\hat{V}_j} + c_{3,\varepsilon}\frac{x +\log m +\hat{h}_{n,x,j} +1}{n}\Bigg]\\
& \hspace*{1cm} \leq \big(6 + 4\big(\log(1+\varepsilon)\big)^{-c_h} \sum_{q\geq 1}q^{-c_h}\big)e^{-x},
\end{align*}
where
\begin{equation*}
\hat{h}_{n,x,j}  = c_h\log\log\Bigg(\frac{2en\hat{V}_j+ 2e(\frac{4}{3}+\varepsilon)(x+\log m)}{ec_0(x+\log m + 1) - 2(\frac{4}{3}+\varepsilon )c_h}\vee e \Bigg).
\end{equation*}
Then, the choice of data-driven weights is given by
\begin{equation*}
  \hat w_j = c_1 \sqrt{\frac{m(x + \log m + \hat h_{n, x, j}) \hat
      V_j}{n}} + c_2\frac{\sqrt{m}(x + 1 + \log m + \hat h_{n,x,j}) }{n},
\end{equation*}
where $c_1= 2c_{1,\varepsilon}$ and  $c_2= 2c_{3,\varepsilon}$  gives  $\P(\Omega_n^\complement) \leq c e^{-x}$. Finally, to get the numerical constants in  Theorem \ref{thm1},  we set $\varepsilon =1, c_h =2,$ and $ c_0 = 28/3e$  in Proposition \ref{prop1}.  
$\hfill \square$

\subsection{Proof of Corollary \ref{cor:slowrate}}

We denote by $\lambda_{0,
  m}$ the projection of $\lambda_0$ onto $\Lambda_m,$ that is $\lambda_{0,m} = \argmin_{\lambda_\beta \in \Lambda_m } \norm{\lambda_\beta -
   \lambda_0}^2.$ Using Pythagoras' theorem, we have 
\begin{equation*}
\norm{\hat{\lambda} - \lambda_0}^2 \leq \norm{  \lambda_{0,m} - \lambda_0}^2 + \norm{\hat{\lambda} - \lambda_{0,m}}^2.
\end{equation*}
By the proof of Theorem \ref{thm1}, we obtain 
\begin{eqnarray*}
\norm{\hat{\lambda} - \lambda_{0,m}}^2 &\leq& 2 \norm{\beta_{0,m}}_{\TV, \hat w}\\
&\leq& 2\norm{\beta_{0,m}}_{\TV} \max_{1\leq j \leq m} \hat{w}_j.
\end{eqnarray*}
Now, the following approximation lemma comes in handy for the control of
the bias term.

\begin{lem}
  \label{lem:approximation}
  Given Assumption~\ref{ass:intensity}, we have
  \begin{equation*}
    \norm{\lambda_{0, m} - \lambda_0}^2 \leq \frac{2 (L_0 - 1)
      \Delta_{\beta, \max}^2}{m},
  \end{equation*}
  where $\Delta_{\beta, \max} = \max\limits_{1 \leq \ell, \ell' \leq L_0}
  |\beta_{0,\ell} - \beta_{0,\ell'}|$.
\end{lem}

The proof of Lemma~\ref{lem:approximation} is given in
Appendix~\ref{app:proof-lem:approximation}.  $\hfill \square$

\subsection{Proof of Theorem~\ref{thm2}}

Using Pythagoras' identity, we obtain the following decomposition
\begin{equation*}
\norm{\lambda_{\hat \beta}- \lambda_{0}}^2  = \norm{\lambda_{ \beta}- \lambda_{0}}^2 + \norm{\lambda_{\hat \beta}- \lambda_{\beta}}^2. 
\end{equation*}
In view of the fact that $\{\lambda_{j,m} : j=1, \ldots, m \}$ is an orthonormal basis of $\Lambda_m$, we have 
\begin{equation*}
\norm{\lambda_{\hat \beta}- \lambda_{\beta}}^2 = \norm{\hat \beta -
  \beta}_2^2,
\end{equation*}
and by the definition of $\hat \beta,$ we get
\begin{align*}
  \norm{\hat \beta - \bN}_2^2 + \sum_{j=2}^{m} \hat w_j | \hat
  \beta_{j,m} - \hat \beta_{j-1,m}| \leq \norm{\beta - \bN}_2^2 +
  \sum_{j=2}^{m} \hat w_j | \beta_{j, m} - \beta_{j-1, m}|.
\end{align*}
Then  
\begin{align*}
&\big\| {\hat{\beta}} - \beta\big\|_{2}^2 \leq  \sum_{j=2}^{m}\hat{w}_{j} \Big(   | \beta_{j,m} - \beta_{j-1,m}| - | \hat{\beta}_{j,m} - \hat{\beta}_{j-1,m}|\Big) 
   + 2\int_0^1\sum_{j=2}^m (\hat{\beta}_{j,m} -  \beta_{j,m})\lambda_{j,m}(t)d\bar{M}_n(t).
\end{align*}
Assume that $ \hat{\beta}$  belongs to a set of dimension at most $L_{\max}$. Let
$ {S} = \big\{j:\, {\beta}_{j,m} \neq {\beta}_{j-1,m}
  \text{ for } j=2, \ldots, m \big\},$ be the support of the discrete gradient of $ \beta.$
 Using the Cauchy–Schwarz inequality, we have
\begin{equation*}
\begin{split}
&\sum_{j=2}^{m}\hat{w}_j \Big(| \beta_{j,m} - \beta_{j-1,m}| - |\hat{\beta}_{j,m} - \hat{\beta}_{j-1,m}| \Big) \\
& \qquad \leq \sum_{j\in \hat S\cup S} \hat{w}_j \Big(| \beta_{j,m} - \hat{\beta}_{j,m}| + |{\beta}_{j-1,m} - \hat{\beta}_{j-1,m}| \Big)\\
&\qquad \leq \sum_{j\in \hat S\cup S } \hat{w}_j \Big(|\beta_{j,m} - \hat{\beta}_{j,m}| \Big) +  \sum_{j\in \hat S\cup S} \hat{w}_j \Big(| \beta_{j-1,m} - \hat{\beta}_{j-1,m}| \Big)\\
& \qquad \leq \sum_{j\in \hat S\cup S \cup (\hat S\cup S+1 )} \hat{w}_j \Big( |\hat{\beta}_{j,m} - {\beta}_{j,m}| \Big)\\
& \qquad \leq  \sqrt{ \big|\hat S\cup S\cup (\hat S+1)\cup( S+1)\big| } \\
&\hspace*{3.5cm} \times \bigg\|\bigg[ {\hat{\beta}}_{j,m} - \beta_{j,m}\bigg]_{j \in \hat S\cup S \cup (\hat S+1)\cup( S+1)}\bigg\|_{2} 
 \times \max_{j\in \hat S\cup S\cup (\hat S+1)\cup( S +1)}\hat{w}_j\\
& \qquad \leq  \sqrt{2} \sqrt{L_{\max}+2(L_0-1)}\,\big\| {\hat{\beta}} - \beta\big\|_{2} \max_{j=1, \ldots,m}\hat{w}_j.
\end{split}
\end{equation*}
Hence
\begin{flalign*}
\nonumber
&\big\| {\hat{\beta}} - \beta\big\|_{2}^2 \leq\sqrt{2}\sqrt{L_{\max}+2(L_0-1)}\,\big\| {\hat{\beta}} - \beta \big\|_{2} \max_{j=1, \ldots, m}\hat{w}_j \\
&\qquad  \qquad \qquad \qquad \qquad + 2 \big\| {\hat{\beta}} - \beta\big\|_{2}\int_0^1\sum_{j=2}^m \frac{(\hat{\beta}_{j,m} -  \beta_{j,m})\lambda_{j,m}(t)}{\big\| {\hat{\beta}} - \beta\big\|_{2}}d\bar{M}_n(t).
\end{flalign*}
Now,  define the functional $G$  for all $\lambda_\beta \in \Lambda_m$ in the following way:
\begin{equation*}
G(\lambda_\beta) =\int_0^1\frac{\lambda_\beta (t)}{\norm{\lambda_\beta}}d\bar{M}_n(t).
\end{equation*}
Therefore, we obtain
\begin{align*}
&\big\| {\hat{\beta}} - \beta\big\|_{2}^2 \leq  \sqrt{2}\sqrt{L_{max}+2(L_0-1)}\,\big\| {\hat{\beta}} - \beta\big\|_{2} \max_{j=1, \ldots,m}\hat{w}_j  + 2\big\| {\hat{\beta}} - \beta\big\|_{2} G( \hat{\beta} - \beta).
\end{align*}
Let 
\begin{equation*}
\mathcal{V} = \bigcup_{L=1}^{L_{\max}} V_L = \bigcup_{L=1}^{L_{\max}} \bigcup_{J \subset \{1, \ldots, m-1\},\, |J| = L} V_{L,J},
\end{equation*}
where 
$\big\{{V}_L: L=1, \ldots, L_{\max}\big\}$ is the collection of the spaces to which $ {\hat{\beta}}$  may belong and $V_{L,J}$ denotes a space of dimension $L$ containing  signals with a  support $J$. \\

\noindent It follows that, 
\begin{equation}
\label{difprjest}
\big\| {\hat{\beta}} - \beta\big\|_{2} \leq \sqrt{2}\sqrt{L_{\max}+2(L_0-1)} \max_{j=1, \ldots, m}\hat{w}_j+ 2 \sup_{\lambda \in \mathcal{V}, \norm{\lambda}=1} G(\lambda).
\end{equation}
Then by Proposition 4 in  \hypertarget{}{\cite{ComGaiGui-11}},
we have for any $z>0$ 
\begin{equation*} 
\P\Bigg[\sup_{\lambda\in {V}_{L,J}, \,\|\lambda\| = 1}G(\lambda) \geq \kappa\bigg(\sqrt{\frac{\|\lambda_0\|_{\infty} (L+z)}{n}} + \frac{2\sqrt{m}(L+z)}{\sqrt{L}n} \bigg)\Bigg]\leq e^{-z}, 
\end{equation*}
where $\kappa= 11.8.$ Then
\begin{align*}
&\sum_{ \substack{L=1,\ldots, L_{\max}\\ {J \subset \{1, \ldots, m-1\},\,  |J| = L}}} \P\Bigg[\sup_{\lambda\in {V}_{L,J},\, \|\lambda\| = 1}G(\lambda)
 \geq \kappa\bigg(\sqrt{\frac{\|\lambda_0\|_{\infty} (L+z)}{n}} + \frac{2\sqrt{m}(L+z)}{\sqrt{L}n} \bigg)\Bigg] \\
&\qquad  \leq  \sum_{ \substack{L=1,\ldots, L_{\max}\\ {J \subset \{1, \ldots, m-1\},\,  |J| = L}}}e^{-z}\\
&\qquad \leq {L_{\max}} m^{L_{\max}}e^{-z}.
\end{align*}
Choosing $ z =  x + L_{\max}\log m$ for $x >0,$ leads to  
 \begin{equation*}
\begin{split}
&\sum_{ \substack{L=1,\ldots, L_{\max}\\ {J \subset \{1, \ldots, m-1\}, \, |J| = L}}} \P\Bigg[\sup_{\lambda\in {V}_{L,J},\, \|\lambda\| = 1}G(\lambda) \geq \kappa\bigg(\sqrt{\frac{\|\lambda_0\|_{\infty} (L+ x + L_{\max}\log m)}{n}} \\
&\hspace*{8cm} + \frac{2\sqrt{m}(L+ x + L_{\max}\log m)}{\sqrt{L}n} \bigg)\Bigg]\\
&\qquad\leq {L_{\max}} e^{-x}.
\end{split}
\end{equation*}
Plugging this in  inequality ($\ref{difprjest}$), we obtain for any $x>0$  and with probability larger than $1 - L_{\max}e^{-x}$ 
\begin{align*}
\big\| {\hat{\beta}} - \beta\big\|_{2} &\leq  \sqrt{2}\sqrt{L_{\max}+ 2(L_0-1)} \max_{j=1, \ldots,m}\hat{w}_j \\
&\qquad + {2}\kappa\sqrt{\frac{\|\lambda_0\|_{\infty} (x +L_{\max} (1+\log m))}{n}} \\
&\qquad  +4\kappa\frac{\sqrt{m}(x + L_{\max} (1+\log m ))}{n},
\end{align*}
and the result follows by using the inequality $(a+b+c)^2 \leq 3(a^2 + b^2 +c^2),$ for all $a, b, c \in \R.$ 
$\hfill \square$

\section{Proof of Theorem \ref{thm3}}
\label{sec:proof-thm3}

Let us give first the overall structure of the proof, which is inspired from~\cite{HarLev-10}. 
In this proof, we repeatedly use the KKT optimality conditions of the optimization problem~\eqref{sparvec}, given by Lemma~\ref{lem:KKT} below.
We use also repeatedly deviation arguments of the data-driven weights $\hat w_j$ and a control of the martingale noise, which are provided by Lemma~\ref{lem:control-martingale} below.
We prove consistency of $\hat \tau_\ell =\frac{\hat j_\ell}{m}$, which is an estimator of the right-hand side boundary $\frac{j_\ell}{m}$ of the interval $I_{j_\ell, m} = (\frac{j_\ell -1}{m}, \frac{j_\ell}{m}],$  by showing that $\P[A_{n,\ell}] \rightarrow 0$ as $n \rightarrow \infty$, where $A_{n,\ell} := \big\{|\hat{j}_\ell - j_\ell| > \frac{m\varepsilon_n}{2}\big\}$, for all $\ell \in \{1, \ldots, L_0 -1\}.$
We treat separately two cases depending on  the positions of $j_\ell$ and $\hat j_\ell,$. In Case~I, we consider $\hat j_\ell < j_\ell$, see Section~\ref{CASE-I} and Figure~\ref{fig:CASE-I}. In Case~II., we consider $\hat j_\ell > j_\ell,$, see Appendix~\ref{app:proof-thm3-case2} and Figure~\ref{fig:CASE-II}.
We decompose even further, using the quantity $\Delta_{j,\min}$ (see Section~\ref{section:change-point-agg}), defining the set $C_n = \big\{ \max_{1\leq \ell\leq L_0 -1}|\hat{j}_\ell - j_\ell| < \frac{\Delta_{j,\min}}{2}\big\}.$ 
We prove that $\P[A_{n,\ell} \cap C_n] \rightarrow 0$ and $\P[A_{n,\ell} \cap C_n ^\complement]\rightarrow 0$ as $ n\rightarrow \infty$ for Case~I in Sections~\ref{subsection1:CASE-I},~\ref{subsection2:CASE-I}, and for Case~II in Appendices~\ref{subsection1:CASE-II},~\ref{subsection2:CASE-II}. \\




\begin{minipage}[b]{0.45\textwidth}
\centering
\begin{tikzpicture}[scale = 1.2]
\tkzInit[xmin=0.02, xmax=0.48, xstep= 0.1]
\tkzDrawX[label=$ t$, right = 8 pt,noticks]
\let\tkzmathstyle\displaystyle
\tkzHTicks[mark=ball, ball color= gray, mark options={mark size=1.1pt, color=gray}]{0.05,0.10,0.15,0.20,0.25,0.30,0.35,0.40,0.45,0.50}
\tikzset{point style/.style={ draw=red,inner sep=0pt,shape=circle,minimum size=4pt,fill = red}}
\tkzHTick[mark=ball, ball color = black, mark size=1.5pt]{0.13}
\tkzDefPoint(0.13,0){tau2}
\tkzLabelPoint[above=8pt](tau2){${\tau_{0,{\ell-1}}}$}
\tkzHTick[mark=ball,mark size=1.5pt,ball color = black]{0.33}
\tkzDefPoint(0.33,0){tau3}
\tkzLabelPoint[above=8pt](tau3){${\tau_{0,\ell}}$}
\tkzHTick[mark= ball,mark size=1.5pt, solid, ball color = black]{0.48}
\tkzDefPoint(0.48,0){tau5}
\tkzLabelPoint[above=8pt](tau5){${\tau_{0,\ell+1}}$}
\draw[] (1.3,0.20) -- (1.3,-0.6);
\draw[] (0.8,0.20) -- (0.8,-0.6);
\draw[,<->](0.8,-0.5)--(1.3,-0.5);
\draw[,->](1.05,-0.8)--(1.05,-0.5);
\tkzDefPoint(0.13,-0.7){labelIj2}
\tkzLabelPoint[below](labelIj2){${I_{j_{\ell-1},m}}$}
\draw[] (2.8,0.20) -- (2.8,-0.6);
\draw[] (3.3,0.20) -- (3.3,-0.6);
\draw[,<->](3.3,-0.5)--(2.8,-0.5);
\draw[,->](3.05,-0.8)--(3.05,-0.5);
\tkzDefPoint(0.33,-0.7){labelIj2}
\tkzLabelPoint[below](labelIj2){${I_{j_{\ell},m}}$}
\draw[] (4.3,0.20) -- (4.3,-0.6);
\draw[] (4.8,0.20) -- (4.8,-0.6);
\draw[,<->](4.3,-0.5)--(4.8,-0.5);
\draw[,->](4.55,-0.8)--(4.55,-0.5);
\tkzDefPoint(0.48,-0.7){labelIj2}
\tkzLabelPoint[below](labelIj2){${I_{j_{\ell +1},m}}$}
\tkzHTick[mark = ball, ball color = gray,  mark size = 1.5pt, ]{0.25}
\draw[gray]((2.3,-0.6)node[above]{{$ \hat{\tau}_\ell$}};
\end{tikzpicture}
\captionsetup{font=small,labelformat=default} 
\captionof{figure}{Case I. $\hat j_\ell < j_\ell$} \label{fig:CASE-I} 
\end{minipage}
\begin{minipage}[b]{0.5\textwidth}
\centering
\begin{tikzpicture}[scale = 1.2]
\tkzInit[xmin=0.02, xmax=0.48, xstep= 0.1]
\tkzDrawX[label=$ t$, right = 8 pt, ,noticks]
\let\tkzmathstyle\displaystyle
\tkzHTicks[mark=ball, ball color= gray, mark options={color=gray,mark size=1.1pt}]{0.05,0.10,0.15,0.20,0.25,0.30,0.35,0.40,0.45,0.50}
\tikzset{point style/.style={ draw=red,inner sep=0pt,shape=circle,minimum size=4pt,fill = red}}
\tkzHTick[mark=ball, ball color = black,mark size=1.5pt]{0.13}
\tkzDefPoint(0.13,0){tau2}
\tkzLabelPoint[above=8pt,](tau2){${\tau_{0,{\ell-1}}}$}
\tkzHTick[mark=ball, ball color = black, mark size=1.5pt]{0.23}
\tkzDefPoint(0.23,0){tau3}
\tkzLabelPoint[above=8pt,](tau3){${\tau_{0,\ell}}$}
\tkzHTick[mark=ball, ball color = black, mark size=1.5pt]{0.48}
\tkzDefPoint(0.48,0){tau5}
\tkzLabelPoint[above=8pt,](tau5){${\tau_{0,\ell+1}}$}
\draw[] (1.3,0.20) -- (1.3,-0.6);
\draw[] (0.8,0.20) -- (0.8,-0.6);
\draw[,<->](0.8,-0.5)--(1.3,-0.5);
\draw[,->](1.05,-0.8)--(1.05,-0.5);
\tkzDefPoint(0.13,-0.7){labelIj2}
\tkzLabelPoint[below,](labelIj2){${I_{j_{\ell-1},m}}$}
\draw[] (1.8,0.20) -- (1.8,-0.6);
\draw[] (2.3,0.20) -- (2.3,-0.6);
\draw[,<->](1.8,-0.5)--(2.3,-0.5);
\draw[,->](2.03,-0.8)--(2.03,-0.5);
\tkzDefPoint(0.23,-0.7){labelIj2}
\tkzLabelPoint[below,](labelIj2){${I_{j_{\ell},m}}$}
\draw[] (4.3,0.20) -- (4.3,-0.6);
\draw[] (4.8,0.20) -- (4.8,-0.6);
\draw[,<->](4.3,-0.5)--(4.8,-0.5);
\draw[,->](4.55,-0.8)--(4.55,-0.5);
\tkzDefPoint(0.48,-0.7){labelIj2}
\tkzLabelPoint[below,](labelIj2){${I_{j_{\ell +1},m}}$}
\tkzHTick[mark = ball, mark size = 1.5pt, ball color  = gray]{0.35}
\draw [gray]((3.4,-0.6)node[above]{\textcolor{gray}{$ \hat{\tau}_\ell$
 }};
\end{tikzpicture}
\captionsetup{font=footnotesize,labelformat=default}
\captionof{figure}{Case II. $\hat j_\ell > j_\ell$ \label{fig:CASE-II}}
 \end{minipage}
\\
 \begin{lem}
\label{lem:KKT}
Consider   the total-variation penalized problems in \eqref{estvec2} and\eqref{sparvec}. Let $\hat{ {\beta}} = [\hat{\beta}_{j,m}]_{1\leq j \leq m}$  and $\hat{ {\mu}} = [\hat{\mu}_{j,m}]_{1\leq j \leq m}$  denote the respective solutions. Then, the latter vectors and the approximate change-points sequence estimators $\hat{j}_1, \ldots, \hat{j}_{|\hat{S}|}$ satisfy for all  $r = 1, \ldots, |\hat{S}|,$
\begin{equation}
 \sum_{j= \hat{j}_r}^m \beta_{0,j,m} - \sum_{j= \hat{j}_r}^m\hat{\beta}_{j,m} +\sqrt{m} \sum_{j= \hat{j}_r}^m\bar M_n(I_{j,m}) = {\hat{w}_{\hat{j}_r}}\sgn(\hat{\mu}_{\hat{j}_r,m}),
\end{equation}
and \textrm{ for all }  $j \in \{1, \ldots, m\},$
\begin{equation}
\label{kkt}
 \bigg|\sum_{q=j}^m \beta_{0,q,m} - \sum_{q= j}^m\hat{\beta}_{q,m} + \sqrt{m} \sum_{q= j}^m\bar M_n(I_{q,m})\bigg| \leq {\hat{w}_j},
\end{equation}
using the convention $\sgn(\hat{\mu}_{\hat{j}_r,m}) = +1,$ if $ \hat{\mu}_{\hat{j}_r,m} >0$ and $-1$ otherwise.  The vectors $\hat{ {\beta}}$  and $\beta_{0,m} = [\beta_{0,j,m}]_{1\leq j \leq m}$  have the following additional properties
\begin{equation}
\label{prpvecest}
\left\{
  \begin{array}{ll}
    \hat{\beta}_{q,m}  = \hat{\beta}_{\hat{j}_r - 1,m}, & \mbox{if } \, \hat{j}_{r-1} + 1 \leq q\leq \hat{j}_r, \textrm{ for  } r=1,\ldots, \hat L, \\
   \beta_{0,q,m} = \beta_{0,j_\ell - 1, m},  & \mbox{if } j_{\ell -1}+1 \leq q \leq j_\ell - 1, \textrm{ for } \ell=1, \ldots, L_0-1.
  \end{array}
\right.
\end{equation}
\end{lem}
The proof of Lemma~\ref{lem:KKT} is given in
Appendix~\ref{app:proof-lem:KKT}.  Let
us now state a lemma which allows us to control the martingale noise
term.

\begin{lem}
\label{lem:control-martingale}
Given two integers $a$ and $b,$  such that  $1 \leq a < b\leq m$, let $\bar M_n(a;b):= \sum_{q=a}^{b} \bar M_n(I_{q,m})$.Then, for all $z>0$ we have 
 \begin{equation}
\label{ctrlmart}
\begin{split}
&\P\Big[\big|\bar M_n(a;b)\big| \geq z \Big] \leq 2\exp\Bigg(-\frac{nz^2}{2\int_{\ind{(\frac{a-1}{m},\frac{b}{m}]}} \lambda_0(t)dt
+\frac{2}{3}z }\Bigg),
\end{split}
\end{equation}
and for all $\xi >0$, the data driven weight $\hat w_a$ satisfies

\begin{equation}
\label{ctrlweigh}
\begin{split}
&\P\Bigg[\hat{w}_a^2 \geq \frac{m\log m }{n}\bigg(\xi- \int_{\ind{(\frac{a-1}{m},1]}} \lambda_0(t)dt
\bigg) \Bigg]\leq2\exp\Bigg(-\frac{n\xi^2}{2\int_{\ind{(\frac{a-1}{m},1]}} \lambda_0(t)dt
+\frac{2}{3}\xi }\Bigg), 
\end{split}
\end{equation}
where $\int_{I} \lambda_0(t) dt = \E[\bar N (I)]$  for any $I \subset [0, 1]$.
\end{lem}

The proof of Lemma~\ref{lem:control-martingale} is given in
Appendix~\ref{app:proof-lem-control-martingale}. Let us now prove Theorem~\ref{thm3}. Recall that the
sequence $(\varepsilon_n)_n$ satisfies $m\varepsilon_n \geq 6,$ for
all $n\geq1$ . An application of the triangle inequality entails that,
\begin{equation*}
  \begin{split}
      & \P\Big[\max_{1\leq \ell\leq L_0 - 1}|\tau_{0,\ell} - \hat{\tau}_\ell| >\varepsilon_n\Big]
     \leq \P\Big[\max_{1\leq \ell\leq L_0 - 1}|\tau_{0,\ell} - \frac{j_\ell}{m}| > \frac{\varepsilon_n}{2}\Big] +  \P\Big[\max_{1\leq \ell\leq L_0 - 1}|\frac{j_\ell}{m} - \hat{\tau}_\ell| >      \frac{\varepsilon_n}{2}\Big].
\end{split}
\end{equation*} 
Moreover, the true change-point $\tau_{0,\ell}$ verifies \eqref{positon-of-true-change} which implies that
\begin{equation*}
  \begin{split}
      & \P\Big[\max_{1\leq \ell\leq L_0 - 1}|\tau_{0,\ell} - \hat{\tau}_\ell| >\varepsilon_n\Big] \leq \P\Big[\max_{1\leq \ell\leq L_0 - 1}|{j_\ell} - \hat{j}_\ell| >\frac{m\varepsilon_n}{2} \Big].
\end{split}
\end{equation*} 
Due to  $$ \P\Big[\max_{1\leq \ell\leq L_0-1}|\hat{j}_\ell - j_\ell| >\frac{m\varepsilon_n}{2}\Big] \leq \sum_{\ell=1}^{L_0-1} \P\big[|\hat{j}_\ell - j_\ell| > \frac{m\varepsilon_n}{2}\big],$$ it suffices to prove that for all $\ell = 1, \ldots, L_0-1$,   $\P[A_{n,\ell}] \rightarrow 0$, as $n$ tending to infinity.

\subsection{Case I}
\label{CASE-I}

Due to the fact that  $m\varepsilon_n \geq 6$ for all $n \geq 1$, it follows  that the event  $ \big\{ \hat{j}_\ell < j_{\ell} -2 \big\}$ a.s.

\subsubsection{Step I.1.  Prove: $ \P[A_{n,\ell} \cap C_n] \rightarrow 0,$ as $n\rightarrow \infty.$ }
\label{subsection1:CASE-I}

By the definition of $C_n$, we have 
\begin{equation}
\label{lochatjl}
 j_{\ell-1} < \hat{j}_\ell < j_{\ell+1}, \, \,\textrm{ for all }\, \, \ell=1, \ldots, L_0-1.
 \end{equation}
Applying $(\ref{kkt})$ in Lemma~\ref{lem:KKT}  with $j = j_\ell$  and $j = \hat{j}_\ell +1$,  we obtain
 \begin{equation*}
-(\hat{w}_{j_\ell}+\hat{w}_{\hat{j}_\ell +1})\leq \sum_{q= \hat{j}_\ell +1}^{j_\ell-1}\bN_q  - \sum_{q= \hat{j}_\ell +1}^{j_\ell -1}\hat{\beta}_{q,m} \leq\hat{w}_{j_\ell}+\hat{w}_{\hat{j}_\ell +1}.
\end{equation*}
Put $\hat{w}_{a,b}:= {\hat{w}_a + \hat{w}_b},$ for any two integers $a$ and $b.$ Thus
\begin{equation*}
\Big| \sum_{q=\hat{j}_\ell +1}^{j_\ell -1}  {\beta_{0,q,m}} -  \hat{\beta}_{q,m}+ \sqrt{m}\bar M_n(I_{q,m})\Big| \leq  \hat{w}_{\hat{j}_\ell +1,j_\ell}.
\end{equation*}
Using the  property of the vector $\hat{ {\beta}}$ in Lemma~\ref{lem:KKT}, we get

\begin{equation*}
 \Big|(j_\ell - \hat{j}_\ell - 2)\big( \beta_{0,j_\ell-1,m} - \hat{\beta}_{\hat{j}_{\ell+1} - 1,m} ) + \sqrt{m}\bar M_n(\hat{j}_\ell + 1;j_\ell-1) \Big|\leq \hat{w}_{\hat{j}_\ell +1,j_\ell}.
\end{equation*}
Therefore, on $C_n \cap \{\hat{j}_\ell  < j_\ell -2\},$ we have
\begin{equation*}
  \begin{split}
    & \Big|( \hat{j}_\ell - j_\ell - 2) ( \hat{\beta}_{\hat{j}_{\ell+1} - 1,m} -     \beta_{0,j_{\ell+1}-1,m}) \\
    &\hspace*{3cm} + (\hat{j}_\ell - j_\ell - 2 )( \beta_{0,j_{\ell+1}-1,m}-   \beta_{0,j_\ell-1,m}) \\
      &\hspace*{5cm}+ \sqrt{m}\bar M_n(\hat{j}_\ell+1;j_\ell-1)\Big|\leq\hat{w}_{\hat{j}_\ell +1,j_\ell}.
    \end{split}
\end{equation*}
\noindent Defining the event
\begin{equation*}
    \begin{split}
       &C_{n,\ell}= \bigg\{\Big|( \hat{j}_\ell - j_\ell - 2) ( \hat{\beta}_{\hat{j}_{\ell+1} - 1,m } -  \beta_{0,j_{\ell+1}-1,m})\\
     &\hspace*{3cm} + (\hat{j}_\ell - j_\ell - 2 )( \beta_{0,j_{\ell+1}-1,m}- \beta_{0,j_\ell-1,m})\\
     &\hspace*{5cm}  + \sqrt{m}\bar M_n(\hat{j}_\ell+1;j_\ell-1)+ \Big|\leq \hat{w}_{\hat{j}_\ell +1,j_\ell}\bigg\},
     \end{split}
\end{equation*}
We observe that $ C_{n,\ell}$ occurs with probability one. In addition,
we remark that for all $n\geq 1$, $m\varepsilon_n \geq 6$ entails $\frac{m\varepsilon_n}{2} - 2 \geq \frac{m\varepsilon_n}{6}$. Then
\begin{equation*}
\Big\{ |\hat{j}_\ell - j_\ell| >\frac{m\varepsilon_n}{2} \Big\} \subset \Big\{ |\hat{j}_\ell - j_\ell - 2| >\frac{m\varepsilon_n}{2} - 2 \Big\} \subset \Big\{ |\hat{j}_\ell - j_\ell - 2| \geq \frac{m\varepsilon_n}{6} \Big\}
\end{equation*}
\noindent Therefore
\begin{equation*}
    \begin{split}
        &\P[A_{n,\ell}\cap C_n \cap C_{n,\ell}]\\
       &\qquad\qquad\leq  \P\Bigg[\bigg\{\frac{\hat{w}_{\hat{j}_\ell +1,j_\ell}}{| \hat{j}_\ell -j_\ell - 2| } \geq \frac{| \beta_{0,j_{\ell+1}-1,m} -  \beta_{0,j_\ell-1,m}|}{3}\bigg\} \cap \Big\{\hat{j}_\ell < j_\ell - 2\Big\}  \Bigg] \\
       &\qquad\qquad\quad  + \P\Bigg[\bigg\{ |\hat{\beta}_{\hat{j}_{\ell+1} - 1,m}-\beta_{0,j_{\ell+1}-1,m}| \geq \frac{|\beta_{0,j_{\ell+1}-1,m} -  \beta_{0,j_\ell-1,m}|}{3}\bigg\} \cap C_n \Bigg]\\
        &  \qquad\qquad \quad+ \P\Bigg[\bigg\{ \bigg| \frac{\sqrt{m}\bar M_n(\hat{j}_\ell + 1 ;j_\ell-1)}{ \hat{j}_\ell -j_\ell - 2 }\Bigg|\geq \frac{|\beta_{0,j_{\ell+1}-1,m} -  \beta_{0,j_\ell-1,m}|}{3}\bigg\}  \Bigg]\\
       & \qquad\qquad := \P[A_{n,\ell,1}] + \P[A_{n,\ell,2}] + \P[A_{n,\ell,3}].
      \end{split}
\end{equation*}
Moreover, we have 
\begin{equation*}
     \begin{split}
       \P[A_{n,\ell,1}] &\leq \P\Big[\hat{w}_{\hat{j}_\ell +1,j_\ell} \geq     \frac{m\varepsilon_n \Delta_{\beta,\min}}{18}\Big]\\
       &\leq \P\Big[\hat{w}_{\hat{j}_\ell +1} \geq \frac{m\varepsilon_n \Delta_{\beta,\min}}{36}\Big]\\
       &\leq \P\Big[\hat{w}^2_{{j}_{\ell-1} +1} \geq \frac{m^2\varepsilon_n^2 \Delta_{\beta,\min}^2}{36^2}\Big].
      \end{split}
\end{equation*}
 By (\ref{ass:consistency-1})  in Assumption~\ref{ass:consistency},  and \eqref{ctrlweigh} in Lemma~\ref{lem:control-martingale} with $\xi = \frac{nm\varepsilon_n^2\Delta_{\beta,\min}^2}{36^2\log m} + \E\big[\bar N_n\big((\frac{j_{\ell-1}}{m},1]\big)\big] ,$
it follows that 
\begin{eqnarray*}
\P[A_{n,\ell,1}] &\leq&  2\exp\bigg(-\frac{n\xi^2}{2 \E\Big[\bar N_n\Big(\big(\frac{j_{\ell-1}}{m},1]\big)\Big) \Big]+\frac{2}{3}\xi }\bigg) \rightarrow 0,
\end{eqnarray*}
as $n \rightarrow \infty.$ Next, consider the event
\begin{eqnarray*}
 A_{n,\ell,3}&=&  \left\{\Bigg| \frac{\sqrt{m}\bar M_n(\hat{j}_\ell+1;j_\ell-1)}{ \hat{j}_\ell - j_\ell - 2}\Bigg|\geq \frac{|\beta_{0,j_{\ell+1}-1,m} -  \beta_{0,j_\ell-1,m}|}{3}\right\}\\
& =& \left\{ \left| \bar M_n(\hat{j}_\ell+1;j_\ell-1)\right|\geq \left| \hat{j}_\ell - j_\ell - 2\right|\frac{\left| \beta_{0,j_{\ell+1}-1,m} -  \beta_{0,j_\ell-1,m}\right|}{3\sqrt{m}}\right\}\\
&\subset&\left\{\left| \bar M_n(\hat{j}_\ell+1;j_\ell-1)\right|\geq \frac{m\varepsilon_n \Delta_{\beta,\min}}{18\sqrt{m}}\right\} \bigcap \bigcup_{q=j_{\ell-1}+1}^{j_\ell -3} \left\{ \hat{j}_\ell = q\right\}\\
&\subset&  \bigcup_{q=j_{\ell-1}+2}^{j_\ell -2}  \left\{\left| \bar M_n(q;j_\ell-1)\right|\geq \frac{ m\varepsilon_n \Delta_{\beta,\min}}{18\sqrt{m}}\right\}.
\end{eqnarray*}
Put $\varphi_n =  \frac{ \sqrt{m}\varepsilon_n \Delta_{\beta,\min}}{18}.$ By \eqref{ctrlmart} in Lemma~\ref{lem:control-martingale}, we have
\begin{eqnarray*}
\P[ A_{n,\ell,3}] &\leq& 2 \sum_{q=j_{\ell-1}+2}^{j_\ell -2} \exp\Bigg(-\frac{n\varphi_n^2}{2 \E\Big[\bar N_n\Big(\big(\frac{q-1}{m},\frac{j_\ell-1}{m}\big]\Big)\bigg] +\frac{2}{3}\varphi_n }\Bigg)\\
&\leq& 2 (j_\ell - j_{\ell-1}-3)\exp\Bigg(-\frac{n\varphi_n^2}{2 \E\Big[\bar N_n\Big(\big(\frac{j_{\ell-1}+1}{m},\frac{j_\ell-1}{m}\big]\Big)\Big]+\frac{2}{3}\varphi_n }\Bigg)\\
&\leq& 2 \exp\Bigg(-\frac{n\varphi_n^2}{2 \E\Big[\bar N_n\Big(\big(\frac{j_{\ell-1}+1}{m},\frac{j_\ell-1}{m}\big]\Big)\Big]+\frac{2}{3}\varphi_n  } + \log m\Bigg).
\end{eqnarray*}
By  (\ref{ass:consistency-1})  in Assumption~\ref{ass:consistency} ,  it implies that $\P[A_{n,\ell,3}]$ goes to zero as $n \rightarrow \infty.$
 We now control $\P[A_{n,\ell,2}]$. Using Lemma~\ref{lem:KKT} with $j = \lceil\frac{j_\ell + j_{\ell+1}}{2} \rceil$ and with $j = j_\ell+1$, and using the triangle inequality, it follows that
\begin{equation*}
\Bigg| \sum_{q=j_\ell+1}^{\lceil\frac{j_\ell + j_{\ell+1}}{2} \rceil -1}\bN_q  - \sum_{q=j_\ell+1}^{\lceil\frac{j_\ell + j_{\ell+1}}{2} \rceil -1}\hat{\beta}_{q,m}\Bigg| \leq \hat{w}_{j_\ell+1,\lceil\frac{j_\ell + j_{\ell+1}}{2} \rceil}.
\end{equation*}
Furthermore, on the event $C_n \cap \{\hat{j}_\ell < j_\ell -2\},$  the following inequalities 
\begin{equation*}
\hat{j}_\ell < j_\ell \leq q \leq \lceil\frac{j_\ell + j_{\ell+1}}{2} \rceil - 1  \leq j_{\ell+1} -1,
\end{equation*}
hold true. Moreover, we note that
$\hat{\beta}_{q,m} = \hat{\beta}_{\hat{j}_{\ell+1}-1,m}$ if  ${j_\ell}\leq q \leq \lceil\frac{j_\ell + j_{\ell+1}}{2} \rceil -1  \leq \hat{j}_{\ell+1} -1.$ Consequently, we have 
\begin{equation*}
\begin{split}
 &\Big| (j_{\ell+1} - j_\ell -2)\frac{(\beta_{0,j_{\ell+1}-1,m} - \hat{\beta}_{\hat{j}_{\ell+1}-1,m})}{2}
 + \sqrt{m}\bar M_n({j}_\ell+1;\lceil\frac{j_\ell + j_{\ell+1}}{2} \rceil -1) \Big| \leq \hat{w}_{j_\ell+1,\lceil\frac{j_\ell + j_{\ell+1}}{2} \rceil},
\end{split}
\end{equation*}
which implies that
\begin{equation*}
\begin{split}
&(j_{\ell+1} - j_\ell-2)\frac{|\hat{\beta}_{\hat{j}_{\ell+1}-1,m} - \beta_{0,j_{\ell+1}-1,m}|}{2}
 \leq \hat{w}_{j_\ell+1,\lceil\frac{j_\ell + j_{\ell+1}}{2} \rceil} +\Big|\sqrt{m}\bar M_n({j}_\ell+1;\lceil\frac{j_\ell + j_{\ell+1}}{2} \rceil -1) \Big|.
\end{split}
\end{equation*}
Therefore, we may upper bound $\P[A_{n,\ell,2}]$ as follows
\begin{equation*}
\begin{split}
  &\P[A_{n,\ell,2}]\\
& = \P\Big[\Big\{ |\hat{\beta}_{\hat{j}_{\ell+1}-1,m}-\beta_{0,j_{\ell+1}-1,m}| \geq \frac{| \beta_{0,j_{\ell+1}-1,m} - \beta_{0,j_\ell-1,m}|}{3}\Big\}\cap C_n\Big] \\
  &= \P\Big[\Big\{(j_{\ell+1} - j_\ell -2) \frac{|\hat{\beta}_{\hat{j}_{\ell+1},m}-\beta_{0,j_{\ell+1}-1,m}|}{2} \\
&\hspace*{5cm} \geq(j_{\ell+1} - j_\ell-2) \frac{| \beta_{0,j_{\ell+1}-1,m} - \beta_{0,j_\ell-1,m}|}{6} \Big\}\cap C_n\Big]\\
   & \leq \P\Big[\Big \{\hat{w}_{j_\ell+1,\lceil\frac{j_\ell + j_{\ell+1}}{2} \rceil} + \big|\sqrt{m}\bar M_n({j}_\ell+1;\lceil\frac{j_\ell + j_{\ell+1}}{2} \rceil -1)\big|  \\
&\hspace*{5cm} \geq(j_{\ell+1} - j_\ell-2) \frac{| \beta_{0,j_{\ell+1}-1,m} - \beta_{0,j_\ell-1,m}|}{6} \Big\}\cap C_n\Big]\\
   &\leq \P\Big[\hat{w}_{j_\ell+1,\lceil\frac{j_\ell + j_{\ell+1}}{2} \rceil} \geq (j_{\ell+1} - j_\ell -2)\frac{| \beta_{0,j_{\ell+1}-1,m} - \beta_{0,j_\ell-1,m}|}{12}\Big] \\
   &   \quad  + \P\Big[\Big|\sqrt{m}\bar M_n({j}_\ell+1;\lceil\frac{j_\ell + j_{\ell+1}}{2} \rceil -1) \Big|\geq (j_{\ell+1} - j_\ell -2)\frac{| \beta_{0,j_{\ell+1}-1,m} - \beta_{0,j_\ell-1,m}|}{12}\Big]\\
 & \leq \P\Big[\hat{w}_{j_\ell+1,\lceil\frac{j_\ell + j_{\ell+1}}{2} \rceil} \geq \frac{(\Delta_{j,\min}-2)\Delta_{\beta,\min}}{12}\Big]\\
 &   \quad   + \P\Big[\Big|\bar M_n({j}_\ell+1;\lceil\frac{j_\ell + j_{\ell+1}}{2} \rceil -1) \Big|\geq \frac{(\Delta_{j,\min}-2) \Delta_{\beta,\min}}{12\, \sqrt{m}}\Big].
\end{split}
\end{equation*}
On the other hand, it is easy to see that \eqref{ass:min-dist-ineq} in   Assumption~\ref{ass:min-dist} yields that $\Delta_{j,\min}-2 \geq \frac{\Delta_{j,\min}}{2}-2 \geq  \frac{\Delta_{j,\min}}{6}$. Thus 
\begin{equation*}
\begin{split}
&\P[A_{n,\ell,2}] \leq \P\Big[\hat{w}_{j_\ell+1,\lceil\frac{j_\ell + j_{\ell+1}}{2} \rceil} \geq \frac{\Delta_{j,\min} \Delta_{\beta,\min}}{72}\Big]\\
&\hspace*{3cm}  + \P\Big[\Big|\bar M_n({j}_\ell+1;\lceil\frac{j_\ell + j_{\ell+1}}{2} \rceil -1) \Big|\geq \frac{\Delta_{j,\min} \Delta_{\beta,\min}}{72\sqrt{m}}\Big]\\
& \hspace*{1.5cm} :=\alpha_{n,\ell,2}^{(1)} + \alpha_{n,\ell,2}^{(2)}.
\end{split}
\end{equation*}
Using the property of the data-driven weights, we remark that 
\begin{equation*}
\alpha_{n,\ell,2}^{(1)} \leq  \P\Big[\hat{w}^2_{j_\ell+1} \geq \frac{\Delta_{j,\min}^2\Delta_{\beta,\min}^2}{{144}^2}\Big].
\end{equation*}
 By  (\ref{ass:consistency-2})  in Assumption~\ref{ass:consistency},   and \eqref{ctrlweigh} in Lemma~\ref{lem:control-martingale} with $\xi = \frac{n\Delta_{j,\min}^2\Delta_{\beta,\min}^2}{{144}^2 m\log m } + \E\big[\bar N_n\big((\frac{j_{\ell}}{m},1]\big)\big],$ it follows that 
\begin{eqnarray*}
\alpha_{n,\ell,2}^{(1)} &\leq&  2\exp\Bigg(-\frac{n\xi^2}{2\E\Big[\bar N_n\Big(\big(\frac{j_{\ell}}{m},1\big]\Big)\Big]+\frac{2}{3}\xi }\Bigg)\rightarrow 0,
\end{eqnarray*}
 as $n \rightarrow \infty.$ Similarly, using   (\ref{ass:consistency-2})  in Assumption~\ref{ass:consistency},  and \eqref{ctrlmart} in Lemma~\ref{lem:control-martingale} with $z=\frac{\Delta_{j,\min} \Delta_{\beta,\min}}{72\sqrt{m}}$, it implies that
\begin{eqnarray*}
\alpha_{n,\ell,2}^{(2)}&\leq& 2\exp\Bigg(-\frac{nz^2}{2\E\Big[\bar N_n\Big(\big(\frac{j_\ell}{m},\frac{\lceil\frac{j_\ell + j_{\ell+1}}{2} \rceil -1}{m}\big]\Big)\Big] +\frac{2}{3}z }\Bigg) \rightarrow 0,
\end{eqnarray*}
 as $n \rightarrow \infty.$ Therefore, we conclude that $\P[A_{n,\ell,2}]\rightarrow 0,$ as  $n\rightarrow \infty.$

\subsubsection{Step I.2.  Prove: $ \P[A_{n,\ell} \cap C^\complement_n] \rightarrow 0,$ as  $n\rightarrow \infty.$ } 
\label{subsection2:CASE-I} 

 Recall that 
$
C_n^\complement= \big\{ \max\limits_{1\leq k\leq L_0-1}|\hat{j}_\ell - j_\ell| \geq \frac{ \Delta_{j,\min}}{2}\big\}. 
$ 
 We split  $\P[A_{n,\ell} \cap C_n^\complement]$ in three terms as following
\begin{equation*} 
\P[A_{n,\ell} \cap C_n^\complement] = \P[A_{n,\ell} \cap D_n^{(l)}]+ \P[A_{n,\ell} \cap D_n^{(m)}]+ \P[A_{n,\ell} \cap D_n^{(r)}],
\end{equation*}
 where
 \begin{eqnarray*}
   D_n^{(l)} &:=& \{\textrm{there exists}\,  \ell\in \{1, \ldots, L_0-1\}: \hat{j}_\ell \leq j_{\ell-1} \} \cap C_n^\complement \\
   D_n^{(m)} &:=& \left\{ \textrm{for all}\,  \ell\in \left\{1, \ldots, L_0-1\right\}: j_{\ell-1} < \hat{j}_\ell < j_{\ell+1}\right\}\cap C_n^\complement, \\
   D_n^{(r)} &:=& \{\textrm{there exists}\,  \ell\in \{1, \ldots, L_0-1\}: \hat{j}_\ell \geq j_{\ell+1}\}\cap C_n^\complement.
 \end{eqnarray*}
 Let us first focus on $\P[A_{n,\ell} \cap D_n^{(m)}]$. Observe that
\begin{equation*}
\begin{split}
& \P[A_{n,\ell} \cap D_n^{(m)}] =  \P\Big[A_{n,\ell} \cap \{ \hat{j}_{\ell+1} - j_\ell  \geq \frac{\Delta_{j,\min}}{2}\}\cap D_n^{(m)}\Big] \\
& \hspace*{4cm}+ \P[A_{n,\ell} \cap \{ \hat{j}_{\ell+1} - j_\ell  <\frac{\Delta_{j,\min}}{2}\}\cap D_n^{(m)}].
\end{split}
\end{equation*}
The fact that $0\leq \hat{j}_{\ell+1} - j_\ell  <\frac{\Delta_{j,\min}}{2}$ yields  $j_{\ell +1} - \hat{j}_{\ell+1} \geq \frac{\Delta_{j,\min}}{2}.$ Then, it is easy to see that 
$j_{\ell +1} - \hat{j}_{\ell+1} = (j_{\ell +1} - j_\ell) - (\hat{j}_{\ell+1} - j_\ell) \geq \Delta_{j,\min} - \frac{\Delta_{j,\min}}{2} \geq \frac{\Delta_{j,\min}}{2}$.\\
Hence
\begin{equation*}
\begin{split}
&\P[A_{n,\ell} \cap D_n^{(m)}]  \leq \P\Big[A_{n,\ell} \cap \{ \hat{j}_{\ell+1} - j_\ell  \geq \frac{\Delta_{j,\min}}{2}\}\cap D_n^{(m)}\Big]\\
& \hspace*{4.5cm}  + \P\Big[A_{n,\ell} \cap \{ j_{\ell +1} -\hat{j}_{\ell+1}  \geq\frac{\Delta_{j,\min}}{2}\}\cap D_n^{(m)}\Big].
\end{split}
\end{equation*}
Moreover, we note that 
\begin{equation*}
\begin{split}
 & A_{n,\ell} \cap \Big\{ j_{\ell +1} -\hat{j}_{\ell+1}  \geq\frac{\Delta_{j,\min}}{2}\Big\}\cap D_n^{(m)}\\
&\hspace*{2.5cm} \subset \bigcup_{r=\ell +1}^{L_0-2}  \Big\{j_r - \hat{j}_r \geq\frac{\Delta_{j,\min}}{2}\Big\}\cap \Big\{\hat{j}_{r+1}-j_r \geq\frac{\Delta_{j,\min}}{2}\Big\} \cap  D_n^{(m)} .
\end{split}
\end{equation*}
Thus, we have 
 \begin{equation}
\label{prdnm}
 \P[A_{n,\ell} \cap D_n^{(m)}]\leq \P[A_{n,\ell}\cap B_{\ell+1,\ell} \cap D_n^{(m)}] + \sum_{s=\ell+1}^{L_0-2} \P[C_{s,s}\cap B_{s+1,s} \cap D_n^{(m)}],
 \end{equation}
 where
\begin{equation*}
 \left\{
      \begin{array}{ll}
        B_{p,q}= \{(\hat{j}_p - j_q) \geq \frac{ \Delta_{j,\min}}{2}\},  \\
        \textrm{with the convention}\, \, B_{L_0, L_0-1} = \{m - j_{L_0-1}\geq \frac{ \Delta_{j,\min}}{2}\} ,  \\
        C_{p,q}= \{( j_p - \hat{j}_q) \geq \frac{ \Delta_{j,\min}}{2}\}.
      \end{array}
    \right.
\end{equation*}
Let us now prove that the first term in the right hand side of $(\ref{prdnm})$ goes to zero as $n$ tends to infinity , the arguments for  the other terms being similar.  Using $(\ref{kkt})$ in Lemma~\ref{lem:KKT} with $j = j_\ell$ and $j = \hat{j}_\ell +1$, on the one hand and $(\ref{kkt})$ in Lemma~\ref{lem:KKT} with $j= j_\ell +1$ and  $j = \hat{j}_{\ell+1}$ on the other hand, we obtain, respectively
\begin{eqnarray}
\label{kktdnm1}
|\hat{j}_\ell - j_\ell -2||\hat{\beta}_{\hat{j}_{\ell+1}-1,m} -\beta_{0,j_\ell-1,m}| \leq \hat{w}_{\hat{j}_\ell +1,j_\ell}+ \big|\sqrt{m}\bar M_n(\hat{j}_\ell+1;j_\ell -1)\big|,
\end{eqnarray}
and
\begin{eqnarray}
\label{kktdnm2}
|\hat{j}_{\ell+1} - j_\ell -2||\hat{\beta}_{\hat{j}_{\ell+1}-1,m} -\beta_{0,j_{\ell+1}-1,m}| \leq \hat{w}_{j_\ell +1,\hat{j}_{\ell+1}} + | \sqrt{m}\bar M_n({j}_\ell+1;\hat{j}_{\ell+1} -1)|.
\end{eqnarray}
 Besides, we have
\begin{equation*}
  \begin{split}
   &|\beta_{0,j_{\ell+1}-1,m} - \beta_{0,j_\ell-1,m}| \\
   & \qquad \qquad= |( \hat{\beta}_{\hat{j}_{\ell+1}-1,m} - \beta_{0,j_\ell-1,m}) - (\hat{\beta}_{\hat{j}_{\ell+1}-1,m} - \beta_{0,j_{\ell+1}-1})|\\
   & \qquad \qquad\leq |\hat{\beta}_{\hat{j}_{\ell+1}-1,m} - \beta_{0,j_\ell-1,m}| + |\hat{\beta}_{\hat{j}_{\ell+1}-1,m} - \beta_{0,j_{\ell+1}-1,m}|\\
                        &\qquad \qquad\leq \frac{\hat{w}_{\hat{j}_\ell +1,j_\ell}}{|\hat{j}_\ell - j_\ell-2|} + \frac{\sqrt{m}\bar M_n(\hat{j}_\ell+1;j_\ell -1)|}{|\hat{j}_\ell - j_\ell-2|} \\
  &\hspace{4.1cm}+ \frac{ \hat{w}_{j_\ell +1,\hat{j}_{\ell+1}}}{|\hat{j}_{\ell+1} - j_\ell-2|} + \frac{| \sqrt{m}\bar M_n({j}_\ell+1;\hat{j}_{\ell+1} -1)|}{|\hat{j}_{\ell+1} - j_\ell-2|}\\
  & \qquad \qquad\leq \frac{\hat{w}_{\hat{j}_\ell +1,j_\ell}}{\frac{m\varepsilon_n}{6}} + \frac{\sqrt{m}\bar M_n(\hat{j}_\ell+1;j_\ell -1)|}{|\hat{j}_\ell - j_\ell-2|} \\
  &\hspace{4.1cm}+ \frac{ \hat{w}_{j_\ell +1,\hat{j}_{\ell+1}}}{\frac{\Delta_{j,\min}}{6}} + \frac{| \sqrt{m}\bar M_n({j}_\ell+1;\hat{j}_{\ell+1} -1)|}{|\hat{j}_{\ell+1} - j_\ell-2|}.
\end{split}
\end{equation*}
Define the event $E_{n,\ell}$ by
\begin{equation*}
   \begin{split}
     \label{eq:40}
      & E_{n,\ell}=\Bigg\{\big| \beta_{0,j_{\ell+1}-1,m} - \beta_{0,j_\ell-1,m}\big| \leq \frac{\hat{w}_{\hat{j}_\ell +1,j_\ell} }{\frac{m\varepsilon_n}{6}} + \frac{\hat{w}_{j_\ell +1,\hat{j}_{\ell+1}}}{\frac{\Delta_{j,\min}}{6}}\\
      &\hspace{6cm} +\Big|\frac{\sqrt{m}\bar M_n(\hat{j}_\ell+1;j_\ell -1)}{\hat{j}_\ell - j_\ell-2 }\Big|\\
       &\hspace{6cm} + \Big|\frac{\sqrt{m}\bar M_n({j}_\ell+1;\hat{j}_{\ell+1} -1)}{\hat{j}_{\ell+1} - j_\ell-2}\Big|\Bigg\}.
     \end{split}
\end{equation*}
We observe that $E_{n,\ell}$ occurs with probability  one. Therefore, we obtain
\begin{flalign*}
      & \P[A_{n,\ell}\cap B_{\ell+1,\ell} \cap D_n^{(m)}]\\
      &\qquad \qquad\leq \P\Big[E_{n,\ell}\cap \{(j_\ell - \hat{j}_\ell) > \frac{m\varepsilon_n}{2}\}\cap\{(\hat{j}_{\ell+1} - j_\ell) \geq \frac{ \Delta_{j,\min}}{2} \}\Big] \\
      &\qquad \qquad\leq \P\Big[\hat{w}_{\hat{j}_\ell +1,j_\ell}\geq \frac{ m\varepsilon_n|\beta_{0,j_{\ell+1}-1,m} -\beta_{0,j_\ell-1,m}|}{24}\Big] \\
      & \qquad \qquad\quad +\P\Big[\hat{w}_{j_\ell +1,\hat{j}_{\ell+1}}  \geq \frac{\Delta_{j,\min}| \beta_{0,j_{\ell+1}-1,m} - \beta_{0,j_\ell-1,m}|}{24}\Big]\\
      &\qquad \qquad\quad  + \P\Bigg[ \Bigg\{\Big|\frac{\sqrt{m}\bar M_n(\hat{j}_\ell+1;j_\ell -1)}{j_\ell - \hat{j}_\ell-2}\Big| \\
     &\hspace{4cm} \geq \frac{| \beta_{0,j_{\ell+1}-1,m} - \beta_{0,j_\ell-1,m}|}{4}\Bigg\}\bigcap \Big\{j_\ell - \hat{j}_\ell -2\geq \frac{m\varepsilon_n}{6}\Big\}\Bigg]\\  
     &\qquad \qquad\quad    + \P\Bigg[\Bigg\{\Big|\frac{\sqrt{m}\bar M_n({j}_\ell+1;\hat{j}_{\ell+1} -1)}{\hat{j}_{\ell+1} - j_\ell-2}\Big| \\
&\hspace{4cm}\geq \frac{| \beta_{0,j_{\ell+1}-1,m} - \beta_{0,j_\ell-1,m}|}{4}\Bigg\} \bigcap \Big\{\hat{j}_{\ell+1} - j_\ell -2 \geq \frac{ \Delta_{j,\min}}{6}\Big\}\Bigg].\\
    &\qquad \qquad := \theta_{n,\ell,1} + \theta_{n,\ell,2} +\theta_{n,\ell,3} +\theta_{n,\ell,4} 
\end{flalign*}
We note
\begin{equation*}
 \theta_{n,\ell,1} \leq \P\Big[\hat{w}_{\hat{j}_\ell +1}\geq \frac{ m\varepsilon_n\Delta_{\beta,\min}}{48}\Big] \leq \P\Big[\hat{w}^2_{j_{\ell -1}+1}\geq \frac{ m^2\varepsilon_n^2\Delta_{\beta,\min}^2}{{48}^2}\Big].
\end{equation*}
Using (\ref{ass:consistency-1})  in Assumption~\ref{ass:consistency}, and  \eqref{ctrlweigh} in Lemma~\ref{lem:control-martingale}  with $\xi = \frac{n m\varepsilon_n^2\Delta_{\beta,\min}^2}{48^2\log m} + \E\Big[\bar N_n\Big(\big(\frac{j_{\ell -1}}{m},1\big]\Big)\Big],$ we get
\begin{eqnarray*}
\theta_{n,\ell,1} &\leq&  2\exp\Bigg(-\frac{n\xi^2}{2 \E\Big[\bar N_n\Big(\big(\frac{j_{\ell -1}}{m},1\big]\Big)\Big] +\frac{2}{3}\xi }\Bigg)\rightarrow 0,
\end{eqnarray*}
 as $n \rightarrow \infty.$
Analogously,  
\begin{equation*}
 \theta_{n,\ell,2} \leq \P\Big[\hat{w}_{\hat{j}_\ell +1}\geq \frac{ \Delta_{j,\min}\Delta_{\beta,\min}}{48}\Big]\leq \P\Big[\hat{w}^2_{j_{\ell}+1}\geq \frac{ \Delta_{j,\min}^2\Delta_{\beta,\min}^2}{48^2}\Big].
\end{equation*}
Using (\ref{ass:consistency-2})  in Assumption~\ref{ass:consistency}, and  \eqref{ctrlweigh} in Lemma~\ref{lem:control-martingale}, with $\xi = \frac{n \Delta_{j,\min}^2\Delta_{\beta,\min}^2}{48^2m\log m } + \E\big[\bar N_n\big((\frac{j_{\ell}}{m},1]\big)\big],$ we have
\begin{eqnarray*}
\theta_{n,\ell,2} &\leq&  2\exp\Bigg(-\frac{n\xi^2}{2\E\Big[\bar N_n\Big(\big(\frac{j_{\ell}}{m},1\big]\Big)\Big] +\frac{2}{3}\xi }\Bigg)\rightarrow 0,
\end{eqnarray*}
 as $n \rightarrow \infty.$
Furthermore, using \eqref{ctrlmart} in Lemma~\ref{lem:control-martingale}, we have 
\begin{eqnarray*}
\theta_{n,\ell,3} &\leq& \P\Big[\left|\bar M_n\left(\hat{j}_\ell+1;j_\ell -1\right)\right| \geq \frac{m\varepsilon_n\Delta_{\beta,\min}}{24\sqrt{m}} \Big]\\
&\leq&2 \exp\Bigg(-\frac{n\psi_n^2}{2\E\Big[\bar N_n\Big(\big(\frac{j_{\ell -1}+1}{m},\frac{j_\ell-1}{m}\big]\Big)\Big] +\frac{2}{3}\psi_n} +\log m\Bigg), 
\end{eqnarray*}
where $\psi_n = \frac{\sqrt{m}\varepsilon_n\Delta_{\beta,\min}}{24}.$ By (\ref{ass:consistency-1})  in Assumption~\ref{ass:consistency}, we get that $\theta_{n,\ell,3}\rightarrow 0,\textrm{ as } n \rightarrow \infty.$
Similarly, using \eqref{ctrlmart} in Lemma~\ref{lem:control-martingale}, we have 
\begin{eqnarray*}
\theta_{n,\ell,4} &\leq& \P\Big[\big|\bar M_n({j}_\ell+1;\hat{j}_{\ell +1}-1)\big| \geq \frac{\Delta_{j,\min}\Delta_{\beta,\min}}{24\sqrt{m}} \Big]\\
&\leq&2 \exp\Bigg(-\frac{n\delta_n^2}{2\E\Big[\bar N_n\Big(\big(\frac{j_\ell}{m},\frac{j_{\ell+2}-2}{m}\big]\Big)\Big] +\frac{2}{3}\delta_n } +\log m\Bigg), 
\end{eqnarray*}
 where $\delta_n = \frac{\Delta_{j,\min}\Delta_{\beta,\min}}{24\sqrt{m}}$. By (\ref{ass:consistency-2})  in Assumption~\ref{ass:consistency}, we get that $\theta_{n,\ell,4}\rightarrow 0,\textrm{as } n \rightarrow \infty.$  Consequently, we obtain $\P[A_{n,\ell}\cap B_{\ell+1,\ell} \cap D_n^{(m)}] \rightarrow 0$ as $n\rightarrow \infty$. Now, we have $\P[A_{n,\ell}\cap {D_n}^{(l)}] \leq \P[{D_n}^{(l)}],$ and  
\begin{eqnarray*}
 \P[{D_n}^{(l)}]&=& \P\Big[\left\{\exists\, \ell\in \{1,\ldots, L_0-1\}: \hat{j}_\ell \leq j_{\ell-1}\right\}\cap C_n^\complement\Big] \\
  &=& \P\Big[\Bigg\{\bigcup_{ \ell=1}^{L_0-1}\max\{1 \leq q \leq L_0-1: \hat{j_q} \leq j_{q-1}\} = \ell\Bigg\} \cap C_n^\complement\Big]\\
  &=& \sum_{\ell=1}^{L_0-1}\P\Big[\left\{\max\{1 \leq q \leq L_0-1: \hat{j_q} \leq j_{q-1}\} = \ell\right\}\cap C_n^\complement\Big].
\end{eqnarray*}
We note that on  the event $\{ \max\{1 \leq q \leq L_0-1; \hat{j_q} \leq j_{q-1}\} = \ell\},$ it is clear to see that $\hat{j}_\ell \leq j_{\ell-1}$ and $\hat{j}_{q+1} > j_q$ for all $q=\ell, \ldots, L_0-1.$ Then,  it follows that 
\begin{eqnarray*}
 \P[{D_n}^{(l)}]&\leq& \sum_{\ell=1}^{L_0-1}2^{\ell -1} \P\Big[\bigcap_{q\geq \ell}^{L_0-1} \{\hat{j}_\ell \leq j_{\ell-1}\} \cap\{\hat{j}_{q+1} > j_q\} \Big].
\end{eqnarray*}
In addition, we note that 
\begin{flalign*}
& \bigcap_{q\geq \ell}^{L_0-1} \{\hat{j}_\ell \leq j_{\ell-1}\}\cap \{\hat{j}_{q+1} > j_q\} \\
&\subset   \{j_\ell \leq \hat{j_\ell} \}\cap \Big( \{\hat{j}_{\ell+1} > \frac{j_\ell+j_{\ell+1}}{2}\} \cup  \{ \hat{j}_{\ell+1} < \frac{j_\ell+j_{\ell+1}}{2}\}\Big)\\
&\qquad  \cap \Big( \{\hat{j}_{\ell+2} > \frac{j_{\ell+1}+j_{\ell+2}}{2}\} \cup  \{ \hat{j}_{\ell+2} < \frac{j_{\ell+2}+j_{\ell+1}}{2}\}\Big)\\
&\qquad \cap \ldots. \cap \Big(\{\hat{j}_{L_0-1} > \frac{j_{L_0-2}+j_{L_0+1}}{2}\} \cup  \{ \hat{j}_{L_0-1} < \frac{j_{L_0-2}+j_{L_0+1}}{2}\}\Big) \\
&\qquad \cap \Big(\{\hat{j}_{L_0} > \frac{j_{L_0-1}+j_{L_0}}{2}\} \cup  \{ \hat{j}_{L_0} < \frac{j_{L_0-1}+j_{L_0}}{2}\}\Big)\\
&\subset  \{j_\ell - \hat{j_\ell} > \frac{\Delta_{j,\min}}{2} \}\cap  \Big( \{\hat{j}_{\ell+1} - j_\ell > \frac{\Delta_{j,\min}}{2}\} \cup \{j_{\ell+1} - \hat{j}_{\ell+1} >  \frac{\Delta_{j,\min}}{2} \}\Big)\\
&\qquad  \cap \Big( \{\hat{j}_{\ell+2} - j_{\ell+1}> \frac{\Delta_{j,\min}}{2} \} \cup  \{ j_{\ell+2}- \hat{j}_{\ell+2} > \frac{\Delta_{j,\min}}{2} \}\Big)\\
&\qquad \cap \ldots. \cap \Big(  \{\hat{j}_{L_0-1} - j_{L_0-2} >  \frac{\Delta_{j,\min}}{2} \} \cup \{j_{L_0-2}-\hat{j}_{L_0-2} > \frac{\Delta_{j,\min}}{2} \}\Big)\\
&\qquad  \cap \Big( \{\hat{j}_{L_0} - j_{L_0 -1} > \frac{\Delta_{j,\min}}{2} \} \cup  \{j_{L_0-1}-\hat{j}_{L_0-1} > \frac{\Delta_{j,\min}}{2} \}\Big)\\
&\subset \bigcup_{q=\ell}^{L_0-2} \Big( \{j_q - \hat{j_q} > \frac{\Delta_{j,\min}}{2} \}\cap\{\hat{j}_{q+1} - j_q > \frac{\Delta_{j,\min}}{2}\}\Big)\cup  \{j_{L_0-1}-\hat{j}_{L_0-1} > \frac{\Delta_{j,\min}}{2}\}.
\end{flalign*}
Hence
\begin{equation}
\label{prdnl}
\begin{split}
& \P[{D_n}^{(l)}]\leq 2^{L_0 -2}\sum_{\ell=1}^{L_0-2}\sum_{q= \ell}^{L_0 -2}\P\Big[\big\{(j_q - \hat{j}_q) >\frac{ \Delta_{j,\min}}{2}\}\cap\{\hat{j}_{q+1}- j_{q} > \frac{ \Delta_{j,\min}}{2}\big\}\Big]\\
   & \hspace{3cm}+ 2^{L_0 -2} \P\Big[\{j_{L_0-1} - \hat{j}_{L_0-1} > \frac{ \Delta_{j,\min}}{2}\}\Big].
\end{split}
\end{equation}
Consider the first term of the sum in the right-hand side of $(\ref{prdnl})$. Using (\ref{kktdnm1}) and $(\ref{kktdnm2})$ with $\ell=q$, we obtain
\begin{flalign*}
&\P\Big[\Big\{ j_{q} - \hat{j}_q > \frac{ \Delta_{j,\min}}{2}  \Big\}\cap \Big\{ \hat{j}_{q+1} - j_q > \frac{ \Delta_{j,\min}}{2}\Big\}\Big]\\
& \qquad \leq  \P\Big[\frac{\hat{w}_{\hat{j}_q +1,j_q}}{\frac{\Delta_{j,\min}}{6}}  \geq \frac{|\beta_{0,j_{q+1}-1,m} - \beta_{0,j_q-1,m}|}{4} \Big] \\
& \qquad\quad + \P\Big[\frac{\hat{w}_{j_q+1,\hat{j}_{q+1}}}{\frac{\Delta_{j,\min}}{6}} \geq \frac{|\beta_{0,j_{q+1}-1,m} - \beta_{0,j_q-1,m}|}{4} \Big] \\
&\qquad \quad + \P\Bigg[\Bigg\{ \Big|\frac{ \sqrt{m}\bar M_n(\hat{j}_q+1; j_q - 1)}{j_q - \hat{j}_q - 2 }\Big|\\
& \hspace{4cm}\geq \frac{|\beta_{0,j_{q+1}-1,m} - \beta_{0,j_q-1,m}|}{4} \Bigg\}\bigcap \Big\{ j_q - \hat{j}_q \geq \frac{ \Delta_{j,\min}}{2}  \Big\}\Bigg] \\
& \qquad\quad + \P\Bigg[ \Bigg\{ \Big| \frac{\sqrt{m}\bar M_n(j_q+1; \hat{j}_{q+1}-1)}{\hat{j}_{q+1} - j_q - 2}\Big|\\
&\hspace{4cm} \geq \frac{|\beta_{0,j_{q+1}-1,m} - \beta_{0,j_q-1,m}|}{4}\Bigg\}\bigcap \Big\{ \hat{j}_{q+1} - {j}_q \geq \frac{ \Delta_{j,\min}}{2}\Big\}\Bigg].\\
&\qquad := \theta_{n,q,1} + \theta_{n,q,2} +\theta_{n,q,3} +\theta_{n,q,4}.
\end{flalign*}
By \eqref{ctrlmart}-\eqref{ctrlweigh} in Lemma~\ref{lem:control-martingale},  and (\ref{ass:consistency-1})-(\ref{ass:consistency-2}) in Assumption~\ref{ass:consistency}, we show that for $s=1, \ldots, 4,\, \theta_{n,q,s }  \rightarrow 0,$ as $n$ tending to infinity. Then
\begin{equation*}
\P\Big[ \Big\{ j_{q} - \hat{j}_q > \frac{ \Delta_{j,\min}}{2}  \Big\}\bigcap \Big\{ \hat{j}_{q+1} - j_q > \frac{ \Delta_{j,\min}}{2}\Big\}\Big] \rightarrow 0.
\end{equation*}
 Let us now consider the last term in the right hand of $(\ref{prdnl})$. Using the observations $(\ref{kktdnm1})$ and $(\ref{kktdnm2})$ with $\ell=L_0-1$ leads to
\begin{flalign*}
       &\P\Big[ \Big\{ j_{L_0-1} - \hat{j}_{L_0-1} > \frac{ \Delta_{j,\min}}{2}  \Big\}\Big]\\
      &\qquad \leq  \P\Bigg[\frac{\hat{w}_{\hat{j}_{L_0-1} +1,j_{L_0-1}}}{\frac{m\varepsilon_n}{6} }  \geq \frac{|\beta_{0,j_{L_0}-1,m} - \beta_{0,j_{L_0-1}-1,m}|}{4} \Bigg] \\
     &\qquad \quad \, + \P\Bigg[\frac{\hat{w}_{j_{L_0-1} +1,m}}{\frac{\Delta_{j,\min}}{6}} \geq \frac{|\beta_{0,j_{L_0}-1,m} - \beta_{0,j_{L_0-1}-1,m}|}{4} \Bigg] \\
     &\qquad \quad \, + \P\Bigg[ \Bigg\{ \Big|\frac{ \sqrt{m}\bar M_n(\hat{j}_{L_0-1}+1; j_{L_0-1} - 1)}{j_{L_0-1} - \hat{j}_{L_0-1} -2 }\Big|\\
&\hspace{3.4cm} \geq \frac{|\beta_{0,j_{L_0}-1,m} - \beta_{0,j_{L_0-1}-1,m}|}{4} \Bigg\}\bigcap \Big\{ j_{L_0-1} - \hat{j}_{L_0-1} \geq \frac{ \Delta_{j,\min}}{2}  \Big\}\Bigg] \\
      &\qquad \quad \, + \P\Bigg[ \Bigg\{ \Big| \frac{\sqrt{m}\bar M_n(j_{L_0-1}+1; m-1)}{m - j_{L_0-1}-2}\Big|  \geq \frac{|\beta_{0,j_{L_0}-1,m} - \beta_{0,j_{L_0-1}-1,m}|}{4}\Bigg\}\Bigg]\\
      &\qquad := \theta_{n,{L_0-1},1} + \theta_{n,{L_0-1},2} +\theta_{n,{L_0-1},3} +\theta_{n,{L_0-1},4}.
\end{flalign*}
By \eqref{ctrlmart}-\eqref{ctrlweigh} in Lemma~\ref{lem:control-martingale},  and (\ref{ass:consistency-1})-(\ref{ass:consistency-2}) in Assumption~\ref{ass:consistency}, we show that for $s=1, \ldots, 4,$ we obtain $ \theta_{n,{L_0 -1} ,s } \rightarrow 0,$ as $n \rightarrow \infty.$ Then
\begin{equation*}
\P\Big[ \Big\{ j_{L_0 -1} - \hat{j}_{L_0 -1}  > \frac{ \Delta_{j,\min}}{2}  \Big\}\bigcap \Big\{ m - j_{L_0 -1}  > \frac{ \Delta_{j,\min}}{2}\Big\}\Big] \rightarrow 0.
\end{equation*}
This implies that $ \P[{D_n}^{(l)}]  \rightarrow 0,$ as $n \rightarrow \infty.$
Similarly, we  prove that $\P[{D_n}^{(r)}]
\rightarrow 0,$ as  $n \rightarrow \infty$ which yields that $\P[A_{n,\ell}
\cap C_n^\complement] \rightarrow 0,$  as   $n \rightarrow \infty.$ This concludes the proof of
Theorem~\ref{thm3}, up to the case $\{\hat j_\ell > j_\ell\}$ for a
fixed $\ell \in \{1, \ldots, L_0-1\}$ which is given in
Appendix~\ref{app:proof-thm3-case2}.

\section{Proof of  Theorem  \ref{thm4}}
\label{sec:proof-thm4}

This proof is based on the same arguments  in the proof of Theorem~\ref{thm3}.
Let $\mathcal{T}^{\textrm{approx}}_{0 } = \big\{\frac{j_1}{m}, \ldots, \frac{j_{L_0-1}}{m}\big\}$ be the set of the true  approximate change-points.
First, we note that 
\begin{equation*}
\P\Big[\cE(\hat{\mathcal{T}}\|\cT_0) > \varepsilon_n \Big]
\leq \P\Big[\cE(\hat{\mathcal{T}}\|\mathcal{T}^{\textrm{approx}}_{0}) > \varepsilon_n\Big] + \P\Big[\cE(\mathcal{T}^{\textrm{approx}}_{0})\| \mathcal{T}_{0})> \varepsilon_n\Big].
\end{equation*}
Obviously, since $ m\varepsilon_n \geq 6,$ we have $\P\Big[\cE(\mathcal{T}^{\textrm{approx}}_{0})\| \mathcal{T}_0)> \varepsilon_n\Big] = 0$. It is clear to remark that the inequality  $\hat{L} \leq m$ holds true.
In order to prove that
\begin{equation*}
\P\Big[\Big\{\cE(\hat{\mathcal{T}}\|\mathcal{T}^{\textrm{approx}}_{0}) > \varepsilon_n \Big\}\bigcap \Big\{\hat{L}\geq L_0-1\Big\}\Big] \rightarrow 0,
\end{equation*}
as $n\rightarrow \infty,$ it is enough to prove that
\begin{equation*}\P\Big[\Big\{\cE(\hat{\mathcal{T}}\|\mathcal{T}^{\textrm{approx}}_{0})  > \varepsilon_n \Big\}\bigcap \Big\{ L_0-1 \leq \hat{L} \leq m\Big\}\Big]\rightarrow 0, 
\end{equation*}
as $n\rightarrow \infty.$
We have that
\begin{align}
\nonumber
 \label{prnbr}
&\P\Big[\Big\{\cE(\hat{\mathcal{T}}\|\mathcal{T}^{\textrm{approx}}_{0}) > \varepsilon_n \Big\}\bigcap \Big\{ L_0-1\leq \hat{L} \leq m\Big\}\Big]\\
\nonumber
&\leq  \P\Big[\Big\{\cE(\hat{\mathcal{T}}\|\mathcal{T}^{\textrm{approx}}_{0}) > \varepsilon_n\Big\} \bigcap \Big\{\ind{\hat{L} = L_0-1} \Big\}  \Big] +\P\Big[\Big\{\cE(\hat{\mathcal{T}}\|\mathcal{T}^{\textrm{approx}}_{0}) > \varepsilon_n\Big\} \bigcap \Big\{\ind{\hat{L} > L_0-1} \Big\}  \Big]\\
\nonumber
&\leq  \P\Big[\Big\{\cE(\hat{\mathcal{T}}\|\mathcal{T}^{\textrm{approx}}_{0}) > \varepsilon_n\Big\} \bigcap \Big\{\ind{\hat{L} = L_0-1} \Big\}  \Big] + \sum_{L = L_0}^{m} \P\Big[\cE(\hat{\mathcal{T}}\|\mathcal{T}^{\textrm{approx}}_{0})  > \varepsilon_n \Big]\\
&\leq \P\Big[\Big\{\cE(\hat{\mathcal{T}}\|\mathcal{T}^{\textrm{approx}}_{0}) > \varepsilon_n\Big\} \bigcap \Big\{\ind{\hat{L} = L_0-1} \Big\}  \Big] \\
\nonumber
&\hspace{1cm}+ \sum_{L = L_0}^{m} \sum_{\ell=1}^{L_0-1} \P\Big[\forall q \in\{1, \ldots, L\}, |\frac{\hat{j}_q}{m} - \frac{j_\ell}{m}|> \varepsilon_n\Big]
\end{align}
 The first term of the right-hand side of $(\ref{prnbr})$ tends to zero as $n$ tends to infinity since it is upper bounded by $\P\Big[\max_{1\leq \ell\leq L_0-1}|\hat{j}_\ell - j_\ell| > m\varepsilon_n\Big]$  which tends to zero by the proof of Theorem~\ref{thm3}. Let us now focus on the second term on the right-hand side of $(\ref{prnbr})$. Note that
\begin{flalign*}
&\sum_{L = L_0}^{m} \sum_{\ell=1}^{L_0-1} \P\Big[\forall\, 1\leq q \leq L,\, |\hat{j}_q - j_\ell| > m\varepsilon_n\Big]
 := \sum_{L = L_0}^{m} \sum_{\ell=1}^{L_0-1} \P\big[R_{n,\ell,1}\big] + \P\big[R_{n,\ell,2}\big] + \P\big[R_{n,\ell,3}\big],
\end{flalign*}
where
\begin{equation*}
\begin{split}
  & R_{n,\ell,1}:= \Big\{\forall\, 1\leq q \leq L:\, |\hat{j}_q - j_\ell| > m\varepsilon_n\textrm{ and } \hat{j}_q < j_\ell\Big\} \\
  & R_{n,\ell,2}:= \Big\{\forall\, 1\leq q \leq L:\, |\hat{j}_q - j_\ell| > m\varepsilon_n\textrm{ and } \hat{j}_q > j_\ell \Big\} \\
 & R_{n,\ell,3}:= \Big\{ \exists\, 1\leq q \leq L-1:\, \big\{|\hat{j}_q - j_\ell| > m\varepsilon_n\big\},
\big\{|\hat{j}_{q+1} - j_\ell|> m\varepsilon_n\big\},\textrm{ and } \big\{\hat{j}_q < j_\ell < \hat{j}_{q+1}\big\}\Big\}.
\end{split}
\end{equation*}
Note that
\begin{equation*}
\P\big[R_{n,\ell,1}\big] = \P\Big[R_{n,\ell,1}\cap \big\{\hat{j}_L > j_{\ell-1}\big\}\Big] + \P\Big[R_{n,\ell,1}\cap \big\{\hat{j}_L\leq j_{\ell-1}\big\}\Big].
\end{equation*}
By applying $(\ref{kkt})$  in Lemma~\ref{lem:KKT} with $j = j_\ell$ and with $j = \hat{j}_L+1$ in the case where $\hat{j}_L > j_{\ell-1},$  it follows that, with probability one,
\begin{equation*}
\begin{split}
&\Big|\big(j_\ell - \hat{j}_L-2\big) \big(\big( \beta_{0,j_\ell-1,m} - \beta_{0,j_{\ell+1}-1,m}\big) \\
&\hspace{3cm}+ \big(\beta_{0,j_{\ell+1}-1,m} - \hat{\beta}_{\hat{j}_{L+1}-1,m}\big) \big) \\
&\hspace{4cm}+ \sqrt{m}\bar M_n(\hat{j}_L+1; j_\ell -1) \Big| \leq \hat{w}_{\hat{j}_L+1,j_\ell}.
\end{split}
\end{equation*}
Thus
\begin{flalign*}
&\P\Big[R_{n,\ell,1}\cap \big\{\hat{j}_L > j_{\ell-1}\big\}\Big]\\
& \leq  \P\bigg[\Big\{\frac{\hat{w}_{\hat{j}_L+1,j_\ell}}{m\varepsilon_n-2} \geq \frac{| \beta_{0,j_{\ell+1}-1,m} - \beta_{0,j_\ell-1,m}|}{3}\Big\} \cap \{\hat{j}_L >  j_{\ell-1}\}\bigg] \\
&\quad  + \P\bigg[\Big\{ |\hat{\beta}_{\hat{j}_{L+1}-1,m} -  \beta_{0,j_{\ell+1}-1,m}| \geq \frac{| \beta_{0,j_{\ell+1}-1,m} - \beta_{0,j_\ell-1,m}|}{3}\Big\}\bigg] \\
& \quad   + \P\bigg[\Big\{\Big|\frac{\bar M_n(\hat{j}_L+1; j_\ell-1)}{j_\ell - \hat{j}_L-2}\Big|\geq \frac{| \beta_{0,j_{\ell+1}-1,m} - \beta_{0,j_\ell-1,m}|}{3\sqrt{m}} \Big\}\cap \Big\{ |j_\ell - \hat{j}_L| \geq m\varepsilon_n\Big\}\bigg] \\
&  := \P\big[R^{(1)}_{n,\ell,1}\big] + \P\big[R^{(2)}_{n,\ell,1}\big] + \P\big[R^{(3)}_{n,\ell,1}\big].
\end{flalign*}
Using (\ref{ass:consistency-1}) in Assumption~\ref{ass:consistency},  and \eqref{ctrlmart}-\eqref{ctrlweigh} in Lemma~\ref{lem:control-martingale}  with $\xi = \frac{nm\varepsilon_n^2\Delta_{\beta,\min}^2}{36^2\log m} + \E\big[\bar N_n\big((\frac{j_{\ell-1}}{m},1]\big)\big],$  we prove that 
 $\sum_{L = L_0}^{m} \sum_{\ell=1}^{L_0-1}\P\big[R^{(3)}_{n,\ell,1}\big] \rightarrow 0, \, as \, n \rightarrow \infty.$ Let us now consider to $\P\big[R^{(2)}_{n,\ell,2}\big]$.
Using $(\ref{kkt})$ in Lemma~\ref{lem:KKT} with $j = j_\ell+1$ and with $j = j_{\ell+1},$ we get
\begin{equation*} 
\big(j_{\ell+1} - j_\ell-2 \big)\big| \hat{\beta}_{\hat{j}_{L+1}-1,m} - \beta_{0,j_{\ell+1}-1,m}\big| \leq \hat{w}_{j_\ell+1,j_{\ell+1}} + \big|\sqrt{m} \bar M_n(j_\ell+1; j_{\ell+1} - 1)\big|.
\end{equation*}
Therefore, we may upper bound $\P\big[R^{(2)}_{n,\ell,2}\big]$ as follows:
\begin{flalign*}
  &\P\bigg[|\hat{\beta}_{\hat{j}_{L+1}-1,m} -  \beta_{0,j_{\ell+1}-1,m}| \geq \frac{| \beta_{0,j_{\ell+1}-1,m} - \beta_{0,j_\ell-1,m}|}{3}\bigg]\\
   &\qquad \qquad \leq \P\bigg[\hat{w}_{j_\ell+1,j_{\ell+1}} \geq (j_{\ell+1} - j_\ell-2)\frac{| \beta_{0,j_{\ell+1}-1,m} - \beta_{0,j_\ell-1,m}|}{6}\bigg] \\
   & \qquad \qquad  \quad+ \P\bigg[\Big|\frac{\bar M_n(j_\ell +1; j_{\ell+1} -1)}{j_{\ell+1} - j_\ell-2}\Big|\geq \frac{| \beta_{0,j_{\ell+1}-1,m} - \beta_{0,j_\ell-1,m}|}{6\sqrt{m}} \bigg].
\end{flalign*}
\noindent By using Lemma~\ref{lem:KKT}, and   \eqref{ass:consistency-1}-\eqref{ass:consistency-2} in Assumption~\ref{ass:consistency}, we conclude  that $\sum_{L = L_0}^{m} \sum_{\ell=1}^{L_0-1} \P\big(R^{(2)}_{n,\ell,1}\big) \rightarrow 0, \, as \, n \rightarrow \infty.$ Analogously, it can be  shown that $\sum_{L = L_0}^{m} \sum_{\ell=1}^{L_0-1}\P\Big[R_{n,\ell,1}\cap \big\{\hat{j}_L \leq j_{\ell-1}\big\}\Big] \rightarrow 0, \, as \, n \rightarrow \infty.$
Moreover, we  prove, similarly, that  $\sum_{L = L_0}^{m} \sum_{\ell=1}^{L_0-1}\P\big[R_{n,\ell,2}\big]\rightarrow 0, \, as \, n \rightarrow \infty.$ Let us now focus on $\sum_{L = L_0}^{m} \sum_{\ell=1}^{L_0-1} \P\big[R_{n,\ell,3}\big].$ Note that $ \P\big[R_{n,\ell,3}\big]$ can be split in four terms as follows:
\begin{equation*}
\P\big[R_{n,\ell,3}\big] = \P\big[R^{(1)}_{n,\ell,3}\big] + \P\big[R^{(2)}_{n,\ell,3}\big] + \P\big[R^{(3)}_{n,\ell,3}\big] + \P\big[R^{(4)}_{n,\ell,3}\big],
\end{equation*}
where
\begin{eqnarray*}
  R^{(1)}_{n,\ell,3}&:=& R_{n,\ell,3} \cap \Big\{j_{\ell-1} < \hat{j}_q < \hat{j}_{q+1} < j_{\ell+1}\Big\} \\
  R^{(2)}_{n,\ell,3}&:=& R_{n,\ell,3} \cap \Big\{j_{\ell-1} < \hat{j}_q < {j}_{\ell+1}, \hat{j}_{q+1} \geq j_{\ell+1} \Big\} \\
  R^{(3)}_{n,\ell,3}&:=& R_{n,\ell,3} \cap \Big\{\hat{j}_q \leq j_{\ell-1}, j_{\ell-1} < \hat{j}_{q+1} < j_{\ell+1}\Big\} \\
  R^{(4)}_{n,\ell,3}&:=& R_{n,\ell,3} \cap \Big\{\hat{j}_q \leq j_{\ell-1}, j_{\ell+1} \leq \hat{j}_{q+1} \Big\}.
\end{eqnarray*}
We have to use Lemma~\ref{lem:KKT} twice. For $\P\big[R^{(1)}_{n,\ell,3}\big]$, we first use $(\ref{kkt})$  in Lemma~\ref{lem:KKT} with $j = j_\ell$ and $j= \hat{j}_q +1$, respectively, which gives with probability one
\begin{equation}
\label{eq:ch1}
\Big|\big(j_\ell - \hat{j}_q - 2\big) \big( \beta_{0,j_\ell-1,m}  - \hat{\beta}_{\hat{j}_{q+1}-1,m}\big) + \sqrt{m}\bar M_n(\hat{j}_q+1; j_\ell -1) \Big| \leq \hat{w}_{\hat{j}_q+1,j_\ell}.
\end{equation}
Thus, 
\begin{equation*}
\big| \beta_{0,j_\ell-1,m}  - \hat{\beta}_{\hat{j}_{q+1}-1,m}\big| \leq \frac{\hat{w}_{\hat{j}_q+1,j_\ell}}{j_\ell - \hat{j}_q -2} +\big| \frac{\sqrt{m}\bar M_n(\hat{j}_q+1; j_\ell -1) }{j_\ell - \hat{j}_q-2}\big|.
\end{equation*}
Second, we use $(\ref{kkt})$  in Lemma~\ref{lem:KKT} with $j = j_\ell +1$ and $j=\hat{j}_{q+1}$, respectively, to get with probability one
\begin{equation*} \Big|\big( \hat{j}_{q+1} - j_\ell - 2\big) \big\{\big( \beta_{0,j_{\ell+1}-1,m} - \hat{\beta}_{\hat{j}_{q+1} -1,m}\big) \big\} + \sqrt{m}\bar M_n(j_\ell +1; \hat{j}_{q+1} -1) \Big| \leq \hat{w}_{j_\ell+1,\hat{j}_{q+1}}.\end{equation*}
Hence
\begin{equation*} 
\big| \beta_{0,j_{\ell+1}-1,m} - \hat{\beta}_{\hat{j}_{q+1}-1,m}\big| \leq  \frac{ \hat{w}_{j_\ell+1,\hat{j}_{q+1}}}{ \hat{j}_{q+1} - j_\ell-2} +\big| \frac{ \sqrt{m}\bar M_n(j_\ell+1; \hat{j}_{q+1} -1) }{ \hat{j}_{q+1} - j_\ell -2}\big|.
\end{equation*}
Define the event
\begin{flalign*}
Q^{(1)}_{n,\ell,3}&=\Bigg\{\big| \beta_{0,j_{\ell+1}-1,m} -  \beta_{0,j_{\ell}-1,m}\big| \leq  \frac{\hat{w}_{\hat{j}_q +1,j_\ell}}{|j_\ell - \hat{j}_q -2|} 
  + \Big| \frac{\sqrt{m}\bar M_n(\hat{j}_q+1; j_\ell -1) }{j_\ell - \hat{j}_q}\Big| \\
&\qquad \qquad \qquad \qquad \qquad+  \frac{ \hat{w}_{j_\ell+1,\hat{j}_{q+1}}}{ |\hat{j}_{q+1} - j_\ell -2|} +\Big| \frac{ \sqrt{m}\bar M_n(j_\ell+1; \hat{j}_{q+1} -1) }{ \hat{j}_{q+1} - j_\ell -2}\Big| \Bigg\}\\
&\subset \Bigg\{ \big| \beta_{0,j_{\ell+1}-1,m} -  \beta_{0,j_{k}-1,m}\big| \leq \frac{\hat{w}_{\hat{j}_q+1,j_\ell}}{m\varepsilon_n -2} 
  + \Big| \frac{\sqrt{m}\bar M_n(\hat{j}_q+1; j_\ell -1) }{m\varepsilon_n -2}\Big| \\
&\qquad \qquad \qquad \qquad \qquad+  \frac{ \hat{w}_{j_\ell +1,\hat{j}_{q+1}}}{ m\varepsilon_n -2} +\Big| \frac{ \sqrt{m}\bar M_n(j_\ell +1; \hat{j}_{q+1} -1) }{ m\varepsilon_n -2}\Big|\Bigg\}.
\end{flalign*}
We observe that the event $Q^{(1)}_{n,\ell,3}$ occurs with probability one, so
\begin{eqnarray*}
\P[R_{n,\ell,3}^{(1)}] &=& \P[R_{n,\ell,3}^{(1)} \cap Q^{(1)}_{n,\ell,3}]\\
&\leq & \P\Big[\frac{\hat{w}_{\hat{j}_q+1,j_\ell} } {m\varepsilon_n -2} \geq  \frac{\big| \beta_{0,j_{\ell+1}-1,m} -  \beta_{0,j_{\ell}-1,m}\big| }{4}\Big]\\
&\qquad  & +  \P\Big[\frac{\hat{w}_{j_\ell+1,\hat{j}_{q+1}} } {m\varepsilon_n -2} \geq  \frac{\big| \beta_{0,j_{\ell+1}-1,m} -  \beta_{0,j_{\ell}-1,m}\big| }{4}\Big]\\
&\qquad& +  \P\Big[ \Big| \frac{\sqrt{m}\bar M_n(\hat{j}_q +1; j_\ell -1) }{m\varepsilon_n -2}\Big| \geq  \frac{\big| \beta_{0,j_{\ell+1}-1,m} -  \beta_{0,j_{\ell}-1,m}\big| }{4}\Big]\\
&\qquad& +  \P\Big[ \Big| \frac{ \sqrt{m}\bar M_n(j_\ell +1; \hat{j}_{q+1} -1) }{ m\varepsilon_n -2}\Big| \geq  \frac{\big| \beta_{0,j_{\ell+1}-1,m} -  \beta_{0,j_{\ell}-1,m}\big| }{4}\Big].
\end{eqnarray*}
Using Lemmas~\ref{lem:KKT} and \ref{lem:control-martingale},  and \eqref{ass:consistency-1}-\eqref{ass:consistency-2} from Assumption~\ref{ass:consistency}, each term of the last inequality goes to zero, as $n \rightarrow \infty$. 
 For $\P[R^{(2)}_{n,\ell,3}]$, we apply  Lemma~\ref{lem:KKT} with $j = j_\ell$ and $j=\hat{j}_q +1$ to obtain  \eqref{eq:ch1} 
 and then with $j=j_\ell +1$ and $j = j_{\ell+1}$ to get
\begin{equation*}
 \Big|\big( {j}_{\ell+1} - j_\ell - 2\big) \big\{\big( \beta_{0,j_{\ell+1}-1,m} - \hat{\beta}_{\hat{j}_{q+1} -1,m}\big) \big\} + \sqrt{m}\bar M_n(j_\ell +1; {j}_{\ell+1} -1) \Big| \leq \hat{w}_{j_\ell+1,{j}_{\ell+1}}.
\end{equation*}
It follows that event $Q_{n,\ell,3}^{(2)}$ occurs with probability one, where
\begin{flalign*}
Q_{n,\ell,3}^{(2)}&=\Bigg\{\big| \beta_{0,j_{\ell+1}-1,m} -  \beta_{0,j_{\ell}-1,m}\big| \leq  \frac{\hat{w}_{\hat{j}_q +1,j_\ell}}{|j_\ell - \hat{j}_q -2|} 
  + \Big| \frac{\sqrt{m}\bar M_n(\hat{j}_q+1; j_\ell -1) }{j_\ell - \hat{j}_q}\Big| \\
&\qquad \qquad \qquad \qquad \qquad+  \frac{ \hat{w}_{j_\ell+1,{j}_{\ell+1}}}{ |{j}_{\ell+1} - j_\ell -2|} +\Big| \frac{ \sqrt{m}\bar M_n(j_\ell+1; {j}_{\ell+1} -1) }{ {j}_{\ell+1} - j_\ell -2}\Big| \Bigg\}\\
&\subset \Bigg\{ \big| \beta_{0,j_{\ell+1}-1,m} -  \beta_{0,j_{k}-1,m}\big| \leq \frac{\hat{w}_{\hat{j}_q+1,j_\ell}}{m\varepsilon_n -2} 
  + \Big| \frac{\sqrt{m}\bar M_n(\hat{j}_q+1; j_\ell -1) }{m\varepsilon_n -2}\Big| \\
&\qquad \qquad \qquad \qquad \qquad+  \frac{ \hat{w}_{j_\ell +1,{j}_{\ell+1}}}{ \Delta_{j,\min}  -2} +\Big| \frac{ \sqrt{m}\bar M_n(j_\ell +1; \hat{j}_{q+1} -1) }{ \Delta_{j,\min}  -2}\Big|\Bigg\}.
\end{flalign*}
Then
\begin{eqnarray*}
\P[R_{n,\ell,3}^{(2)}] &=& \P[R_{n,\ell,3}^{(2)} \cap Q^{(2)}_{n,\ell,3}]\\
&\leq & \P\Big[\frac{\hat{w}_{\hat{j}_q+1,j_\ell} } {m\varepsilon_n -2} \geq  \frac{\big| \beta_{0,j_{\ell+1}-1,m} -  \beta_{0,j_{\ell}-1,m}\big| }{4}\Big]\\
&\qquad  & +  \P\Big[\frac{\hat{w}_{j_\ell+1,{j}_{\ell+1}} } {\Delta_{j,\min} -2} \geq  \frac{\big| \beta_{0,j_{\ell+1}-1,m} -  \beta_{0,j_{\ell}-1,m}\big| }{4}\Big]\\
&\qquad& +  \P\Big[ \Big| \frac{\sqrt{m}\bar M_n(\hat{j}_q +1; j_\ell -1) }{m\varepsilon_n -2}\Big| \geq  \frac{\big| \beta_{0,j_{\ell+1}-1,m} -  \beta_{0,j_{\ell}-1,m}\big| }{4}\Big]\\
&\qquad& +  \P\Big[\Big| \frac{ \sqrt{m}\bar M_n(j_\ell +1; {j}_{\ell+1} -1) }{ \Delta_{j,\min} -2}\Big| \geq  \frac{\big| \beta_{0,j_{\ell+1}-1,m} -  \beta_{0,j_{\ell}-1,m}\big| }{4}\Big].
\end{eqnarray*}
 Using Lemmas~\ref{lem:KKT} and ~\ref{lem:control-martingale},  \eqref{ass:consistency-1}-\eqref{ass:consistency-2} in Assumption~\ref{ass:consistency},  each term of the last inequality tends to zero as $ n \rightarrow + \infty.$  For $\P[R^{(3)}_{n,\ell,3}]$, we first use Lemma~\ref{lem:KKT} with $j = j_{\ell-1} +1$ and $j = j_\ell$ to get 
\begin{equation*}
\Big|\big(j_\ell - {j}_{\ell-1} - 2\big)\big( \beta_{0,j_\ell-1,m}  - \hat{\beta}_{\hat{j}_{q+1}-1,m}\big)  + \sqrt{m}\bar M_n({j}_{\ell-1}+1; j_\ell -1) \Big| \leq \hat{w}_{{j}_{\ell-1}+1,j_\ell}.
\end{equation*}
And then with $j = j_\ell +1$ and $j = \hat{j}_{q + 1}$, to obtain
\begin{equation*}
\Big|\big(\hat{j}_{q+1}- {j}_{\ell} - 2\big) \big( \beta_{0,j_{\ell+1}-1,m}  - \hat{\beta}_{\hat{j}_{q+1}-1,m}\big)  + \sqrt{m}\bar M_n({j}_{\ell}+1; \hat{j}_{q+1} -1) \Big| \leq \hat{w}_{{j}_{\ell}+1,\hat{j}_{q+1}}.
\end{equation*}
Hence the event $Q_{n,\ell,3}^{(3)}$ occurs with probability one, where
\begin{flalign*}
Q_{n,\ell,3}^{(3)}&=\Bigg\{\big| \beta_{0,j_{\ell+1}-1,m} -  \beta_{0,j_{\ell}-1,m}\big| \leq  \frac{\hat{w}_{{j}_{\ell-1} +1,j_\ell}}{|j_\ell - {j}_{\ell-1}-2|} 
  + \Big| \frac{\sqrt{m}\bar M_n({j}_{\ell-1}+1; j_\ell -1) }{j_\ell - {j}_{\ell-1} -2}\Big| \\
&\qquad \qquad \qquad \qquad \qquad+  \frac{ \hat{w}_{j_\ell+1,\hat{j}_{q+1}}}{ |\hat{j}_{q+1} - j_\ell -2|} +\Big| \frac{ \sqrt{m}\bar M_n(j_\ell+1; \hat{j}_{q+1} -1) }{ \hat{j}_{q+1}- j_\ell -2}\Big| \Bigg\}\\
&\subset \Bigg\{ \big| \beta_{0,j_{\ell+1}-1,m} -  \beta_{0,j_{k}-1,m}\big| \leq \frac{\hat{w}_{{j}_{\ell-1} +1,j_\ell}}{\Delta_{j,\min} -2} 
  + \Big| \frac{\sqrt{m}\bar M_n({j}_{\ell-1}+1; j_\ell -1) }{\Delta_{j,\min}-2}\Big| \\
&\qquad \qquad \qquad \qquad \qquad+  \frac{ \hat{w}_{j_\ell+1,\hat{j}_{q+1}}}{ m\varepsilon_n  -2} +\Big| \frac{ \sqrt{m}\bar M_n(j_\ell+1; \hat{j}_{q+1} -1)}{ m\varepsilon_n   -2}\Big|\Bigg\}.
\end{flalign*}
Then
\begin{eqnarray*}
\P[R_{n,\ell,3}^{(3)}]&=& \P[R_{n,\ell,3}^{(3)} \cap Q^{(3)}_{n,\ell,3}]\\
&\leq & \P\Big[\frac{\hat{w}_{{j}_{\ell-1} +1,j_\ell}}{\Delta_{j,\min} -2}  \geq  \frac{\big| \beta_{0,j_{\ell+1}-1,m} -  \beta_{0,j_{\ell}-1,m}\big| }{4}\Big]\\
&\qquad  & +  \P\Big[\frac{ \hat{w}_{j_\ell+1,\hat{j}_{q+1}}}{ m\varepsilon_n  -2} \geq  \frac{\big| \beta_{0,j_{\ell+1}-1,m} -  \beta_{0,j_{\ell}-1,m}\big| }{4}\Big]\\
&\qquad& +  \P\Big[\Big| \frac{\sqrt{m}\bar M_n({j}_{\ell-1}+1; j_\ell -1) }{\Delta_{j,\min}-2}\Big|\geq  \frac{\big| \beta_{0,j_{\ell+1}-1,m} -  \beta_{0,j_{\ell}-1,m}\big| }{4}\Big]\\
&\qquad& +  \P\Big[\Big| \frac{ \sqrt{m}\bar M_n(j_\ell+1; \hat{j}_{q+1} -1)}{ m\varepsilon_n   -2}\Big| \geq  \frac{\big| \beta_{0,j_{\ell+1}-1,m} -  \beta_{0,j_{\ell}-1,m}\big| }{4}\Big].
\end{eqnarray*}
By Lemmas~\ref{lem:KKT} and \ref{lem:control-martingale}, and \eqref{ass:consistency-1}-\eqref{ass:consistency-2} in Assumption~\ref{ass:consistency}, it  implies  that each term of the last inequality tends to zero as $ n \rightarrow + \infty.$.  Finally, for $\P[R^{(4)}_{n,\ell,3}],$ we first use Lemma~\ref{lem:KKT} with $j = j_{\ell-1} +1$ and $j = j_\ell$ to obtain 
\begin{equation*}
\Big|\big(j_\ell - {j}_{\ell-1} - 2\big) \big\{\big( \beta_{0,j_\ell-1,m}  - \hat{\beta}_{\hat{j}_{q+1}-1,m}\big) \big\} + \sqrt{m}\bar M_n({j}_{\ell-1}+1; j_\ell -1) \Big| \leq \hat{w}_{{j}_{\ell-1}+1,j_\ell}.
\end{equation*}
Second, we use Lemma~\ref{lem:KKT} with $j = j_\ell +1$ and $j = j_{\ell+1}$ to obtain 
\begin{equation*}
\Big|\big({j}_{\ell+1}- {j}_{\ell} - 2\big) \big( \beta_{0,j_{\ell+1}-1,m}  - \hat{\beta}_{\hat{j}_{q+1}-1,m}\big) + \sqrt{m}\bar M_n({j}_{\ell}+1; {j}_{\ell+1} -1) \Big| \leq \hat{w}_{{j}_{\ell}+1,{j}_{\ell+1}}.
\end{equation*}
It follows that the event $Q_{n,\ell,3}^{(4)}$ occurs with probability one, where\\

\begin{flalign*}
Q_{n,\ell,3}^{(4)}&=\Bigg\{\big| \beta_{0,j_{\ell+1}-1,m} -  \beta_{0,j_{\ell}-1,m}\big| \leq  \frac{\hat{w}_{{j}_{\ell-1} +1,j_\ell}}{|j_\ell - {j}_{\ell-1}-2|} 
  + \Big| \frac{\sqrt{m}\bar M_n({j}_{\ell-1}+1; j_\ell -1) }{j_\ell - {j}_{\ell-1} -2}\Big| \\
&\hspace{4.7cm}+  \frac{ \hat{w}_{j_\ell+1,{j}_{\ell+1}}}{ |j_{\ell+1} - j_\ell -2|} +\Big| \frac{ \sqrt{m}\bar M_n(j_\ell+1; {j}_{\ell+1} -1) }{ {j}_{\ell+1}- j_\ell -2}\Big| \Bigg\}\\
&\subset \Bigg\{ \big| \beta_{0,j_{\ell+1}-1,m} -  \beta_{0,j_{k}-1,m}\big| \leq \frac{\hat{w}_{{j}_{\ell-1} +1,j_\ell}}{\Delta_{j,\min} -2} 
  + \Big| \frac{\sqrt{m}\bar M_n({j}_{\ell-1}+1; j_\ell -1) }{\Delta_{j,\min}-2}\Big| \\
&\hspace{5.3cm}+  \frac{ \hat{w}_{j_\ell+1,{j}_{\ell+1}}}{ \Delta_{j,\min} -2} +\Big| \frac{ \sqrt{m}\bar M_n(j_\ell+1; {j}_{\ell+1} -1)}{ \Delta_{j,\min} -2}\Big|\Bigg\}.
\end{flalign*}
Then
\begin{eqnarray*}
\P[R_{n,\ell,3}^{(4)}]&=& \P[R_{n,\ell,3}^{(4)} \cap Q^{(4)}_{n,\ell,3}]\\
&\leq & \P\Bigg[ \frac{\hat{w}_{{j}_{\ell-1} +1,j_\ell}}{\Delta_{j,\min} -2}  \geq  \frac{\big| \beta_{0,j_{\ell+1}-1,m} -  \beta_{0,j_{\ell}-1,m}\big| }{4}\Bigg]\\
&\qquad  & +  \P\Bigg[ \frac{ \hat{w}_{j_\ell+1,{j}_{\ell+1}}}{ \Delta_{j,\min} -2} \geq  \frac{\big| \beta_{0,j_{\ell+1}-1,m} -  \beta_{0,j_{\ell}-1,m}\big| }{4}\Bigg]\\
&\qquad& +  \P\Bigg[\bigg| \frac{\sqrt{m}\bar M_n({j}_{\ell-1}+1; j_\ell -1) }{\Delta_{j,\min}-2}\bigg|\geq  \frac{\big| \beta_{0,j_{\ell+1}-1,m} -  \beta_{0,j_{\ell}-1,m}\big| }{4}\Bigg]\\
&\qquad& +  \P\Bigg[\bigg| \frac{ \sqrt{m}\bar M_n(j_\ell+1; {j}_{\ell+1} -1)}{ \Delta_{j,\min} -2}\bigg| \geq  \frac{\big| \beta_{0,j_{\ell+1}-1,m} -  \beta_{0,j_{\ell}-1,m}\big| }{4}\Bigg]\\
&\rightarrow& 0.
\end{eqnarray*}
as $n\rightarrow \infty.$
This concludes the proof of Theorem~\ref{thm4}.
$\hfill \square$

\begin{appendices}
\section{ }

\renewcommand{\theequation}{\thesection.\arabic{equation}}

Here we prove Proposition~\ref{prop1} and Lemmas
\ref{lem:approximation}, \ref{lem:KKT} and \ref
{lem:control-martingale}

\setcounter{equation}{0}
\subsection{ Proof of Proposition~\ref{prop1}}
\label{app:proof-prop1}

Fix $ j \in \{1, \ldots, m\}.$  We have
\begin{equation*}
  \begin{split}
&U_j= \frac{1}{n}\sum_{i=1}^n  \int_0^1 \ind{(\frac{j-1}{m},1]} (s)dM_i(s),\\
&V_j= n \inr{U_j}= \frac{1}{n}  \int_0^1 \ind{(\frac{j-1}{m},1]} (s)\lambda_0(s)d(s),\\
&\hat{V}_j=n [U_j]= \frac{1}{n}\sum_{i=1}^n  \int_0^1 \ind{(\frac{j-1}{m},1]} (s)dN_i(s).\\
\end{split}
\end{equation*}
 Classical  Bernstein deviation inequality applied to $U_j$, see \hypertarget{}{\cite{Van-95}},  yields that 
\begin{equation}
  \label{bern}
\begin{aligned}
  \P\Big[|U_j| \geq \sqrt{2\theta z} + \frac{z}{3n}, \frac{1}{n}  \int_0^1 \ind{(\frac{j-1}{m},1]} (s)\lambda_0(s)d(s) \leq \theta\Big] \leq 2e^{-z}.
\end{aligned}
\end{equation}
for all $\theta > 0,$ and $z> 0.$ It follows that
\begin{equation}
\label{bernbis}
 \P\Big[|U_j| \geq \sqrt{\frac{2\theta z}{n}} + \frac{z}{3n}, V_j \leq \theta \Big] \leq 2e^{-z}.
\end{equation}
By choosing $\theta = c_0(z+1)/n$, this gives 
\begin{equation}
  \label{bersubtet}
   \P\Big[|U_j| \geq \Big(\sqrt{2c_0} + \frac{1}{3}\Big)\frac{z+1}{n}, V_j\leq \frac{c_0(z+1)}{n} \Big] \leq 2e^{-z}.
\end{equation}
For any $0 < \eta < \theta < \infty$, we have
\begin{equation*}
\Big\{| U_j| \geq \sqrt{\frac{2\theta V_j z}{\eta n}} + \frac{z}{3n}\Big\}\cap\Big\{\eta< V_j \leq \theta\Big\} \subset \Big\{| U_j| \geq \sqrt{\frac{2\theta  z}{ n}} + \frac{z}{3n}\Big\}\cap\Big\{\eta < V_j\leq\theta\Big\}.
 \end{equation*}
Together with \eqref{bernbis}, we obtain
\begin{equation}
  \label{bernvj}
 \P\Big[|U_j| \geq \sqrt{\frac{2\theta V_j z}{\eta n}} + \frac{z}{3n},  \eta < V_j \leq \theta \Big] \leq 2e^{-z}.
\end{equation}
Now we want to replace $V_j$ by the observable $\hat{V}_j$ in the deviation (\ref{bernbis}). Define  $\widetilde{U}_j$  by
\begin{eqnarray*}
\widetilde{U}_j &=& \hat{V}_j - V_j\\
&=&  \frac{1}{n}\sum_{i=1}^n  \int_0^1\ind{(\frac{j-1}{m},1]} (s)dM_i(s).
\end{eqnarray*}
Now writing again  (\ref{bernvj}) and using the same argument as before, we arrive at
\begin{equation}
  \label{berdiffvj}
   \P\Big[|\widetilde{U}_j| \geq \sqrt{\frac{2\theta V_j z}{\eta n}} + \frac{z}{3n},  \eta < V_j \leq \theta \Big] \leq 2e^{-z}.
\end{equation}
But, if $V_j$   satisfies
\begin{equation*}
|\widetilde{U}_j| \leq \sqrt{\frac{2\theta V_j z}{\eta n}} + \frac{z}{3n},
\end{equation*}
then it satisfies
\begin{equation*}
V_j \leq 2\hat{V}_j + 2\big( \frac{\theta}{\eta} + \frac{1}{3}\big)\frac{z}{n},
\end{equation*}
and
\begin{equation*}
\hat{V}_j \leq 2V_j + 2\Big(\frac{1}{3} + 2\sqrt{\frac{\theta}{\eta}\big( \frac{\theta}{\eta} + \frac{1}{3}\big)} + 2\frac{\theta}{\eta} \Big)\frac{z}{n},
\end{equation*}
simply using the fact that  $A \leq b + \sqrt{aA}$ entails $A \leq a + 2b$ for any $a, A, b > 0$. This proves that
\begin{equation}
  \label{absuj}
\begin{split}
  &\Big\{| U_j| \leq \sqrt{\frac{2\theta V_j z}{\eta n}} + \frac{z}{3n}\Big\}\cap\Big\{ \Big|\widetilde{U}_j\Big| \leq \sqrt{\frac{2\theta V_j z}{\eta n}} + \frac{z}{3n}\Big\}\\
&\qquad\subset \Big\{| U_j| \leq 2\sqrt{\frac{\theta z}{\eta n} \hat{V}_j} +  \Big(\frac{1}{3} + 2\sqrt{\frac{\theta}{\eta}\big( \frac{\theta}{\eta} + \frac{1}{3}\big)} + 2\frac{\theta}{\eta} \Big)\frac{z}{n}\Big\}.
\end{split}
\end{equation}
So, using (\ref{bernvj}) and (\ref{berdiffvj}), we obtain
\begin{equation*}
\P\Big[| U_j| \geq 2\sqrt{\frac{\theta z}{ n} \hat{V}_j} +  \Big(\frac{1}{3} + 2\sqrt{\frac{\theta}{\eta}\big( \frac{\theta}{\eta} + \frac{1}{3}\big)}  \Big)\frac{z}{n}, \eta< V_j \leq \theta \Big] \leq 4e^{-z}.
\end{equation*}
The inequality is similar to (\ref{bernvj}), where we replaced $V_j$ by the observable $\hat{V}_j$. It remains to remove the event $\big\{\eta< V_j \leq \theta \big\}$ from this inequality. First, recall that (\ref{bersubtet}) holds, so we can work on the event $\big\{V_j > c_0(z+1)/n \big\}$ from now on. We use a peeling argument. Define, for $q \geq 0$:
\begin{equation*}
 \theta_q = c_0\frac{(z+1)}{n}(1+\varepsilon)^q,
\end{equation*}
and use the following decomposition into disjoint sets:
\begin{equation*}
\{V_j> \theta_0\} = \bigcup_{q\geq 0}\big\{\theta_q < V_j\leq \theta_{q+1}\big\}.
\end{equation*}
We have
\begin{equation*}
\P\Big[|U_j|\geq c_{1,\varepsilon} \sqrt{\frac{z}{n}\hat{V}_j} + c_{2,\varepsilon}\frac{z}{n}, \theta_q < V_j \leq \theta_{q+1}\Big] \leq 4e^{-z}
\end{equation*}
where
\begin{equation*}
c_{1,\varepsilon}= 2\sqrt{1+\varepsilon}\quad  \textrm{and} \quad c_{2,\varepsilon}=  2\sqrt{ (1+\varepsilon)\big(\frac{4}{3}+\varepsilon\big)} + \frac{1}{3}.
\end{equation*}
Let 
\begin{equation*}
h_j= c_h \log\log\big(\frac{V_j}{\theta_0}\vee e\big).
\end{equation*}
 On the event
\begin{equation*}
\Big\{|\widetilde{U}_j| \leq \sqrt{\frac{2(1+\varepsilon)V_j(z+h_j)}{n}}+ \frac{(z+h_j)}{3n}\Big\}
\end{equation*}
we have
\begin{equation*}
V_j \leq 2\hat{V}_j + 2(\frac{4}{3}+\varepsilon)\frac{z}{n} + 2\frac{\frac{4}{3}+\varepsilon}{n}c_h \log\log\big(\frac{V_j}{\theta_0}\vee e\big),
\end{equation*}
which entails, assuming that $ec_0 > 2((1+\varepsilon) +\frac{1}{3})c_h,$
\begin{equation*}
V_j \leq \frac{ec_0(z+1)}{ec_0(z+1) - 2(\frac{4}{3}+\varepsilon)c_h}\Big(2\hat{V}_j + 2(\frac{4}{3}+\varepsilon)\frac{z}{n}\Big),
\end{equation*}
where we used the fact that $\log\log z \leq z/e -1$ for any $z\geq e$. This entails, together with (\ref{absuj}), the following embedding:
\begin{equation*}
\begin{split}
&\Big\{| U_j| \leq \sqrt{\frac{2(1+\varepsilon)(z+ h_j)V_j}{n}}+ \frac{z+h_j}{3n}\Big\}\cap\Big\{|\widetilde{U}_j| \leq \sqrt{\frac{2(1+\varepsilon)(z+h_j)V_j}{n}}+ \frac{(z+h_j)}{3n}\Big\}\\
&\qquad \subset \Big\{ |U_j|\geq c_{1,\varepsilon} \sqrt{\frac{z+\hat{h}_{n,z,j}}{n}\hat{V}_j} + c_{2,\varepsilon}\frac{z+\hat{h}_{n,z,j}}{n}\Big\},
\end{split}
\end{equation*}
where
\begin{equation*}
\hat{h}_{n,z,j}= c_h\log\log\Big(\frac{2en\hat{V}_j+ 2e(\frac{4}{3}+\varepsilon)z}{ec_0(z+1) - 2(\frac{4}{3}+\varepsilon )c_h}\vee e \Big).
\end{equation*}
Now, using the previous embeddings together with (\ref{bernvj})  and \eqref{berdiffvj} we obtain
\begin{flalign*}
&\P\Big[|U_j|\geq c_{1,\varepsilon} \sqrt{\frac{z+\hat{h}_{n,z,j}}{n}\hat{V}_j} + c_{2,\varepsilon}\frac{z+\hat{h}_{n,z,j}}{n}, V_j > \theta_0\Big] \\
&\qquad \leq  \sum_{q\geq0}\P\Big[|U_j| \geq \sqrt{\frac{2(1+\varepsilon)V_j(z+h_j)}{n}}+ \frac{z+h_j}{3n}, \theta_q < V_j \leq \theta _{q+1}\Big]\\
&\qquad \qquad + \sum_{q\geq 0}\P\Big[|\widetilde{U}_j| \geq \sqrt{\frac{2(1+\varepsilon)V_j(z+h_j)}{n}}+ \frac{z+h_j}{3n}, \theta_q < V_j \leq \theta _{q+1} \Big]\\
&\qquad \leq   4\big(e^{-z} + \sum_{q\geq 1}e^{-\big(z+c_h\log\log(\frac{V_j}{\theta_0})\big)}\big)\\
&\qquad =  4\big(1 + \big(\log(1+\varepsilon)\big)^{-c_h} \sum_{q\geq 1}q^{-c_h}\big)e^{-z}.
\end{flalign*}
Then with \eqref{bersubtet}, it implies that
\begin{equation*}
\P\Big[ |U_j|\geq c_{1,\varepsilon} \sqrt{\frac{z+\hat{h}_{n,z,j}}{n}\hat{V}_j} + c_{3,\varepsilon}\frac{z+1+\hat{h}_{n,z,j}}{n}\Big] \leq \big(6 + 4\big(\log(1+\varepsilon)\big)^{-c_h} \sum_{q\geq 1}q^{-c_h}\big)e^{-z},
\end{equation*}
where $c_{3,\varepsilon}= \sqrt{2\max\big(c_0, 2(1+\varepsilon)(\frac{4}{3}+\varepsilon)\big)} + \frac{1}{3}.$
$\hfill \square$

\subsection{Proof of Lemma~\ref{lem:approximation}}
\label{app:proof-lem:approximation}

Using the fact that the functions $\left\{\lambda_{j,m}: j=1, \ldots, m\right\}$ form a basis of $\Lambda_m$, and under Assumption~\ref{ass:intensity},  one can give the explicit form of $\lambda_{0,m}$ as following
\begin{equation}
  \lambda_{0,m} = m \sum_{j=1}^m \sum_{\ell=1}^{L_0} \beta_{0,\ell} |J_\ell\cap I_{j,m}|\ind{I_{j,m}} =  \sum_{j=1}^m \beta_{0,j,m}\lambda_{j,m},
\end{equation}
here $|A|$ is the Lebesgue measure of the set $A$ and  $ \beta_{0,j,m} = \sqrt{m} \sum_{\ell=1}^{L_0} \beta_{0,\ell} |J_\ell\cap I_{j,m}|$.
We remark  that the intervals $J_\ell$ do not share the same
boundaries as the smaller intervals $I_{j, m}$.
Setting the sequence $(\bar{\ell}_j)_{j=0,\ldots,m}$ be the sequence defining by:
\begin{equation*}
\bar{\ell}_0 = 1,\textrm{ and } 
\bar{\ell}_j = \max\{\ell=1, \ldots, L_0: J_\ell \cap I_{j,m} \neq \emptyset\}, \textrm{ for } j=1, \ldots, m.
\end{equation*}
Using the sequence  $(\bar{\ell}_j)_{j=0,\ldots,m}$, one  has the expression of the functions $\lambda_0$ and  $\lambda_{0,m}$ as follows:
 \begin{equation*}
 \lambda_0 =  \sum_{j=1}^m \sum_{\ell=\bar{\ell}_{j-1}}^{\bar{\ell}_j} \beta_{0,\ell} \ind{{J_\ell}\cap I_{j,m}},
 \end{equation*}
and
\begin{equation*}
\begin{split}
 \lambda_{0,m} &=  m \sum_{j=1}^m \sum_{\ell=\bar{\ell}_{j-1}}^{\bar{\ell}_j} \beta_{0,\ell} |J_\ell\cap I_{j,m}|\ind{ I_{j,m}}
= \sum_{j=1}^m \alpha_{0,j,m}\ind{ I_{j,m}} 
= \sum_{j=1}^m \alpha_{0,j,m}\sum_{\ell=\bar{\ell}_{j-1}}^{\bar{\ell}_j}\ind{J_\ell\cap I_{j,m}},
\end{split}
\end{equation*}
where $ \alpha_{0,j,m}= m\sum_{\ell=\bar{\ell}_{j-1}}^{\bar{\ell}_j} \beta_{0,\ell} |J_\ell\cap I_{j,m}|.$
From the fact that  $\big\{ \ind{J_\ell\cap I_{j,m}}: j=1, \ldots, m\textrm{ and }\ell= 1, \ldots, L_0\big\}$ is an orthogonal basis of $\Lambda_m$ (with respect to the  $\mathbb{L}^2$-norm), we  obtain
\begin{equation*}
\begin{split}
\norm{\lambda_0 - \lambda_{0,m} }^2 &= \Bigg\|\sum_{j=1}^m \sum_{\ell=\bar{\ell}_{j-1}}^{\bar{\ell}_j} \Big( \beta_{0,\ell} - m\sum_{\ell'=\bar{\ell}_{j-1}}^{\bar{\ell}_j} \beta_{0,\ell'} |J_{\ell'}\cap I_{j,m}|\Big)\ind{J_\ell\cap I_{j,m}}\Bigg\|^2\\
 &= \sum_{j=1}^m \sum_{\ell=\bar{\ell}_{j-1}}^{\bar{\ell}_j}\Big( \beta_{0,\ell} - m\sum_{\ell'=\bar{\ell}_{j-1}}^{\bar{\ell}_j} \beta_{0,\ell'} |J_{\ell'}\cap I_{j,m}|\Big)^2 |J_\ell\cap I_{j,m}|\\
 &= \sum_{j=1}^m \ind{[\bar{\ell}_j - \bar{\ell}_{j-1} > 0]}\sum_{\ell=\bar{\ell}_{j-1}}^{\bar{\ell}_j}\Big( \beta_{0,\ell} - m\sum_{\ell'=\bar{\ell}_{j-1}}^{\bar{\ell}_j} \beta_{0,\ell'} |J_{\ell'}\cap I_{j,m}|\Big)^2 |J_\ell\cap I_{j,m}|\\
&\leq \sum_{j=1}^m \ind{[\bar{\ell}_j - \bar{\ell}_{j-1} > 0]}\sum_{\ell=\bar{\ell}_{j-1}}^{\bar{\ell}_j} (\bar{\ell}_j - \bar{\ell}_{j-1} +1) \max_{\ell, \ell' = k_{j-1}, \ldots, k_j} \Big( \beta_{0,\ell} -  \beta_{0,\ell'}\Big)^2|J_\ell\cap I_{j,m}| \\
&\leq  \sum_{j=1}^m\ind{[\bar{\ell}_j - \bar{\ell}_{j-1} > 0]}  (\bar{\ell}_j - \bar{\ell}_{j-1} +1) \max_{\ell,\ell'\in\{ \bar{\ell}_{j-1}, \ldots, \bar{\ell}_j\}} \Big( \beta_{0,\ell} -  \beta_{0,\ell'}\Big)^2| I_{j,m}|\\
&\leq\frac{{2(L_0-1)\Delta_{\beta,\max}^2}}{m}.
\end{split}
\end{equation*}
This proves Lemma~\ref{lem:approximation}.
$\hfill \square$

\subsection{Proof of Lemma~\ref{lem:KKT}}
\label{app:proof-lem:KKT}

To prove Lemma~\ref{lem:KKT}, we invoke subdifferential calculus, see \cite{Ber-99}. We first write our objective functional as  
\begin{equation*}\Phi(\mu) = \frac{1}{2}\Sum_{j =1}^m (\bN_j - ({\bT}\mu)_j)^2 +  \sum_{j =1}^m \hat{w}_j|\mu_{j,m}|.\end{equation*}
So a necessary and sufficient condition for a vector $\hat{ {\mu}}$ in $\R^m$ to minimize the function $\Phi$ 
is that the zero vector in $\R^m$ belongs to the sub-differential of $\Phi(\mu)$ at the point $\hat{ {\mu}}$, that is, the following  optimality condition holds:
\begin{equation*}
\textrm{for all } j =1,\ldots, m \left\{
  \begin{array}{ll}
    \Big({\bT}^\top\big( \bN- {\bT}\hat{ {\mu}}\big)\Big)_j = {\hat{w}_j}\textrm{sign}(\hat{\mu}_{j,m}),\textrm{ if }\hat{\mu}_{j,m} \neq 0, \\
    \Big|\Big( {\bT}^\top\big(\bN - {\bT}\hat{ {\mu}}\big) \Big)_j \Big| \leq {\hat{w}_j}\textrm{sign}(\hat{\mu}_{j,m}),\textrm{ if }\hat{\mu}_{j,m} = 0.
  \end{array}
\right.
\end{equation*}
Using that $({\bT}^\top \bN)_j = \sum_{q=j}^m \bN_q$ and that $({\bT}^\top\hat{ {\beta}})_j = \sum_{q=j}^m\hat{\beta}_{q,m},$ since ${\bT}$ is a $m \times m$ lower triangular matrix having all its nonzero elements equal to one.
\noindent Now,  for $q=1, \ldots, m$, we observe that
\begin{equation*}
\begin{split}
\bN_q &= \sqrt{m}\int_{I_{q,m}}\lambda_0(t)dt +  \sqrt{m}\bar M_n(I_{q,m})\\
&= \sqrt{m}\int_{I_{q,m}}(\lambda_0(t) - \lambda_{0,m}(t))dt + \sqrt{m}\int_{I_{q,m}} \lambda_{0,m}dt +  \sqrt{m}\bar M_n(I_{q,m})\\
&= \sqrt{m}\int_{I_{q,m}}\Big(\sum_{\ell=1}^{L_0}\beta_{0,\ell}\ind{J_\ell} - \sqrt{m}\sum_{j=1}^m\beta_{0,j,m}\ind{I_{j,m}}\Big)dt\\
& \qquad \quad \,\, \,\, \,+  \sqrt{m}\sum_{j=1}^m\beta_{0,j,m}\int_{I_{q,m}}\ind{I_{j,m}}(t)dt + \sqrt{m}\bar M_n(I_{q,m})\\
&= \beta_{0,{q,m}} + \sqrt{m}\bar M_n(I_{q,m}),
\end{split}
\end{equation*}
and we get the desired result.

$\hfill \square$

\subsection{Proof of Lemma~\ref{lem:control-martingale}}
\label{app:proof-lem-control-martingale}

For the first statement, we have  by definition,
\begin{equation*}
\begin{split}
\Big|\bar M_n(a;b) \Big|  &= \Big|\sum_{q=a}^{b}  \bar M_n (I_{q,m})\Big|  \\
&=  \Big|\frac{1}{n}\sum_{i=1}^n\sum_{q=a}^{b}  \int_0^1{I_{q,m}}(t)dM_i(t)\Big| \\
&=  \Big|\frac{1}{n}\sum_{i=1}^n  \int_0^1{\textbf{1}_{(\frac{a-1}{m},\frac{b}{m}]}}(t)dM_i(t)\Big|.
\end{split}
\end{equation*} 
Moreover, using  Bernstein's  inequality, it follows that, for any $z, \alpha >0,$ 
 \begin{equation*}
\begin{split}
&\P\Big[\Big|\frac{1}{n}\sum_{i=1}^n  \int_0^1{\textbf{1}_{(\frac{a-1}{m},\frac{b}{m}]}}(t)dM_i(t)\Big| \geq z,\\
& \qquad \qquad  \qquad \quad \inr{\frac{1}{n}\sum_{i=1}^n  \int_0^1{\textbf{1}_{(\frac{a-1}{m},\frac{b}{m}]}}(t)dM_i(t)} \leq  \alpha\Big]
   \leq 2\exp\Big\{-\frac{z^2}{2\alpha + \frac{2}{3}\rho z }\Big\}, 
\end{split}
\end{equation*}
where $\rho$ is a upper bound of $\Big\| \frac{1}{n}\textbf{1}_{(\frac{a-1}{m},\frac{b}{m}]}\Big\|_\infty$. Here we can choose $\rho = \frac{1}{n}$ and
\begin{equation*}
\alpha = n^{-1}{\int_{\frac{a-1}{m}}^{\frac{b}{m}}\lambda_0(t) dt}  = n^{-1}{ \E\Big[\bar N_n\Big(\big(\frac{a-1}{m},\frac{b}{m}\big]\Big)\Big]}.
\end{equation*}
Hence, we obtain the  first statement. For the second one,  recall that for any $a=2, \ldots,m$  we have 
\begin{equation*}
\hat{w}_a= c_1  \sqrt{\frac{m(x + \log m + \hat{h}_{n,x,a})\hat{V}_a}{n}} + c_2\frac{\sqrt{m}(x+1+ \log m +\hat{h}_{n,x,a})}{n}.
\end{equation*}
Since each term of $\hat{w}_a$ is positive and taking in account the
dominant one, we have
\begin{equation*}
  \bigg\{\hat{w}_a^2 \geq \frac{m\log m}{n}\Big(\xi-
  \int_{\ind{(\frac{a-1}{m},1]}} \lambda_0(t)dt\Big) \bigg\} \subset
  \bigg\{ \hat{V}_a\geq \xi- \int_{\ind{(\frac{a-1}{m},1]}}
  \lambda_0(t)dt\bigg\},
\end{equation*}
for all $\xi >0$.
By the Doob-Meyer decomposition theorem, we get
\begin{equation*}
  \bigg\{\hat{w}_a^2 \geq \frac{m\log m}{n}\Big(\xi-
  \int_{\ind{(\frac{a-1}{m},1]}} \lambda_0(t)dt\Big) \bigg\} \subset
  \big\{ \bar M_n(a;1)\geq \xi \big\},
\end{equation*}
Finally, by applying the first statement, see \eqref{ctrlmart}, to $\bar M_n(a;1),$  we concludes the proof of Lemma~\ref{lem:control-martingale}. 
 $\hfill \square$

\section{}
\label{app:proof-thm3-case2}
\setcounter{equation}{0}

Here we prove the second case of the proof of Theorem~\ref{thm3} which is quite similar to the first one  with a careful choice of the bounded terms in the approximate change-points sequence while applying the KKT optimality conditions.
As $m\varepsilon_n \geq 6$ for all $n \geq 1$, it  yields
that the event $ \big\{ j_{\ell} +2 < \hat{j}_\ell\big\}$ a.s.

\subsection{Step II.1. Prove: $ \P[A_{n,\ell} \cap C_n] \rightarrow 0,$  as  $ n\rightarrow \infty.$ }
\label{subsection1:CASE-II}

 Applying $(\ref{kkt})$ in Lemma~\ref{lem:KKT}  with $j = j_\ell+1$  and $j=\hat{j}_\ell,$  we get
\begin{equation*}
-{\hat{w}_{j_\ell}}\leq \sum_{q=j_\ell+1}^m \bN_q - \sum_{q= j_\ell}^m\hat{\beta}_{q,m}  \leq \hat{w}_{j_\ell+1},
\end{equation*}
and
\begin{equation*}
 -{\hat{w}_{\hat{j}_\ell}}\leq \sum_{q=\hat{j}_\ell }^m \bN_q - \sum_{q=\hat{j}_\ell}^m\hat{\beta}_{q,m} \leq {\hat{w}_{\hat{j}_\ell}}.
\end{equation*}
It follows that
\begin{equation*}
\Big| \sum_{q= {j}_\ell +1}^{\hat{j}_\ell -1}  {\beta_{0,q,m}} + \sqrt{m}\bar M_n(I_{q,m}) - \sum_{q={j}_\ell + 1 }^{\hat{j}_\ell -1}\hat{\beta}_{q,m}\Big| \leq  \hat{w}_{j_\ell+1,\hat{j}_\ell} .
\end{equation*}
The property of the vector $\hat{ {\beta}}$ in Lemma~\ref{lem:KKT}  yields that
\begin{equation*}
\Big|( \hat{j}_\ell - j_\ell - 2) (\beta_{0,j_{\ell+1}-1,m} -  \hat{\beta}_{\hat{j}_{\ell}-1,m})  + \sqrt{m}\bar M_n({j}_\ell+1;\hat{j}_\ell-1) \Big|\leq  \hat{w}_{j_\ell+1,\hat{j}_\ell} .
\end{equation*}
Therefore, on $C_n \cap \{\hat{j}_\ell > j_\ell +2\}$ we have
\begin{equation*}
\begin{split}
&\Big|( \hat{j}_\ell - j_\ell -2)(\beta_{0,j_{\ell}-1,m} -  \hat{\beta}_{\hat{j}_{\ell}-1,m})\\
 & \hspace{4cm}+ \sqrt{m}\bar M_n({j}_\ell+1;\hat{j}_\ell-1)\\
&\hspace{5cm} +(\hat{j}_\ell - j_\ell-2 )( \beta_{0,j_{\ell+1}-1,m}- \beta_{0,j_\ell-1,m})\Big|\leq \hat{w}_{j_\ell+1,\hat{j}_\ell}.
\end{split}
\end{equation*}
\noindent Define the event
\begin{equation*}
\begin{split}
C'_{n,\ell}&= \bigg\{\Big|( \hat{j}_\ell - j_\ell-2)(\beta_{0,j_{\ell}-1,m} -  \hat{\beta}_{\hat{j}_{\ell}-1,m}) \\
& \hspace{3.7cm} + \sqrt{m}\bar M_n({j}_\ell+1;\hat{j}_\ell-1)\\
& \hspace{4cm}+  (\hat{j}_\ell - j_\ell -2 )( \beta_{0,j_{\ell+1}-1,m}- \beta_{0,j_\ell-1,m})\Big|\leq \hat{w}_{j_\ell+1,\hat{j}_\ell}\bigg\}.
\end{split}
\end{equation*}
It follows that $ C'_{n,\ell}$ occurs with probability one. 
We observe  that, $m\varepsilon_n \geq 6$  for all $n$, entails that $\frac{m\varepsilon_n}{2} - 2 \geq \frac{m\varepsilon_n}{6}$. Then
\begin{equation*}
\big\{ |\hat{j}_\ell - j_\ell| >\frac{m\varepsilon_n}{2} \big\} \subset \big\{ |\hat{j}_\ell - j_\ell - 2| >\frac{m\varepsilon_n}{2} - 2 \big\} \subset \big\{ |\hat{j}_\ell - j_\ell - 2| \geq \frac{m\varepsilon_n}{6} \big\}.
\end{equation*}
\noindent Therefore
\begin{equation*}
\begin{split}
  &\P[A_{n,\ell}\cap C_n] \\
&= \P\bigg[\Big\{\Big|( \hat{j}_\ell - j_\ell-2)(\beta_{0,j_{\ell}-1,m} -  \hat{\beta}_{\hat{j}_{\ell}-1,m})  + \sqrt{m}\bar M_n({j}_\ell+1;\hat{j}_\ell-1)\\
   &   \hspace{3.5cm}+(\hat{j}_\ell - j_\ell -2 )( \beta_{0,j_{\ell+1}-1,m}- \beta_{0,j_\ell-1,m})\Big|\leq \hat{w}_{j_\ell+1,\hat{j}_\ell}\Big\}\\
& \hspace{8cm}\cap A_{n,\ell}\cap C_n \cap \{\hat{j}_\ell > j_\ell +2\}\bigg] \\
   &\leq \P\Big[\Big\{\frac{ \hat{w}_{j_\ell+1,\hat{j}_\ell}}{\frac{m\varepsilon_n}{6} } \geq \frac{| \beta_{0,j_{\ell+1}-1,m}- \beta_{0,j_\ell-1,m}|}{3}\Big\} \cap \Big\{\hat{j}_\ell > j_\ell +2\Big\}\Big] \\
   &  \qquad+ \P\Big[\Big\{ | \hat{\beta}_{\hat{j}_{\ell}-1,m}-\beta_{0,j_{\ell}-1,m}| \geq \frac{| \beta_{0,j_{\ell+1}-1,m}- \beta_{0,j_\ell-1,m}|}{3}\Big\}\cap C_n\Big] \\
    &   \qquad+ \P\Big[ \Big| \frac{\sqrt{m}\bar M_n({j}_\ell;\hat{j}_\ell-1)}{\hat{j}_\ell - {j}_\ell -2}\Big|\geq \frac{| \beta_{0,j_{\ell+1}-1,m}- \beta_{0,j_\ell-1,m}|}{3} \Big] \\
    &:= \P[A'_{n,\ell,1}] + \P[A'_{n,\ell,2}] + \P[A'_{n,\ell,3}].
\end{split}
\end{equation*}
We have 
\begin{eqnarray*}
\P[A'_{n,\ell,1}]&\leq& \P\Big[\hat{w}_{j_\ell+1,\hat{j}_\ell} \geq \frac{m\varepsilon_n \Delta_{\beta,\min}}{18}\Big]\\
 &\leq& \P\Big[\hat{w}_{{j}_\ell +1} \geq \frac{m\varepsilon_n \Delta_{\beta,\min}}{36}\Big]\\
 &=& \P\Big[\hat{w}^2_{{j}_{\ell} +1} \geq \frac{m^2\varepsilon_n^2 \Delta_{\beta,\min}^2}{36^2}\Big].
\end{eqnarray*}
 By \eqref{ass:consistency-1} in Assumption~\ref{ass:consistency}, and \eqref{ctrlweigh} in Lemma~\ref{lem:control-martingale} with $\xi = \frac{nm\varepsilon_n^2\Delta_{\beta,\min}^2}{36^2\log m} + \E\big[\bar N_n\big((\frac{j_{\ell}}{m},1]\big)\big]$, it follows that 
\begin{eqnarray*}
\P(A'_{n,\ell,1}) &\leq&  2\exp\Bigg(-\frac{n\xi^2}{2 \E\Big[\bar N_n\Big(\big(\frac{j_{\ell}}{m},1\big]\Big)\Big] +\frac{2}{3}\xi }\Bigg) \rightarrow 0,
\end{eqnarray*}
as $n \rightarrow \infty.$
\noindent Now, we remark that 
\begin{equation*}
 A'_{n,\ell,3}\subset\Big\{\Big| \bar M_n({j}_\ell;\hat{j}_\ell-1)\Big|\geq \frac{ m\varepsilon_n\Delta_{\beta,\min} }{18\sqrt{m}}\Big\}
\subset \bigcup_{q = j_\ell+2}^{j_{\ell +1}-2}\Big\{\Big| \bar M_n({j}_\ell;q)\Big|\geq \frac{ m\varepsilon_n\Delta_{\beta,\min} }{18\sqrt{m}}\Big\}
\end{equation*}
Let $\varphi'_n :=  \frac{ \sqrt{m}\varepsilon_n \Delta_{\beta,\min}}{18}$. Hence, by \eqref{ctrlmart} in Lemma~\ref{lem:control-martingale} we  get
\begin{eqnarray*}
\P[A'_{n,\ell,3}]&\leq& 2 \sum_{q=j_{\ell}+2}^{j_{\ell +1}-2} \exp\Bigg(-\frac{n{\varphi'}_n^2}{2 \E\Big[\bar N_n\Big(\big(\frac{j_{\ell}-1}{m},\frac{q}{m}\big]\Big)\Big] +\frac{2}{3}\varphi'_n }\Bigg)\\
&\leq& 2 (j_{\ell+1} - j_{\ell}-3)\exp\Bigg(-\frac{n{\varphi'}_n^2}{2 \E\Big[\bar N_n\Big(\big(\frac{j_{\ell}-1}{m},\frac{j_{\ell +1}-2}{m}\big]\Big)\Big]+\frac{2}{3}\varphi'_n }\Bigg)\\
&\leq& 2 \exp\Bigg(-\frac{n{\varphi'}_n^2}{2 \E\Big[\bar N_n\Big(\big(\frac{j_{\ell}-1}{m},\frac{j_{\ell +1}-2}{m}\big]\Big)\Big]+\frac{2}{3}\varphi'_n } + \log m\Bigg).
\end{eqnarray*}
By ~\eqref{ass:consistency-1} in  Assumption~\ref{ass:consistency},  we get  $\P[A'_{n,\ell,3}]$ tends to zero as $n$ tends to infinity. Let us now address $\P[A'_{n,\ell,2}]$. Using $(\ref{kkt}$) in Lemma~\ref{lem:KKT} with $j = j_\ell$ and with $j=\lceil\frac{j_\ell + j_{\ell-1}}{2}\rceil$, and using the triangle inequality, it follows that
\begin{equation*}
\Big| \sum_{q=\lceil\frac{j_\ell + j_{\ell-1}}{2}\rceil}^{j_\ell-1}\bN_q  - \sum_{q=\lceil\frac{j_\ell + j_{\ell-1}}{2}\rceil}^{j_\ell -1}\hat{\beta}_{q,m}\Big| \leq \hat{w}_{\lceil\frac{j_\ell + j_{\ell-1}}{2}\rceil, j_\ell}. 
\end{equation*}
On the event $C_n \cap \{\hat{j}_\ell > j_\ell\},$  we get
\begin{equation*}
 \Big| \frac{{j}_{\ell} - j_{\ell-1} }{2}(\beta_{0,j_{\ell}-1,m} - \hat{\beta}_{\hat{j}_{\ell}-1,m}) + \sqrt{m}\bar M_n(\lceil\frac{j_\ell + j_{\ell-1}}{2}\rceil;{j}_\ell-1) \Big| \leq  \hat{w}_{\lceil\frac{j_\ell + j_{\ell-1}}{2}\rceil, j_\ell}.
\end{equation*}
This implies 
\begin{equation*}
| \frac{{j}_{\ell} - j_{\ell-1} }{2}| |\hat{\beta}_{\hat{j}_{\ell}-1,m} - \beta_{0,j_{\ell}-1,m}| \leq \hat{w}_{\lceil\frac{j_\ell + j_{\ell-1}}{2}\rceil, j_\ell} +\Big|\sqrt{m}\bar M_n(\lceil\frac{j_\ell + j_{\ell-1}}{2}\rceil;{j}_\ell -1) \Big|.
\end{equation*}
Therefore, we may upper bound $\P[A'_{n,\ell,2}]$ as follows
\begin{equation*}
\begin{split}
  &\P[A'_{n,\ell,2}]\\
&\qquad = \P\Big[\Big\{ | \hat{\beta}_{\hat{j}_{\ell}-1,m}-\beta_{0,j_{\ell}-1,m}| \geq \frac{| \beta_{0,j_{\ell+1}-1,m}- \beta_{0,j_\ell-1,m}|}{3}\Big\}\cap C_n \cap\{\hat{j}_\ell > j_\ell\}\Big] \\
  &\qquad = \P\Big[\Big\{| \frac{{j}_{\ell} - j_{\ell-1} }{2}| | \hat{\beta}_{\hat{j}_{\ell}-1,m}-\beta_{0,j_{\ell}-1,m}|\\
&\hspace{5cm} \geq | \frac{{j}_{\ell} - j_{\ell-1} }{2}| \frac{| \beta_{0,j_{\ell+1}-1,m}- \beta_{0,j_\ell-1,m}|}{3} \Big\}\cap C_n\Big]\\
   &\qquad \leq \P\Big[\Big \{ \hat{w}_{\lceil\frac{j_\ell + j_{\ell-1}}{2}\rceil, j_\ell} +\Big|\sqrt{m}\bar M_n(\lceil\frac{j_\ell + j_{\ell-1}}{2}\rceil;{j}_\ell -1) \Big|\\
&\hspace{5cm} \geq | \frac{{j}_{\ell} - j_{\ell-1} }{2}|  \frac{| \beta_{0,j_{\ell+1}-1,m}- \beta_{0,j_\ell-1,m}|}{3} \Big\}\cap C_n\Big]\\
   &\qquad \leq \P\Big[\hat{w}_{\lceil\frac{j_\ell + j_{\ell-1}}{2}\rceil, j_\ell} \geq ({j}_{\ell} - j_{\ell-1}) \frac{| \beta_{0,j_{\ell+1}-1,m}- \beta_{0,j_\ell-1,m}|}{6}\Big] \\
   & \qquad \qquad+ \P\Big[\Big|\sqrt{m}\bar M_n(\lceil\frac{j_\ell + j_{\ell-1}}{2}\rceil;{j}_\ell -1) \Big|\geq \frac{| \beta_{0,j_{\ell+1}-1,m}- \beta_{0,j_\ell-1,m}|}{6}\Big]\\
&\qquad \leq \P\Big[\hat{w}_{\lceil\frac{j_\ell + j_{\ell-1}}{2}\rceil, j_\ell} \geq \frac{\Delta_{j,\min}\Delta_{\beta,\min}}{12}\Big]  +\P\Big[\Big|\bar M_n(\lceil\frac{j_\ell + j_{\ell-1}}{2}\rceil;{j}_\ell -1) \Big|\geq \frac{\Delta_{j,\min} \Delta_{\beta,\min}}{12\sqrt{m}}\Big]\\
&\qquad := \alpha'_{n,\ell,2}{^{(1)}} + \alpha'_{n,\ell,2}{^{(2)}} .
\end{split}
\end{equation*}
We observe that 
\begin{equation*}
\alpha'_{n,\ell,2}{^{(1)}} \leq \P\Big[\hat{w}^2_{\lceil\frac{j_\ell + j_{\ell-1}}{2}\rceil} \geq \frac{\Delta_{j,\min}^2\Delta_{\beta,\min}^2}{24^2}\Big] \leq \P\Big[\hat{w}^2_{ j_{\ell-1}} \geq \frac{\Delta_{j,\min}^2\Delta_{\beta,\min}^2}{24^2}\Big].
\end{equation*}
By  \eqref{ass:consistency-2} in Assumption~\ref{ass:consistency}, \eqref{ctrlweigh} in  Lemma~\ref{lem:control-martingale} with $\xi = \frac{n \Delta_{j,\min}^2\Delta_{\beta,\min}^2}{24^2m\log m} + \E\big[\bar N_n\big((\frac{j_{\ell-1}-1}{m},1]\big)\big],$ it follows that 
\begin{eqnarray*}
\alpha'_{n,\ell,2}{^{(1)}} &\leq&   2\exp\Bigg(-\frac{n\xi^2}{2\E\Big[\bar N_n\Big(\big(\frac{j_{\ell-1}-1}{m},1\big]\Big)\Big] +\frac{2}{3}\xi }\Bigg)\rightarrow 0,
\end{eqnarray*}
as $n \rightarrow \infty.$
By \eqref{ctrlmart} in  Lemma~\ref{lem:control-martingale} with $ z=\frac{\Delta_{j,\min} \Delta_{\beta,\min}}{12\sqrt{m}}$  and \eqref{ass:consistency-2} in Assumption~\ref{ass:consistency}, we obtain
\begin{eqnarray*}
\alpha'_{n,\ell,2}{^{(2)}} &\leq&  2\exp\Bigg(-\frac{nz^2}{2\E\Big[\bar N_n\Big(\big(\frac{\lceil\frac{j_\ell + j_{\ell-1}}{2}\rceil-1}{m},\frac{j_\ell-1}{m}\big]\Big)\Big] +\frac{2}{3}z }\Bigg)\rightarrow 0,
\end{eqnarray*}
as $n \rightarrow \infty.$
Therefore, we conclude that $\P[A'_{n,\ell,2}]\rightarrow 0, \, as\, n\rightarrow \infty.$

 \subsection{Step II.2. Prove: $ \P[A_{n,\ell} \cap C^\complement_n]\rightarrow 0, \, as \,  n\rightarrow \infty.$ }
\label{subsection2:CASE-II}

 As in Case~I from~\ref{CASE-I}, we split $\P[A_{n,\ell} \cap
 C_n^\complement]$ into
 \begin{equation*}
\P[A_{n,\ell} \cap C_n^\complement] = \P[A_{n,\ell} \cap D_n^{(l)}] + \P[A_{n,\ell} \cap D_n^{(m)}]+ \P[A_{n,\ell} \cap D_n^{(r)}].
\end{equation*}
\noindent Let us first focus on $\P[A_{n,\ell} \cap D_n^{(m)}]$. Note that
 \begin{equation}
\label{eq:37a}
\begin{split}
 \P[A_{n,\ell} \cap D_n^{(m)}\cap \{\hat{j}_\ell > j_\ell\}] &\leq \P[A_{n,\ell}\cap B_{\ell+1,\ell} \cap D_n^{(m)}] \\
& \qquad + \sum_{l=\ell+1}^{L_0-2} \P[C_{l,l}\cap B_{l+1,l} \cap D_n^{(m)}].
\end{split}
 \end{equation}
Let us now prove that the first term in the right hand side of (\ref{eq:37a}) goes to zero as $n$ tends to infinity.
\noindent Using $(\ref{kkt})$  in Lemma~\ref{lem:KKT} with $j = \lceil \frac{j_{\ell+1}+j_\ell}{2}\rceil$ and $j = j_{\ell+1}$, on the first hand and $(\ref{kkt})$ in Lemma~\ref{lem:KKT}  with $j= j_{\ell+1} +1$ and  with $j = j_{\ell+2}$ on the other hand, we obtain, respectively
\begin{equation}
\label{eq:38a}
\begin{split}
&|\frac{{j}_{\ell+1} - j_{\ell}}{2}||\hat{\beta}_{\hat{j}_{\ell+1}-1,m} -\beta_{0,j_{\ell+1}-1,m}| \leq \hat{w}_{\lceil \frac{j_{\ell+1}+j_\ell}{2}\rceil,j_{\ell+1}}\\
&\hspace{7cm}+ |\sqrt{m}\bar M_n(\lceil \frac{j_{\ell+1}+j_\ell}{2}\rceil;{j}_{\ell+1} -1)|,
\end{split}
\end{equation}
and
\begin{equation}
\label{eq:39a}
\begin{split}
&| {j}_{\ell+2} - j_{\ell+1} - 2||\hat{\beta}_{\hat{j}_{\ell+1}-1,m} -\beta_{0,j_{\ell+2}-1,m}| \leq \hat{w}_{{j}_{\ell+1}+1,j_{\ell+2}} \\
&\hspace{7.5cm}+| \sqrt{m}\bar M_n({j}_{\ell+1}+1;{j}_{\ell+2} -1)|.
\end{split}
\end{equation}
In addition, we have
\begin{flalign*}
&|\beta_{0,j_{\ell+2}-1,m} - \beta_{0,j_{\ell+1}-1,m}| \\
&\qquad \qquad= |( \hat{\beta}_{\hat{j}_{\ell+1}-1,m} - \beta_{0,j_{\ell+1}-1,m}) - (\hat{\beta}_{\hat{j}_{\ell+1}-1,m} - \beta_{0,j_{\ell+2}-1,m})|\\
                        & \qquad\qquad\leq |\hat{\beta}_{\hat{j}_{\ell+1}-1,m} - \beta_{0,j_{\ell+1}-1,m}| + |\hat{\beta}_{\hat{j}_{\ell+1}-1,m} - \beta_{0,j_{\ell+2}-1,m}|\\
                        &\qquad\qquad\leq \frac{ \hat{w}_{\lceil \frac{j_{\ell+1}+j_\ell}{2}\rceil,j_{\ell+1}}}{|\frac{{j}_{\ell+1} - j_{\ell}}{2}|} + \frac{\sqrt{m}\bar M_n(\lceil \frac{j_{\ell+1}+j_\ell}{2}\rceil;{j}_{\ell+1} -1)|}{|\frac{{j}_{\ell+1} - j_{\ell}}{2}|} \\
&\hspace{2.5cm}+ \frac{ \hat{w}_{{j}_{\ell+1}+1,j_{\ell+2}} }{| {j}_{\ell+2} - j_{\ell+1} -2|} + \frac{| \sqrt{m}\bar M_n({j}_{\ell+1}+1;{j}_{\ell+2} -1)|}{|{j}_{\ell+2} - j_{\ell+1}-2|}\\
&\qquad \qquad\leq 2\frac{\hat{w}_{\lceil \frac{j_{\ell+1}+j_\ell}{2}\rceil,j_{\ell+1}}}{\Delta_{j,\min}} +  2\frac{\sqrt{m}\bar M_n(\lceil \frac{j_{\ell+1}+j_\ell}{2}\rceil;{j}_{\ell+1} -1)|}{\Delta_{j,\min}}\\
&\hspace{2.5cm}+ \frac{\hat{w}_{{j}_{\ell+1}+1,j_{\ell+2}}}{|\Delta_{j,\min}-2|} + \frac{| \sqrt{m}\bar M_n({j}_{\ell+1}+1;{j}_{\ell+2} -1)|}{|\Delta_{j,\min}-2  |}.
\end{flalign*}
Define the event $E'_{n,\ell}$ by
\begin{flalign*}
\nonumber
   &E'_{n,\ell} = \Bigg\{|\beta_{0,j_{\ell+2}-1,m} - \beta_{0,j_{\ell+1}-1,m}| \leq \frac{2\hat{w}_{\lceil \frac{j_{\ell+1}+j_\ell}{2}\rceil,j_{\ell+1}}}{\Delta_{j,\min}} + \frac{6\hat{w}_{{j}_{\ell+1}+1,j_{\ell+2}}}{\Delta_{j,\min}} \\ 
 &\hspace{6.5cm} +  \frac{2\sqrt{m}\bar M_n(\lceil \frac{j_{\ell+1}+j_\ell}{2}\rceil;{j}_{\ell+1} -1)|}{\Delta_{j,\min}}\\
&\hspace{6.5cm} + \frac{|6 \sqrt{m}\bar M_n({j}_{\ell+1}+1;{j}_{\ell+2} -1)|}{\Delta_{j,\min}}\Bigg\}.
\end{flalign*}
$E'_{n,\ell}$ occurs with probability  one. Therefore, we obtain
\begin{equation*}
\begin{split}
  &\P[A_{n,\ell}\cap B_{\ell+1,\ell} \cap D_n^{(m)}]\\
&\qquad \leq   \P[A_{n,\ell}\cap B_{\ell+1,\ell} \cap D_n^{(m)}\cap E'_{n,\ell}]\\
&\qquad \leq \P\Big[\hat{w}_{\lceil \frac{j_{\ell+1}+j_\ell}{2}\rceil,j_{\ell+1}} \geq \frac{\Delta_{j,\min}| \beta_{0,j_{\ell+2}-1,m} - \beta_{0,j_{\ell+1}-1,m}|}{8}\Big] \\
  & \qquad \quad +\P\Big[\hat{w}_{{j}_{\ell+1}+1,j_{\ell+2}} \geq \frac{|\beta_{0,j_{\ell+2}-1,m} - \beta_{0,j_{\ell+1}-1,m}|}{24}\Big]\\
  & \qquad \quad + \P\Big[|\bar M_n(\lceil \frac{j_{\ell+1}+j_\ell}{2}\rceil;{j}_{\ell+1} -1)| \geq \Delta_{j,\min}\frac{| \beta_{0,j_{\ell+2}-1,m} - \beta_{0,j_{\ell+1}-1,m}|}{8\sqrt{m}} \Big]\\
  & \qquad \quad + \P\Big[|\bar M_n({j}_{\ell+1}+1;{j}_{\ell+2} -1)| \geq  \Delta_{j,\min}\frac{| \beta_{0,j_{\ell+2}-1,m} - \beta_{0,j_{\ell+1}-1,m}|}{24 \sqrt{m}}\Big]\\
&\qquad := \theta'_{n,\ell,1} + \theta'_{n,\ell,2} + \theta'_{n,\ell,3} + \theta'_{n,\ell,4}.
\end{split}
\end{equation*}
By \eqref{ctrlmart}-\eqref{ctrlweigh} in
Lemma~\ref{lem:control-martingale}, and \eqref{ass:consistency-1}-\eqref{ass:consistency-2}  in Assumption~\ref{ass:consistency}, we show
that for $s=1, \ldots, 4, \theta'_{n,\ell,s } \rightarrow 0,\,
\P[A_{n,\ell}\cap B_{\ell+1,\ell} \cap D_n^{(m)}] \rightarrow 0,$ as  $n\rightarrow \infty.$
Recall that in  Case~I from Section~\ref{CASE-I}, we
proved  $\P[A_{n,\ell}\cap D_n^{(l)}] \rightarrow 0,$  as  $n\rightarrow \infty$ and in a
similar way $\P[A_{n,\ell}\cap D_n^{(r)}]\rightarrow 0,$ as  $n\rightarrow \infty.$ This
concludes the proof of Theorem~\ref{thm3}.  $\hfill \square$

\end{appendices}

\bibliographystyle{plain}

{\small

\bibliography{biblio}
}

\end{document}